\documentclass{article}

\usepackage{authblk}
\usepackage{amsrefs}
\usepackage{amsmath}
\usepackage{graphics}
\usepackage{epsfig}
\usepackage{mathtools}
\usepackage{nccmath}
\usepackage{amssymb}
\usepackage{xfrac}
\usepackage{amscd}
\usepackage[all]{xy}
\usepackage{latexsym}
\usepackage{graphicx}
\usepackage{comment}
\usepackage{multirow}
\usepackage{geometry}
\usepackage{color}
\usepackage[usenames,dvipsnames]{pstricks}
\usepackage{epsfig}
\usepackage{pst-grad} 

\usepackage{pst-plot} 
\usepackage{subfig}
\usepackage{float}
\usepackage[utf8]{inputenc}
\usepackage{verbatim}

\usepackage{hyperref}

\newtheorem{thm}{Theorem}[section]

\newtheorem{cor}[thm]{Corollary}
\newtheorem{lem}[thm]{Lemma}
\newtheorem{prop}[thm]{Proposition}

\newtheorem{Def}[thm]{Definition}
\newtheorem{rem}[thm]{Remark}

\newtheorem{ex}{Example}[section]

\title{An Introduction to \(p\)-adic Hodge Theory}
\author[1]{Olivia Dumitrescu}
\author[2]{Kexuan Yang} 

\affil[1,2]{University of North Carolina at Chapel Hill Chapel Hill, NC 27599}

\begin{document}

\maketitle
\begin{abstract}
This survey gives an introductory review of the classical Hodge theory, arithmetic geometry and algebraic number theory that inspires \(p\)-adic Hodge theory. Basic results in this new area developed by Tate, Faltings and Scholze are presented with corresponding mathematical constructions. Some background in algebraic geometry and algebraic topology is required.\footnote{
 Mathematics Subject Classification [2020] Primary: 11E95, 11-01. Secondary:  14F85, 11F80\\
 Email: [1] dolivia@unc.edu and Email [2] ykx13579@gmail.com}

\end{abstract}



\tableofcontents

\section{Introduction}
\subsection{General Introduction}

The authors wrote this note in their journey of exploring the transition from algebraic geometry to arithmetic geometry. This note represents a survey of many papers, textbooks and lecture notes that were used by the authors for building the basic foundation of their knowledge. In algebraic topology, cohomology groups form a category of significant topological invariants. In the case of smooth, complex manifolds, de Rham cohomology theory stands out as the cohomology theory of particular analytic interest, and is natural to attempt to find out a "representative" of each element of a cohomology group, a quotient group of differential forms, that bears some analytical significance. Fortunately, in complex Hodge theory, the space of harmonic functions can serve this purpose through the decomposition of the cohomology group known as Hodge Decomposition. By identifying the cohomology groups with the spaces of harmonic functions, a series of results can be eluded. Especially, the degeneracy of Hodge-de Rham spectral sequences for K\"ahlerian manifolds yield a beautiful result known as Hodge symmetry.

The development of algebraic geometry gives a purely algebraic definition of de Rham cohomology that can be generalised to non-smooth geometric results. Under the context of arithmetic geometry there arises the question of characterisation of de Rham cohomology on \(p\)-adic varieties. An analogue of Hodge Decomposition was proved on \(p\)-adic varieties by Faltings. On the analytical side, one can generate the theory to a kind of \(p\)-adic geometrical object known as rigid analytic space. On a rigid analytic space, Faltings' result can be extended without assuming K\"ahlerian conditions. That is established by Peter Scholze. However, the Hodge Symmetry part is more trivial and does not hold in general, although David Hansen and Shizhang Li established that for lower cohomology groups with some qualifications analogous to K\"ahlerian conditions.

The works of Scholze involve techniques as perfectoid spaces to deal with spaces with complicated topologies. These techniques are generalized by the same author and Dustin Clausen to a program called condensed mathematics aiming at uniting topology, complex analysis and functional analysis. It is breifly covered but not explained in detail in this note.

The study of Hodge theory also gives inspiration to the study of mirror symmetry. Indeed, Hodge numbers are significant indicators of mirror symmetry. In particular, if the complex manifold \(V\) is the mirror dual to the symplectic manifold \(V^{\circ}\), then the Hodge number satisfies the equation \(h^{p,q}(V)=h^{3-p,q}(V^{\circ})\). Maxim Kontsevich furthered on this direction (see \cite{Kont}).

At the beginning of each chapter we give the major references that were used; we leave the readers to consult the original material. There are many resources the reader can use to learn about this theory, we apologize if some were unintentionally omitted by the authors. Since there are already many resources available, this work is meant to overview related topics and give students a map for the foundations of these theories. The authors intend to continue their work incorporating more recent related topics in the areas of $p$ adic geometry and mathematical physics.

\subsection{Preliminaries}
In this note we presume that the readers have basic knowledge on algebraic geometry (for instance, the definition of sheaves, schemes, algebraic varieties) and algebraic topology (for instance, de Rham cohomology and \v{C}ech cohomology on real and complex manifolds, Poincare duality, K\"unneth formula). One can consult \cite{Bo1} and \cite{Has} for explanations. For de Rham cohomology theory one can consult \cite{Hat}. Basic knowledge of complex geometry (for instance, the definition of complex manifolds, K\"ahleran condition) is also assumed and \cite{bri} is a good source on it. One is also expected to understand basic homological algebra and category theory (for instance, concepts like derived functors, Abelian categories, direct and inverse limits). The best resources to consult are  \cite{Lee} and \cite{Vos2}. No knowledge on analysis is assumed besides the very basics.

\subsection*{Acknowledgements} The authors would like to express their gratitude to the Max Planck Institute of Mathematics, MPI Leipzig; Institute of Mathematics of the Romanian Academy IMAR Bucharest, and RIMS Institute, Kyoto for their long stays in the summers of 2026, 2025 and 2024 respectively. In particular, they are grateful to Takuro Mochizuki and Masa-Hiko Saito for the ample daily discussions, seminars and conferences organized during the program "Integrable Systems and Mathematical Physics" 2024, when this work was initiated. 

The authors would like to express their gratitude to Max Planck Institute of Mathematics Bonn and IHES, Bures sur Yvette and for their generosity, hospitality and stimulating environment during the long stay of the first author's in 2023. The first author is also affiliated to the Institute of Mathematics of the Romanian Academy, IMAR, Bucharest. We also thank Bernd Sturmfels for our discussions in MPIM, Leipzig, where the authors finalize the note.

The author's research was partially supported by the NSF-FRG DMS 2152130 grant and by the UNC JFDA Award 2022. 

\section{Classical Hodge Theory}\label{Hodge theory}
The major reference of this chapter is \cite{GrHa}. Throughout this chapter we will assume that \(M\) is a connected, compact, complex manifold of complex dimension \(n\). 

Define \(A^{p,q}(M)\) as the module of \(p,q\)-forms on a complex manifold \(M\). Since we can  choose a hermitian metric \(ds^{2}\) on \(M\), it is therefore possible to extend the given metric to an inner product of any two \(p\), \(q\) forms and to associate a volume form on \(M\), denoted by $\Phi$. Thus, for \(\psi,\eta \in A^{p,q}(M)\), we define the global inner product

\begin{equation}
(\psi,\eta)=\int_{M}(\psi(z),\eta(z))\Phi(z)\nonumber
\end{equation}
that makes \(A^{p,q}(M)\) a pre-Hilbert space. For a heuristic argument, assume that it is complete with respect to the given norm and is thus a Hilbert space and also assume that \({\overline{\partial}}\) is bounded as a linear operator.

First define the Hodge star \(*:A^{p,q}(M)\to A^{n-p,n-q}(M)\). Let \(\{\phi_{1},\phi_{2}...\}\) be the unitary coframe under the given Hermitian metric and let \(\phi_{I}=\wedge_{i\in I}\phi_{i}\) for all \(I\subseteq\{1,2,...,n\}\). Locally, \(\eta\in A^{p,q}(M)\) can be expressed as \(\eta=\sum_{I,J}\eta_{I\bar{J}}\phi_{I}\wedge \overline{\phi_{J}}\), where \(\eta_{I\bar{J}}\in\mathbb{C}\). Let $\overline{\eta_{I\overline{J}}}$ denote the complex conjugate of $\eta_{I\overline{J}}$, and define

\begin{equation}
*\eta=2^{p+p-n}\Sigma_{I,J} \epsilon_{IJ} \cdot \overline{\eta_{I\overline{J}}} \cdot \phi_{I^{0}}\wedge \overline{\phi_{J^{0}}}\nonumber
\end{equation}
where $I^{0}$ stands for the complement of the set $I$,  \(I^{0}:=\left\{1,..,n\right\}-I\), while $\overline{\phi_{J^{0}}}$ is the complex conjugate of \(\phi_{J^{0}}\), and \(\epsilon_{IJ}\) is the sign of the permutation 
\begin{equation*}
\left(1,...,n,1^{'},...,n^{'}\right)\to\left(i_{1},...,i_{p},j_{1},...,j_{q},i^{0}_{1},...,i^{0}_{n-p},j^{0}_{1},...,j^{0}_{n-q}\right).\nonumber
\end{equation*}

It can then be verified by direct calculation that
\begin{equation}
**\eta=(-1)^{p+q}\eta\nonumber
\end{equation}
and that for any two forms $\psi$ and $\eta$ in $A^{p,q}(M)$ and the volume form $\Phi$ then
\begin{equation}
(\psi(z),\eta(z))\Phi(z)=\psi(z)\wedge *\eta(z).\nonumber
\end{equation}

Define the adjoint operator \(\overline{\partial}^{*}\) as
\begin{equation}
\overline{\partial}^{*}=-*\overline{\partial}*.\nonumber
\end{equation}
It can be calculated that
\begin{equation}
(\overline{\partial}\psi,\rho)=(\psi,\overline{\partial}\rho)\nonumber
\end{equation}
for all \(\rho\in A^{p,q-1}\left(M\right)\).

For each Dolbeault cohomology group \(H^{p,q}_{\bar{\partial}}(M)\) there arises the problem to find a representative of the group with the minimal norm given the Hermitian metric. In fact, this problem is solved by the introduction of \(\overline{\partial}^{*}\) operator.

\begin{lem}\label{minimal norm} Let  $M$ be a
connected, compact, complex manifold. A \(\overline{\partial}\)-closed form \(\psi\in Z^{p,q}_{\overline{\partial}}(M)\) is of minimal norm in \(\psi+\overline{\partial}A^{p,q-1}(M)\) if and only if \(\overline{\partial}^{*}\psi=0\), for any non-negative integers $p$ and $q$.
\end{lem}

Thus, the Dolbeault cohomology group  
$$H^{p,q}_{\overline{\partial}}(M)=\frac{Z^{p,q}_{\overline{\partial}}(M)}{\overline{\partial}A^{p,q-1}(M)}$$
is represented by the first-order equations
\begin{equation}
\centering
\left\{\begin{split}
\overline{\partial}\psi=0 \\
\overline{\partial}^{*}\psi=0. \\ 
\end{split}\right.\nonumber
\end{equation}
The system of equations is equivalent to the equation
\begin{equation}\label{delta}
\Delta_{\overline{}{\partial}}\psi =\left(\overline{\partial}\overline{\partial}^{*}+\overline{\partial}^{*}\overline{\partial}\right)\psi=0
\end{equation}
where \(\Delta_{\overline{\partial}}:A^{p,q}(M)\to A^{p,q}(M)\) is called the \(\overline{\partial}\)-Laplacian, or simply the Laplacian (written as \(\Delta\)). Differential forms satisfying the Laplace equation
\begin{equation}
\Delta\psi =0\nonumber
\end{equation}
are called harmonic forms. We will denote the space of harmonic forms of type \(\left(p,q\right)\), as \(\mathcal{H}^{p,q}(M)\nonumber\) and we refer to it as the \textbf{harmonic space}. Equation \ref{isomorphism} justifies that in fact we have the following relation
\begin{equation}
\mathcal{H}^{p,q}(M)\cong H^{p,q}_{\overline{\partial}}(M)\nonumber
\end{equation}
since every Dolbeault cohomology group is associated with a harmonic form.

We have the following famous result, known as Hodge Theorem. This result assumes knowledge of functional analysis and harmonic analysis:

\begin{thm}\label{dimensionality theorem}
Let $M$ be a
connected, compact, complex manifold.
The following statements hold:
\begin{enumerate}
\item $\dim \mathcal{H}^{p,q}(M)<\infty$;

\item The orthogonal projection \(\mathcal{H}:A^{p,q}\left(M\right)\to \mathcal{H}^{p,q}\left(M\right)\) is well-defined. And there exists a unique operator (the Green's operator)
\begin{equation}
G:A^{p,q}\left(M\right)\to A^{p,q}\left(M\right),\nonumber
\end{equation}
with \(G\left(\mathcal{H}^{p,q}\left(M\right)\right)=0\), \(\overline{\partial}G=G\overline{\partial}\), \(\overline{\partial}^{*}G=G\overline{\partial}^{*}\). Moreover, on \(A^{p,q}\left(M\right)\) the following holds
\begin{equation}\label{I=H}
I=\mathcal{H}+\Delta G.
\end{equation}

\end{enumerate}
for any non negative integers $p$ and $q$.
\end{thm}

Especially, for all \(\psi\in A^{p,q}(M)\) there is
\begin{equation}
\psi=\mathcal{H}(\psi)+\overline{\partial}\left(\overline{\partial}^{*}G\psi\right)+\overline{\partial}^{*}\left(\overline{\partial}G\psi\right),\nonumber
\end{equation}
which is called the Hodge decomposition on forms.

Consider the equation
\begin{equation}
\Delta\psi=\gamma\nonumber
\end{equation}
for some \(\gamma\in A^{p,q}(M)\)

If \(\mathcal{H}(\gamma)=0\), by Equation \eqref{I=H}
 this equation has a solution
\begin{equation}\label{Image of G}
\psi=G\left(\gamma\right).
\end{equation}
Conversely, if the equation has a solution \(\psi\), then there is \(\mathcal{H}\left(\gamma\right)=\Delta\left(\psi-G\left(\gamma\right)\right)\) by Equation \eqref{Image of G}. It further implies

$$\mathcal{H}\left(\mathcal{H}\left(\gamma\right)\right)=\mathcal{H}\left(\Delta\left(\phi\right)\right)=\Delta\left(\mathcal{H}\left(\phi\right)\right)=0 \text{ if } \phi=\psi-G\left(\gamma\right).$$ Since \(\mathcal{H}\) is a projection, we have  
$$\mathcal{H}\left(\gamma\right)=\mathcal{H}\left(\mathcal{H}\left(\gamma\right)\right)=0.$$ If \(\psi_{1}\) and \(\psi_{2}\) are two different solutions, then \(\Delta\left(\psi_{1}-\psi_{2}\right)=0\) and \(\psi_{1}-\psi_{2}\in\mathcal{H}^{p,q}(M)\), so \(\psi_{1}-\psi_{2}=0\).

\begin{cor}\label{unique solution} In the hypothesis of Theorem \ref{dimensionality theorem}, then 
equation \eqref{Image of G} has a unique solution
\begin{equation}
\psi=G(\gamma)\nonumber
\end{equation}
if and only if \(\mathcal{H}(\gamma)=0\).
\end{cor}

Indeed, for a \(\overline{\partial}\)-closed form \(\psi\in Z_{\overline{\partial}}^{p,q}(M)\)
, there is a Hodge decomposition
\begin{equation}
\psi=\mathcal{H}(\psi)+\overline{\partial}(\overline{\partial}^{*}G
\psi),\nonumber
\end{equation}
since \(\overline{\partial}G\psi=G\overline{\partial}\psi=0\). Therefore, every \(\overline{\partial}\)-closed form corresponds to a harmonic form in the same Dolbeault cohomology class. This correspondence, combining with the trivial inverse map, establishes an isomorphism 

\begin{equation}\label{isomorphism}
\mathcal{H}^{p,q}(M)\cong H_{\bar{\partial}}^{p,q}(M)
\end{equation}
between the harmonic space and Dolbeault cohomology groups. Combining with the Dolbeault isomorphism, there is \(\mathcal{H}^{p,q}(M)\cong H^{q}(M,\Omega^{p})\). By the first statement of Theorem \ref{dimensionality theorem}  we obtain the following consequence
\begin{cor}\label{finite dimensionality}
For $M$ a
connected, compact, complex manifold, then
$$\dim H^{q} (M, \Omega^{p})<\infty,$$
for non-negative integers $p$ and $q$.
\end{cor}

Recall the definition of the Laplacian \(\Delta\), since \(\overline{\partial}^{*}=-*\overline{\partial}*\) there is \(*\Delta=\Delta *\). Therefore \(*\) maps \(\mathcal{H}^{p,q}(M)\) to \(\mathcal{H}^{n-p,n-q}(M)\). There is thus an isomorphism between the two spaces \(\mathcal{H}^{p,q}(M)\) and \(\mathcal{H}^{n-p,n-q}(M)\).

In particular,
\(
\mathcal{H}^{n,n}(M)\cong \mathbb{C}\cdot\Phi
\) where \(\Phi\) is the volume form of the metric. By Equation \eqref{isomorphism}, this isomorphism is independent of the metric and
\begin{equation}\label{complementary isomorphism}
H^p(M,\Omega^{p})\cong H^{n-p}(M,\Omega^{n-p}).
\end{equation}

Consider the complex manifolds \(M\) and \(N\) with forms \(\psi_{M}\) and \(\psi_{N}\) on \(M\), \(N\) respectively. Define 
$$(\psi_{M}\otimes\psi_{N})(z,w)=\psi_{M}(z)\wedge\psi_{N}(w).$$ These  forms are known to be decomposable. By a technical lemma, the decomposable forms are \(L^{2}\)-dense in all the forms on \(M\times N\). It can also be proved technically that $$\Delta_{M\times N}=\Delta_{M}+\Delta_{N}$$ while \(\Delta_{M}\) and \(\Delta_{N}\) are pullbacks of correspondent Laplacian operator on \(M\times N\) through the projective map.

We have the following equality for decomposable forms

$$\Delta_{M\times N}(\psi_{M}\otimes\psi_{N})=(\Delta_{M}\psi_{M})\otimes\psi_{N}+\psi_{M}\otimes(\Delta_{N}\psi_{N}).$$ In particular, if 
$\{(\psi_{M})_{i}\}_{i\geq 1}$ is a complete set of eigenforms corresponding to eigenvalues $\{(\lambda_{i})\}_{i\geq 1}$ for \(\Delta_{M}\) and $\{(\psi_{N})_{j}\}_{j\geq 1}$ is a complete set of eigenforms corresponding to eigenvalues $\{(\mu_{j})\}_{j\geq 1}$  for \(\Delta_{N}\) then by the density of decomposable forms \(\{(\psi_{M})_{i}\otimes(\psi_{N})_{j}\}\) are a complete sets of eigenforms of \(\Delta_{M\times N}\) corresponding to eigenvalues \(\lambda_{i}+\mu_{j}\). Moreover, \(\Delta\) can be proved to be a real and non-negative operator, so \(\lambda_{i}+\mu_{j}=0\) is equivalent to that \(\lambda_{i}=\mu_{j}=0\). Therefore 
\begin{equation}
\mathcal{H}^{u,v}(M\times N)\cong\bigoplus_{p+r=u,q+r=s}\mathcal{H}^{p,q}(M)\otimes\mathcal{H}^{r,s}(N).\nonumber
\end{equation}
Consider the corresponding Dolbeault cohomology. The projection maps induce two injective maps \(H^{*}(M,\Omega^{p}_{M})\to H^{*}(M\times N,\Omega^{p}_{M\times N})\) and \(H^{*}(N,\Omega^{p}_{N})\to H^{*}(M\times N,\Omega^{p}_{M\times N})\).

 The K\"unneth formula can be derived by combining the previous results on the harmonic space.
\begin{thm}\label{cup product}
   Let $M$ and $N$ be two
connected, compact, complex manifolds. The cup product induces a group isomorphism
    \begin{equation}
        H^{*}(M,\Omega^{*}_{M})\otimes H^{*}(N,\Omega^{*}_{N})\cong H^{*}(M\times N,\Omega^{*}_{M\times N})
    \end{equation}

\end{thm}

Moreover, let \(M\) be, in addition, a K\"ahleran manifold. Define \(\Pi^{p,q}:A^{*}(M)\to A^{p,q}(M)\) and \(\Pi^{r}:A^{*}(M)\to A^{r}(M)\) to be the projections. We further define the \(\overline{\partial}\) - operators \(d^{*}\),\(\partial ^{*}\), \(\Delta_{d}\) and \(\Delta_{\partial}\) as  above via equations \ref{delta}. By local calculation, we get

\begin{prop}\label{delta operator}  On  a connected, compact, complex manifold $M$ the operator
    \(\Delta_{d}\)  satisfies:
    \begin{enumerate}
    \item $\Delta_{d}=2\Delta_{\overline{\partial}}=2\Delta_{\partial}$
    
    \item $[\Delta_{d},\Pi^{p,q}]=0.$

\end{enumerate}
\end{prop}
Let \(\psi\in\mathcal{H}^{r}(M)\), and \(p+q=r\), we have

$$\Delta_{\bar{\partial}}(\Pi^{p,q}\psi)=\frac{1}{2}\Delta_{d}(\Pi^{p,q}\psi)=\frac{1}{2}\Pi^{p,q}(\Delta_{d}\psi)=0.$$ This implies that its projection to the subspace of \(p,q\) forms is also harmonic. Threfore, 

$$\mathcal{H}^{r}(M)\cong\bigoplus_{p+q=r}\mathcal{H}^{p,q}(M).$$ Similarly, it can be proved that \(\Delta_{d}\) is real, therefore $$\Delta_{\overline{\partial}}(\overline{\psi})=\frac{1}{2}\Delta_{d}(\overline{\psi})=\frac{1}{2}\overline{\Delta_{d}(\psi)}=0.$$ That is, \(\overline{\psi}\) is also a harmonic form. So \(\mathcal{H}^{p,q}(M)\cong\overline{\mathcal{H}^{q,p}(M)}\). Denote by 

$$H^{p,q}(M)=\frac{Z^{p,q}_{d}(M)}{dA^{*}(M)\cap Z^{p,q}_{d}(M)}.$$ By part \(1\) of Proposition \refeq{delta operator} and a similar argument as that of Corollary \eqref{unique solution}, $$H^{p,q}(M)\cong \mathcal{H}^{p,q}(M).$$ There is similarly a Hodge theorem for \(\Delta_{d}\) concluding that \(H^{r}(M,\mathbb{C})\cong\mathcal{H}^{r}(M)\). 

Combining this with the previous results, there is the Hodge Decomposition:
\begin{thm}\label{decomposition thm}
    For a compact K\"ahler manifold \(M\), the complex cohomology satisfies
    \begin{equation}
        H^{r}(M,\mathbb{C})\cong\bigoplus_{p+q=r}H^{p,q}(M)\nonumber
    \end{equation}
    and
    \begin{equation}
        H^{p,q}(M)=\overline{H^{q,p}(M)}.\nonumber
    \end{equation}
\end{thm}

Finally define Hodge numbers
\begin{equation}\label{defHod}
    h^{p,q}(M)= \mathrm{dim}\ H^{q}(M,\Omega^{p}).
\end{equation}
We arrive to the following conclusion:
\begin{cor}\label{hodge numbers}
    If \(M\) is a complex compact manifold, then
    \begin{equation}
        \centering
\left\{\begin{split}
h^{p,q}(M)<\infty\\
h^{n,n}(M)=1\\
h^{p,q}(M)=h^{n-p,n-q}(M)\\
h^{u,v}(M\times N)=\sum_{p+r=u,q+s=v}h^{p,q}(M)h^{r,s}(N).
\end{split}\right.\nonumber
    \end{equation}
    If \(M\) is K\"ahler, then in addition
    \begin{equation}
        \centering
\left\{\begin{split}
h^{p,q}(M)=h^{q,p}(M)\\
b_{r}(M)=\sum_{p+q=r}h^{p,q}(M).
\end{split}\right.\nonumber
    \end{equation}
\end{cor}
\section{Hodge-de Rham  Spectral Sequence}\label{HdR spectral sequence}
The major reference of this section is \cite{Vos1} and \cite{Vos2}.

We recall that we denote a complex K\"ahler manifold by \(M\). Let \((A^{\bullet,\bullet},\partial,\overline{\partial})\) be the double complex of smooth differential forms on a complex manifold. We encode the conclusion \(b_{r}(M)=\sum_{p+q=r}h^{p,q}(M)\) in Corollary \eqref{hodge numbers}  with the language of the \textit{Hodge de-Rham Spectral Sequence}. Denote \(X\) to be any complex manifold. We will not apply the K\"aherian condition in this section before Corollary \eqref{DegenerationSequence2}. Note that \((A^{\bullet,\bullet},\partial,\overline{\partial})\) forms a double chain complex supported in positive degree with the columns \((A^{p,\bullet},\overline{\partial})\) being Doubeault complexes. For each double complex \((I^{\bullet,\bullet},D_{1},D_{2})\), we define the simple complex \textbf{associated to} it \((I^{\bullet},D)\) by setting $$I^{n}=\bigoplus_{p+q=n}I^{p,q}  \text{ and }  D=D_{1}+(-1)^{p}D_{2} \text{ on } I^{p,q}\subset I^{n}.$$ It can be easily verified that \((I^{\bullet},D)\) satisfies that \(D\circ D=0\). In our case, \((A^{\bullet,\bullet},\partial,\overline{\partial})\) is associated to de Rham complex \((\Omega^{\bullet},d)\).

For a double complex \((I^{\bullet,\bullet},D_{1},D_{2})\) supported in positive degree with associated simple complex \((I^{\bullet},D)\), set

$$F^{p}I^{k}=\bigoplus_{r+s=k,r\geq p}I^{r,s}.$$ 

Notice that \(F^{p}\) defines a \textit{decreasing filtration} on \(I^{k}\) with a family of sub-objects
\begin{equation}
    \cdots F^{p}I^{k}\xhookrightarrow{}F^{p-1}I^{k}\xhookrightarrow{}\cdots\xhookrightarrow{}F^{0}I^{k}=I^{k}.\nonumber
\end{equation}

Such a family of decreasing filtrations satisfies that each \(D\) sends \(F^{p}I^{k}\) to \(F^{p}I^{k+1}.\) Such a family defines a \textit{decreasing filtration on the complex} \(I^{\bullet}\). Recall that on a simple complex \(I^{\bullet}\) we use \(H^{k}(I^{\bullet})\) to denote the group \(\mathrm{ker}\ d/\mathrm{im}\ d\) at the \(I^{k}\). It naturally induces a filtration \(F^{p}\) on the cohomology of \(I^{\bullet}\) by setting $$F^{p}H^{i}(I^{\bullet})=\mathrm{Im}(H^{i}(F^{p}I^{\bullet})\to H^{i}(I^{\bullet})).$$  We further go back to the context of complex geometry and apply this to \((I^{\bullet,\bullet},\partial,\overline{\partial})\) and there will be a decreasing filtration on de Rham complex. We call such a filtration \textbf{Hodge filtration}. Note that on a complex manifold of finite dimension \(n\), for any \(k\), and for any \(l>n\) it gives \(F^{l}A^{k}=0\) in Hodge filtration.

Let \((I^{\bullet},d)\) be a simple complex and let \(F^{p}I^{\bullet}\) be a decreasing filtration on that complex, set $$Z^{p,q}_{r}=\{x\in F^{p}I^{p+q}\ |\ dx\in F^{p+r}I^{p+q+1}\}.$$ Note that \(Z^{p,q}_{r}\) naturally contains both sets \(Z^{p+1,q-1}_{r-1}\) and \(dZ^{p-r+1,q+r-2}_{r-1}\). Set $$B^{p,q}_{r}=Z^{p+1,q-1}_{r-1}+dZ^{p-r+1,q+r-2}_{r-1}\subseteq Z^{p,q}_{r}.$$ and \(E^{p,q}_{r}=Z^{p,q}_{r}/B^{p,q}_{r}\). Here \(d\) induces a differential \(d_{r}:E^{p,q}_{r}\to E^{p+r,q-r+1}_{r}\) that satisfies \(d_{r}^{2}=0\). We call the family of complexes \((E^{p,q}_{r},d_{r})\) the \textbf{spectral sequence} associated with the filtration \((I^{\bullet},F^{p}I^{\bullet})\). For a de Rham complex with the aforementioned filtration, we therefore have a \textbf{Hodge-de Rham spectral sequence}.

Set $$\mathrm{Gr}^{F}_{p}H^{i}(I^{\bullet})=F^{p}H^{i}(I^{\bullet})/F^{p+1}H^{i}(I^{\bullet}).$$ We present without a proof the following results.

\begin{thm}\label{chain complex}
    Let \((I^{\bullet},D)\) be a chain complex with a decreasing filtration and let \((E^{p,q}_{r},d_{r})\) be its associated spectral sequence. The complexes \((E^{p,q}_{r},d_{r})\) satisfies the followings:
\begin{enumerate}
    \item \(E^{p,q}_{0}=\mathrm{Gr}^{F}_{p}I^{p+q}=F^{p}I^{p+q}/F^{p+1}I^{p+q}\) and \(d_{0}\) is induced by \(d\).

    \item \(E^{p,q}_{r+1}\) can be identified with the cohomology of \((E_{r}^{p,q},d_{r})\).

    \item For \(p+q\) fixed and \(r\) sufficiently large, we have \(E^{p,q}_{r}\) stabilizes
    \begin{equation}\label{stabilizing}
        E^{p,q}_{r}=\mathrm{Gr}^{F}_{p}H^{p+q}(I^{\bullet}).
    \end{equation}
    \end{enumerate}
\end{thm}

We can therefore consider the spectral sequence as a machinery to get information about the \((p+q)^{th}\) cohomology group. Part \(1\) gives us an "initial page" and Part \(2\) gives us a method of "turning to a new page" by calculation a cohomology. Part \(3\) gives us the final result: after finitely many calculations, the \((p,q)^{th}\) entry on a page stabilizes and gives some information about \(H^{p+q}\). We give a concept measuring the place of stabilization.

\begin{Def}\label{degeneration}
    Let \((I^{\bullet},F)\) be a filtered complex with spectral sequence \((E^{p,q}_{r},d_{r})\), we say that \((E^{p,q}_{r},d_{r})\) \textbf{degenerates} at \(E_{r}\) if for any \(k\geq r\), \(d_{k}=0\).
\end{Def}

Note that if the spectral sequence degenerates at \(r\), we have 
$$E^{p,q}_{r}=E^{p,q}_{\infty}=\mathrm{Gr}^{F}_{p}H^{p+q}(I).$$

Consider de Rham complexes and Hodge-de Rham spectral sequences. Theorem \ref{chain complex} implies \(E_{0}^{p,q}=A^{p,q}\). By some routine calculations we get 
\begin{equation}\label{EHH1}
E_{1}^{p,q}(X)=H^{q}(X,\Omega^{p}_{X})=H^{p,q}_{\bar{\partial}}(X)
\end{equation}
for a complex manifold \(X\). Note that by Theorem \(4.1(3)\) there exists an \(r\) such that 
\begin{equation}\label{EHH2}
E^{p,q}_{r}(X)=\mathrm{Gr}^{F}_{p}H^{p+q}(A^{\bullet})=F^{p}H^{p+q}(X)/F^{p+1}H^{p+q}(X).
\end{equation}
Note that we compute \(E_{i+1}^{p,q}\) as the cohomology \(\mathrm{Ker}\ d_{i}/\mathrm{Im}\ d_{i}\) at \(E_{i}^{p,q}\) by Theorem \ref{chain complex}. Therefore, we obtain

$$\mathrm{dim}\ E^{p,q}_{i+1}\leq\mathrm{dim}\ E^{p,q}_{i}$$

with equality for every \((p,q)\) if and only if \(d_{i}=0\) for all \(1\leq i\leq r\). Equations \eqref{EHH1} and \eqref{EHH2} imply

$$H^{p,q}_{\bar{\partial}}(X)=F^{p}H^{p+q}(X)/F^{p+1}H^{p+q}(X)$$ if and only if \(d_{i}=0\) for any \(i\geq1\), that is, if and only if the Hodge de-Rham spectral sequence degenerates at \(E_{1}\). For each \(k\geq0\), by summing over \(p,q\geq0\) such that \(p+q=k\), we have that the degeneracy of the Hodge de-Rham spectral sequence at \(E_{1}\) implies that \(b_{k}=\sum_{p+q=k}h^{p,q}\). Conversely, by the previous argument, if the Hodge de-Rham sequence does not degenerate at \(E_{1}\), we will have 

$$F^{a}H^{k}(X)/F^{a+1}H^{k}(X)<H^{a,b}_{\bar{\partial}}(X)$$

for some \(a,b\geq0\) and \(k=a+b\). For any other \(p,q\geq0\) and \(p+q=k\), we further obtain

$$F^{p}H^{k}(X)/F^{p+1}H^{k}(X)\leq H^{p,q}_{\bar{\partial}}(X).$$ We obtain \(b_{k}<\sum_{p+q=k}h^{p,q}\). We have the following result.
\begin{thm}\label{DegenerationSequence}
    Let \(X\) be a complex manifold, its Hodge de Rham spectral sequence degenerates at \(E_{1}\) if and only if \(b_{k}=\sum_{p+q=k}h^{p,q}\) for any \(k\geq0\).
\end{thm}
\begin{rem} The part \(E_{1}^{p,q}(X)=H^{q}(X,\Omega^{p}_{X})\) of \ref{EHH1} can be generalized to \(E^{p,q}_{1}=H^{q}(I^{p,\bullet})\) for a double space that is associated to \((I^{\bullet,\bullet},D_{1},D_{2})\). We can still apply our previous argument to get that the associated spectral sequence degenerates at \(1\) if and only if \(H^{q}(I^{p,\bullet})=F^{p}H^{p+q}/F^{p+1}H^{p+q}\).
\end{rem}

Specifically, we have the corollary from Corollary \ref{hodge numbers}.
\begin{cor}\label{DegenerationSequence2}
    The Hodge de Rham sequence degenerates at \(E_{1}\) for any K\"ahler manifold \(M\).
\end{cor}
\begin{rem} For any complex, compact, connected manifold $M$, the degeneracy of \(E^{p,q}_{1}\) is very close to, but still weaker than the Hodge decomposition of Theorem \ref{decomposition thm}. However, neither does imply the symmetry \(h^{p,q}(M)=h^{q,p}(M)\) in Corollary \ref{hodge numbers}.
\end{rem}
\begin{rem}
    The technique of spectral sequence is widely used when we need to calculate a cohomology group "is determined by" two different factors (here, holomorphic and anti-holomorphic factors) and \((E_{1}^{p,q},E^{p,q}_{2},...)\) can be treated as successive approximations. Spectral sequences can be used in the theory of vector bundles, in which we need to consider how we can combine the information of the base manifold and the information of fibers to get the information of the total space. This gives rise to \textit{Serre Spectral Sequences}.
\end{rem}

\section{Basic Arithmetic Geometry}\label{Arithmetic Geometry}

Our major references in this section are \cite{aty} and \cite{Po} for Sections \ref{p adics},  \ref{basics}
and \cite{sil} for Section \ref{elliptic curves}.

\subsection{\(p\)-adic Integers}\label{p adics}
In algebraic geometry, we treat as geometrical objects zero points of polynomial equations over a field. The most straightforward choice of the base field is \(\mathbb{C}\), since a great many geometrical objects can be represented thus (for instance, \(z_{1}^{2}-z_{2}=0\) represents a Riemann surface in \(\mathbb{C}^2\)). In number theory, we can also transfer certain problems as finding solutions to polynomial equations. For instance, finding whether there is a perfect square that is equivalent to \(q\) modulo \(p\), is equivalent to finding zeros of the polynomial \(x^2-q\) in the field \(\mathbb{Z}/p\mathbb{Z}\): these problems are studied in \(\textit{Arithmetic Geometry}\). Regarding the ground field, we choose \(\mathbb{Q}\) or \(\mathbb{Q}_{p}\) instead of \(\mathbb{R}\) or \(\mathbb{C}\) to develop this theory.

The construction of the field \(\mathbb{Q}_{p}\) and its norm and topology (which gives it a geometrical significance) follows from the necessity to define a number-theoretic "distance". We can say that, for instance, we can consider an integer divisible by \(p^2\) to be closer to \(0\) than an integer only divisible by \(p\) since the former integer is not only \(0\ \mathrm{mod}\ p\), but also \(0\ \mathrm{mod}\ p^{2}\). For a prime integer \(p\), we define a topology on \(\mathbb{Z}\) with a topological basis \(\{\mathbb{Z},p\mathbb{Z},p^2\mathbb{Z}...\}\). Suppose that there is a Cauchy sequence \((x_{n})_{n}\) defined in \(\mathbb{Z}\) with respect to this topology; then the image of \(x_{n}\) in \(\mathbb{Z}/p^{m}\mathbb{Z}\) will be constant as \(\zeta_{m}\) for \(n\) large enough. The canonical projective map from \(\theta_{m+1}\) from \(\mathbb{Z}/p^{m+1}\mathbb{Z}\) to \(\mathbb{Z}/p^{m}\mathbb{Z}\) will satisfy that
\begin{equation}\label{Cauchy}
    \zeta_{m}=\theta_{m+1}\zeta_{m+1}.
\end{equation}
We obsertve that every Cauchy sequence defines a sequence \((\zeta_{m})_{m}\) such that Equation \eqref{Cauchy} holds. It is obvious that two equivalent Cauchy sequences define the same sequence that satisfies \eqref{Cauchy} and that each sequence satisfying \eqref{Cauchy} corresponds to an equivalent class of Cauchy sequences. Term-wise addition defines an Abelian group structure on equivalence classes of Cauchy sequences. It can also be interpreted as the completion of \(\mathbb{Z}\) with respect to the topology given at the beginning. The Abelian group of the equivalence classes of Cauchy sequences is called \(\textbf{p-adic integers}\) and is denoted \(\mathbb{Z}_{p}\). We can consider an inverse system with groups of the forms \(\mathbb{Z}/p^{n}\mathbb{Z}\) as objects and \(\theta_{m}\) as maps,

$$\mathbb{Z}_{p}=\underset{n}{\varprojlim}\mathbb{Z}/p^{n}\mathbb{Z}.$$

The field of fraction of \(\mathbb{Z}_{p}\) is then denoted as \(\mathbb{Q}_p\). Note that \(\mathbb{Q}_{p}\) is not algebraically closed (for instance, \(x^2+p=0\) is irreducible over \(\mathbb{Q}_{p}\)).

Since each \(x\in\mathbb{Z}_{p}\) is derived from a sequence \((\zeta_{m})_{m}\) satisfying \eqref{Cauchy}, at each point \(\zeta_{m}\) it is possible (with a small abuse of notation) to write $$\zeta_{m+1}=\zeta_{m}+p^{m}s_{m+1}$$ for some \(s_{m+1}\in\{0,1,2,...,p-1\}\). Conventionally set \(\zeta_{0}=0\). Then the \(p\)-adic integer \(x\) can be represented by a string \((s_{1},s_{2},s_{3},...)\).

\begin{ex}\label{string rep}
Let \(p=3\), then the string \((2,0,1,2,1,0,...)\) while \(s_{m}=0\) for all \(m\ge 6\) will represent \(2\times1+1\times9+2\times27+1\times81=146\).
\end{ex}

Note that given \(a\in\mathbb{Z}_{p}\), a solution of equation \(ax=1\) gives an inverse of \(a\) in \(\mathbb{Z}_{p}\). By induction, it can be shown that the inverse exists if and only if \(s_{1}\) corresponding to \(a\) is non-zero. We denote the set of invertible elements in \(\mathbb{Z}_{p}\) to be \(U^{\times}\). It is a group under multiplication. By abuse of notation, we also use \(U^{\times}\) to denote the same set in \(\mathbb{Q}_{p}\).

We further aim to grant an absolute value on \(\mathbb{Q}_{p}\), a norm that conforms to the topology defined at the beginning of this session.

\begin{Def}\label{valuation}
    For an element \(x\in\mathbb{Q}_{p}\), if \(x/{p^n}\in U^{\times}\) for some integer \(n\), then \(n\) is called \(\mathbf{valuation}\) of \(x\) in \(\mathbb{Q}_{p}\). We denote it by \(\mathrm{val}(x)=n\). By convention, \(\mathrm{val}(0)=\infty\).
\end{Def}

The concept is well defined. In particular, \(\mathrm{val}(x)>0\) if and only if \(x\) is not invertible in \(\mathbb{Z}_{p}\) and \(x\in U^{\times}\) if and only if \(\mathrm{val}(x)=0\). Therefore the valuation \(n\) is unique. We define an absolute value \(||\ ||:\mathbb{Q}_{p}\to\mathbb{R}\) by setting \(||x||=\mathrm{exp}(-\mathrm{val}(x))\).

\begin{prop}\label{norm on Qp}
    The norm \(||\ ||\) on \(\mathbb{Q}_{p}\) satisfies three conditions:
    \begin{enumerate}
    \item \(||x||\geq0\) while \(||x||=0\) if and only if \(x=0\);
    
    \item \(||x+y||\leq \mathrm{max}(||x||,||y||)\);

    \item \(||xy||=||x||\cdot||y||\).
    \end{enumerate}
\end{prop}

\begin{Def}\label{non-archimedean}
    Let K be an arbitrary field. A function \(||\ ||:K\to\mathbb{R}\) is called a \(\textbf{non-archimedean}\) \(\textbf{absolute}\) \(\textbf{value}\) if it satisfies all three conditions in Proposition \ref{norm on Qp}.
\end{Def}

Note that condition \(2\) implies that a non-archimedean absolute value satisfies the triangular inequality \(||x+y||\leq||x||+||y||\), which defines an absolute value in general.

The following theorem is an equivalence of the Newton-Raphson method in the \(p\)-adic world.

\begin{thm}\label{Hensel}
    (Hensel's lemma) Let \(f\in\mathbb{Z}_{p}[x]\). Suppose that \(f(a)\equiv 0\ (\mathrm{mod}\ p)\) and \(f'(a)\not\equiv 0\ (\mathrm{mod}\ p)\). Then there exists a unique \(b\in\mathbb{Z}_{p}\) with \(b\equiv a\ (\mathrm{mod}\ p)\) such that \(f(b)=0\).
\end{thm}

\subsection{Basic Concepts of Arithmetic Geometry}\label{basics}
In an introductory course in algebraic geometry, we start with affine space and affine varieties over an algebraically closed field. Now we consider fields that are not necessarily algebraically closed. Let \(K\) be a perfect field and let \(\overline{K}\) be the algebraic closure of \(K\).

\begin{Def}\label{affinespace}
    For \(n\in \mathbb{Z}\), \(n\geq0\), \(\mathbb{A}_{K}^{n}=\overline{K}^{n}\) is called \(n\)-dimensional \(\textbf{affine space}\) over \(K\).
\end{Def}

\begin{Def}\label{affinevariety}
    If \(T\subseteq K[x_{1},x_{2},...,x_{n}]\), then \(Z_{T}=\{P\in \mathbb{A}^{n}_{K}:f(P)=0\ for\ all\ f\in T\}\) is called an \textbf{affine variety} over \(K\).
\end{Def}
We generalize the concept of rational number to that of a rational point.
\begin{Def}\label{rationalpoint}
    An element of \(Z_{T}\) is called a \(\mathbf{K}\textbf{-rational point}\) in \(Z_{T}\) if it is an element in \(K^{n}\subseteq\overline{K}^{n}\).
\end{Def}
\begin{ex}
    The \(\mathbb{Q}\)-rational points of the equation \(x^2+y^2=1\) in \(\mathbb{Q}^{2}\) are pairs \((\frac{a}{c},\frac{b}{c})\) where \(a^2+b^2=c^2\).
\end{ex}
Any ideal \(I\) defines the same set of zeros as its \(\textbf{radical}\)
\begin{equation}
    \sqrt{I}=\{f\in K[x_{1},...,x_{n}]:f^m\in I\ for\ some\ m\geq0\}.\nonumber
\end{equation}
We have an arithmetic version of Hilbert Nullstellensatz.
\begin{thm}\label{nullstellensatz}
    There is a bijection between the set of radical ideals of \(K[x_{1},...,x_{n}]\) and the set of affine varieties \(Z_{L}\) in \(\mathbb{A}^{n}_{k}\) given by mapping a radical ideal \(I\) to the set of common zeros of \(f\in I\) and inversely, by mapping \(Z_{L}\) to the minimal radical in \(K[x_{1},...,x_{n}]\) containing \(L\).
\end{thm}
Other concepts like irreducible variety, singularity, projective space and projective variety are analogously defined in arithmetic geometry as in introductory algebraic geometry.

As in algebraic geometry in general, the morphism in arithmetic geometry is defined with respect to rational maps.
\begin{Def}\label{function field}
    If \(Z_{T}=\mathrm{Spec}\ K[x_{1},...,x_{n}]/I\) is irreducible, the \(\textbf{function field}\) \(\kappa(Z_{T})\) of \(Z_{T}\) is the fraction field of \(\mathrm{Frac}\ K[x_1,...,x_{n}]/I\).
\end{Def}
    \vskip 0.2in
    The function field of \(\mathbb{A}^{n}_{K}\) in Definition \ref{affinespace} is the called the \textbf{rational function field}.
    \vskip 0.2in
\begin{Def}\label{equivalence class}
    Let \(X\) be an irreducible variety and let \(Y\) be a projective variety in \(\mathbb{P}^{n}\). A \(\textbf{rational map}\) \(f:X\to Y\) is an equivalence class of \((n+1)\)-tuples \((f_{0}:f_{1}:...:f_{n})\) such that: 
    \begin{enumerate}
    \item \(f_{i}\in\kappa(X)\) for all \(i\) and the \(f_{i}\) are not identically \(0\);
    
    \item For any algebraic field extension \(L\supseteq K\) and any \(P\in X(L)\) such that the \(f_{i}(P)\) are all defined and not all \(0\) and that \((f_{0}(P):...:f_{n}(P))\in Y(L)\).
    \end{enumerate}
\end{Def}
Here \(X(L)\), \(Y(L)\) refer to the set defined by the same vanishing polynomials in \(L\) as \(X\) and \(Y\). Note that if \(K=\mathbb{C}\), then condition \(2\) will be redundant.

Obviously \((f_{0},f_{1},...,f_{n})\) is equivalent to \((\lambda f_{0},\lambda f_{1},...,\lambda f_{n})\) for any \(\lambda\in\kappa(X)^{\times}\). We say that \(f\) is \(\textbf{defined}\) at \(P\in X(L)\) if there exists \(\lambda\in\kappa(X)^{\times}\) such that \(\lambda f_{0}(P),\lambda f_{1}(P),...,\lambda f_{n}(P)\) are all defined.

\begin{Def}\label{mor}
    A rational map \(f:X\to Y\) that is defined at every \(P\in X(L)\) for all algebraic extension \(L\) of \(K\) is called a \(\textbf{morphism}\). \(f\) is called an \(\textbf{isomorphism}\) if it is a morphism and there exists a morphism \(g:Y \to X\) such that \(f\circ g\) and \(g\circ f\) are identities maps on \(Y\) and \(X\), respectively. 
\end{Def}

\begin{ex}
    The map \((x:y)\to (x^2:xy:y^2)\) from \(\mathbb{CP}^{2}\) to \(\mathbb{CP}^{3}\) is a morphism.
\end{ex}
\subsection{Elliptic Curves}\label{elliptic curves}
Let us introduce briefly a type of variety that appears generally in algebraic geometry but is especially important in arithmetic geometry. In this subsection let us assume \(K\) to be a field and \(\overline{K}\) its algebraic closure. Consider \(\mathbb{P}^{2}\) to be the projective curve over \(\overline{K}\) and \(O=[0,1,0]\) be its base point. We define \textbf{Weierstrass equations} by formulas of the form
\begin{equation}\label{weierstrass}
    Y^{2}Z+a_{1}XYZ+a_{3}YZ^{2}=X^{3}+a_{2}X^{2}Z+a_{4}XZ^{2}+a_{6}Z^{3}
\end{equation}
such that \(a_{1},...,a_{6}\in K\). We call the zero locus of Weierstrass equations \textbf{elliptic curves}. To ease the notation, we can write \(x=X/Z\) and \(y=Y/Z\) and keep in mind that there is an extra point of infinity. If \(\mathrm{char}(K)\neq2\) and \(\mathrm{char}(K)\neq3\), then we can further complete the squares and reduce our equation to the form 
\begin{equation}\label{simplifiedform}
    y^{2}=x^{3}-27c_{4}x-54c_{6}
\end{equation}
where we define quantities \(b_{2}=a_{1}^{2}+4a_{4}\), \(b_{4}=2a_{4}+a_{1}a_{3}\), \(b_{6}=a_{3}^{2}+4a_{6}\), \(c_{4}=b_{2}^{2}-24b_{4}\) and \(c_{6}=-b_{2}^{3}+36b_{2}b_{4}-216b_{6}\).

We now define \(b_{8}=a_{1}^{2}a_{6}+4a_{2}a_{6}-a_{1}a_{3}a_{4}+a_{2}a_{3}^{2}-a_{4}^{2}\) and \(\Delta=-b_{2}^{2}b_{8}-8b_{4}^{3}-27b_{6}^{2}+9b_{2}b_{4}b_{6}\). We call \(\Delta\) the \textbf{discriminant} of the Weierstrass equation. We have the following proposition by calculation.
\begin{prop}\label{smoothcurve}
    An elliptic curve given by Equation \eqref{simplifiedform} is smooth if and only if its discriminant is \(0\).
\end{prop}

The concept of smoothness for an algebraic variety will be elaborated on at Section \ref{smooth and etale}. Here we can just get an intuitive idea by considering \(K=\mathbb{C}\) and regarding its smoothness as the differential-geometric smoothness of the curve. By genus-degree formula we can calculate the arithmetic genus of an elliptic curve as \(g=\frac{1}{2}(d-1)(d-2)=1\) by plugging in \(d=3\).

An astonishing fact about elliptic curves is that we can give it a group structure. Let \(E\) be a smooth elliptic curve and let \(P,Q\in E\), and let \(L\) be the line through \(P\) and \(Q\). If \(P=Q\), let \(L\) be the tangent line to \(E\) at \(P\). Let \(R\) be the third point of intersection of \(L\) with \(E\). Let \(L'\) be the line through \(R\) and \(O\). Then \(L'\) intersects \(E\) at \(R\), \(O\) and a third point. We denote the third point by \(P\oplus Q\). It is easy to check that the operation \(\oplus\) gives the elliptic curve an abelian group structure with \(O=[0,1,0]\) as the identity element.

There naturally arises the question whether all smooth algebraic curves over \(K\) of genus \(1\) is an elliptic curve, that is, whether all smooth algebraic curves over \(K\) of genus \(1\) can be embedded into \(\mathbb{P}^{2}\). Indeed, it is the case.
\begin{prop}\label{embedding}
    Every smooth algebraic curve of genus \(1\) can be embedded into \(\mathbb{P}^{2}\).
\end{prop}
\begin{rem}
    The name "elliptic" comes from the relation between elliptic curves and elliptic equations. Let \(K=\mathbb{C}\), and there is an almost everywhere smooth function \(\phi\) on defined any complex torus called \textit{Weierstrass function}. It is an elliptic function and \((\phi,\phi')\) will satisfy the Weierstrass Equation \ref{embedding}. Indeed this defines a bijection between complex torus and smooth elliptic curves over \(\mathbb{C}\).
\end{rem}
\begin{rem}
    Elliptic curves are applied in number theory and arithmetic geometry. Intuitively \(\mathbb{Z}/p\mathbb{Z}\)-rational points on the elliptic curves over \(\mathbb{Z}/p\mathbb{Z}\) are solutions modulo \(p\) of certain Diophantine equations. More generally, let \(K\) be a field with \(q\) elements and let \(E\) be an elliptic curve over \(K\). Let \(K_{n}\) be the extension of \(K\) of degree \(n\). Moreover, \(\overline{K}\) will also be the algebraic closure of \(K_{n}\). The generating function of the number of \(K_{n}\)-rational points of \(E\) gives rise to a series conjectures called \(\textit{Weil Conjecture}\) that significantly inspired \(p\)-adic Hodge theory. A cohomology theory that satisfies certain axioms that can be used to prove the Weil Conjecture is called a \(\textit{Weil cohomology theory}\). Our cohomology theories in Chapter \ref{etale cohomology} and \ref{Witt vectors} are Weil cohomology theories when applied to finite fields. Especially, we can consider the number of \(K_{n}\) rational points as the number of intersections between the the diagonal of \(E\times E\) and the graph of \(q^{th}\)-power map from \(E\) to itself. This can be represented by the cup product of related Weil cohomology groups.
\end{rem}
\begin{rem}
    Consider the elliptic curve \(E\) over the field \(\mathbb{Z}/p\mathbb{Z}\) $$y^2=x^3-27c_{4}x-54c_{6}$$  and write \(a_{p}\) as the number of points within this elliptic curve. We need to consider the notion of a \textit{modular form} (a function from \(\mathbb{H}^{2}\) to \(\mathbb{C}\) satisfying certain restrictions) which bridges analysis with number theory. Especially, we can break a modular form \(f\) into Fourier series $$f(z)=\Sigma_{n=1}^{\infty}b_{n}e^{2\pi inz}.$$ Construct the generating function \(L(E,s)=\Sigma_{n=1}^{\infty}\frac{1}{n^{s}}a_{n}\) and another generating function \(L(f,s)\) called \textit{Hesse-Weil} \(L\)\textit{-function}. \textit{Taniyama-Shimura theorem} implies that there is a one-to-one correspondence between elliptic curves \(E\) and modular forms \(f\) such that for a corresponding pair \((E,f)\), $$L(E,s)=L(f,s).$$

    This theorem has a significant relationship with \textit{Fermat's last theorem}. Indeed, in \(1986\) Gerhard Frey showed that if there were a solution to Fermat's equation \(a^{n}+b^{n}=c^{n}\) for \(n>2\), there would be an elliptic curve called \textit{Frey's curve} that had no corresponding modular forms. According to Taniyama-Shimura theorem, if Frey's curve existed, than Fermat's last theorem (that the Fermat's equation has no solution for \(n>1\)) must also be true. While the full Taniyama-Shimura theorem was not proved until 2001 (by Christophe Breuil, Brian Conrad, Fred Diamond, and Richard Taylor) (see \cite{bre}), Andrew Wiles proved a strong enough version of it to establish Fermat's last theorem in \(1994\) (see \cite{wil}).
\end{rem}
To pave the road for our introduction of the cohomology theories, we state the following Theorem of Mordell-Weil without proof.
\begin{thm}[Mordell-Weil]\label{MordellWeil}
    Let \(K\) be a field that is a finite extension of \(\mathbb{Q}\) and \(E\) be an elliptic curve over \(K\), then \(E\) is finitely generated as an Abelian group.
\end{thm}

\section{Inseparable Extensions and Galois Representation}\label{Galois representation}
We use \cite{rom} as the major reference in this section.

We recall that a field extension is Galois if and only if it is finite, normal, and separable. However, when we use the geometrical method to deal with number-theoretic problems, we will encounter fields with finite characteristics that are not perfect, therefore extensions that are not separable. For example, the field \(\mathbb{F}_{2}(T)\) is not perfect since the polynomial \(f(x)=x^2-T\) is irreducible with a multiple root.
\subsection{Extension of Non-Perfect Fields}\label{non-perfect fields}
We need some qualifications to apply Galois Theory to non-perfect fields.

\begin{Def}\label{separableclosure}
    Let \(K\) be a field, the \(\textbf{separable closure}\) of \(K\) is the maximal field extension \(L\) of \(K\) such that the extension \(L/K\) is separable. We write \(L=\overline{K}^{sep}\).
\end{Def}
In fact, such an extension exists and is unique up to isomorphism. It contains every separable element (that is, elements whose minimal polynomials over \(K\) are separable). Obviously, the separable closure of a field coincides with the algebraic closure if and only if the field is perfect. We use the following concepts to characterize how far a field is from being perfect.
\begin{Def}\label{degree}
    Let \(K\) be a field, and \(F/K\) be a finite extension. Let \(F^{sep}/F\) be the maximal subextension of \(F/K\) that is separable. Then \(\textbf{separable degree}\) \([F:K]_{s}\) is defined as the integer \([F^{sep}:K]\). The \(\textbf{inseparable degree}\) \([F:K]_{i}\) is defined as the integer \([F:F^{sep}]\).
\end{Def}

It follows that $$[F:K]=[F:K]_{s}[F:K]_{i}$$ and that if the extension is separable, \([F:K]_{i}=1\) and \([F:K]=[F:K]_{s}\). However, if \([F:K]_{s}=1\), then we say that the extension is \(\textbf{purely inseparable}\). Note that \(F/K\) is both separable and purely inseparable if and only if \(F=K\).

We can examine the problem more closely. A field can be non-perfect only if it has a non-zero characteristics. Let \(\mathrm{char}(K)=p\) for a field \(K\). By the finite field theory, if the minimal polynomial \(f(x)\) of \(x\) has a multiple root, then \(f(x)=q(x^{p})\) for some polynomial \(q\). By induction, it can be found that there exists an integer \(d\geq 1\) and an irreducible polynomial \(r\) such that \(f(x)=r(x^{p^{d}})\). Then every root of \(f\) has a multiplicity of \(p^{d}\). In this case, \(x^{p^d}\) is a separable element, \([K(x):K]_{s}=\mathrm{deg}(r)\) and \([K(x):K]_{i}=p^d\). In general, if \(F/K\) is a finite extension, then \(F=K(x_{1},..,x_{n})\) for some \(x_{1},..,x_{n}\in K\).

\begin{ex}
    The extension \(\mathbb{F}_{2}(T)(\alpha)/\mathbb{F}_{2}(T)\) where \(\alpha\) is a root of the polynomial \(x^2-T\) is purely inseparable with an inseparable degree \(2\). Indeed \(\alpha^2=T\), \(x^2-T=x^2-\alpha^2=(x-\alpha)^2\). We obtain $$[\mathbb{F}_{2}(T)(\alpha):\mathbb{F}_{2}(T)]_{i}=2=[\mathbb{F}_{2}(T)(\alpha):\mathbb{F}_{2}(T)]$$ and thus \([\mathbb{F}_{2}(T)(\alpha):\mathbb{F}_{2}(T)]_{s}=1\).
\end{ex}

\subsection{Absolute Galois Group and Galois Representation}
\label{Absolute Galois}
We can further apply Galois theory to even non-perfect fields. By definition \(\overline{K}^{sep}/K\) is separable. If \(\alpha\) is separable, so do all conjugates of \(\alpha\) in the algebraic closure of \(K\). Then all conjugates of \(\alpha\) will be in \(\overline{K}^{sep}\). So the extension \(\overline{K}^{sep}/K\) is also normal. It is not typically finite (for instance, \(\overline{\mathbb{Q}}^{sep}=\mathbb{C}\) is not finite over \(\mathbb{Q}\)), but we can generalize the concept of finiteness.

\begin{Def}\label{profinitegroup}
    A \(\textbf{profinite group}\) is a group that is an inverse limit of finite groups.
\end{Def}
\begin{ex}
    Recall that the group of \(p\)-adic integers \(\mathbb{Z}_{p}=\underset{n}{\varprojlim}\mathbb{Z}/p^{n}\mathbb{Z}\), so \(\mathbb{Z}_{p}\) is a profinite group.
\end{ex}

We will also give the profinite group a topology. In a profinite group \(G=\underset{i\in I}{\varprojlim}G_{i}\) where each \(G_{i}\) is finite, we give each \(G_{i}\) a discrete topology and then \(\underset{i\in I}{\Pi}G_{i}\) a product topology. \(G\) is then a closed subset of \(\underset{i\in I}{\Pi}G_{i}\). So \(G\) can be granted a subspace topology. By Tychonoff's theorem \(\underset{i\in I}{\Pi}G_{i}\) is compact and so does \(G\) as a closed subspace of a compact space. We call this topology \(\textbf{profinite topology}\).

Consider the group \(\mathrm{Aut}(\overline{K}^{sep}/K)\), it might not be finite, but there is an inverse system composed of the groups \(\mathrm{Aut}(K'/K)\) (where \(K'/K\) is a finite subextension of \(\overline{K}^{sep}/K\)). If \(K'_{1}/K\) is a subsextension of \(K'_{2}/K\), the restriction map \(\mathrm{Aut}(K'_{2}/K)\to \mathrm{Aut}(K'_{1}/K)\) will give an arrow of this inverse system. Then \(\mathrm{Aut}(\overline{K}^{sep}/K)\) will be the inverse limit of this system. Therefore, it is profinite. We can generalize the definition of Galois group to profinite extensions.

\begin{Def}\label{absolutegalois}
    If \(K\) is a field, then \(\textbf{absolute Galois group}\) of \(K\) is the group \(\mathrm{Gal}(\overline{K}^{sep}/K)\). That is, the group of automorphisms of \(\overline{K}^{sep}\) that fix \(K\).
\end{Def}

\(\textbf{Galois representation}\) is a linear representation of the absolute Galois group. Let the representation be a linear map \(\mathrm{Gal}(\overline{K}^{sep}/K)\to \mathrm{GL}_{n}(F)\) where \(F\) is a field. Assume that \(F\) is a topological field. Observe that \(\mathrm{GL}_{n}(F)\) inherits a topology as a subset of \(F^{n^2}\). From then on we assume that the representation is a continuous one.

\begin{ex}
    If \(F=\mathbb{C}\), then the Galois representation is called an \(\textbf{Artin representation}\).
\end{ex}
\begin{ex}
    If \(F=\mathbb{Q}_{l}\) for a prime integer \(l\), then the Galois representation is called an \(\textbf{l-adic}\) \(\textbf{representation}\).
\end{ex}

If \(G\) is an abelian group, and \(l\) is a prime integer, then the \(\mathbf{l}\)\(\textbf{-adic torsion}\) of \(G\), \(G[l^{n}]\) is the kernel of the multiplication-by-\(l^{n}\) map. Then all groups of the form \(G[l^{n}]\) together with the multiplication-by-\(l\) map \(G[l^{n+1}]\to G[l^n]\) define an inverse system. Define the \(\textbf{l-adic Tate}\) \(\textbf{module}\) of \(G\) to be its inverse limit:
\begin{equation}
    T_{l}(G)=\underset{n}{\varprojlim}G[l^{n}].\nonumber
\end{equation}

Assume that \(G\) is the group of all roots of unity in a separable closure \(\overline{K}^{sep}\) of a field \(K\). Note that it has an Abelian group structure inherited from the multiplicative group of \(\overline{K}^{sep}\). The \(l\)-adic Tate module of \(G\) is called \(\textbf{the}\) Tate module. For \(z\in\mathbb{Z}_{l}\), set \(\gamma_{n}\in G[l^n]\) to be the \(l^{n}\)-th root of unity \(T_{l}(G)\). Then a structure of a \(\mathbb{Z}_{l}\)-module is defined on \(T_{l}(G)\) by \(z\gamma_{n}=\gamma_{n}^{z\ \mathrm{mod}\ l^{n}}\); it is a free module of rank \(1\).

The inherited additive structure on \(T_{l}(G)\) makes every element of the absolute Galois group of \(K\) acting linearly on \(T_{l}(G)\) as a \(\mathbb{Z}_{l}\)-module (one can verify this by checking how these two inverse limits are constructed). Therefore the Tate module gives a Galois representation, more specifically, an \(l\)-adic representation. We call it the \(\mathbf{l}\)\(\textbf{-adic cyclotomic character}\) of \(K\).

Let \(K\) be a field and \(E\) be an elliptic curve over \(K\). It implies that \(E\) is an Abelian group, as we explained in Section \ref{elliptic curves}. We write the \(l\)-adic torsion of \(E\) as \(E[l^{n}]\). We can subsequently define \(T_{l}(E)\). Note that Theorem \ref{MordellWeil} implies that \(E\) is a direct sum of its non-torsion part and its torsion part, so the study of \(l\)-adic torsion of \(E\) is of special interests.
\section{De Rham and \'Etale Cohomology}\label{etale cohomology}
Our major references for this section are \cite{Bo1} for Sections \ref{differentials} and \ref{smooth and etale}, \cite{Bo2} for Section \ref{Grothendieck} and \cite{Vos2} for Section \ref{cohomology theories}.
\subsection{Algebraic Definition of Differentials}\label{differentials}
The algebraic and number-theoretic aspect of \(p\)-adic Hodge Theory is the study of Galois Representations while the geometric aspect of \(p\)-adic Hodge Theory involves the relationship of various cohomology theories. As we see in Proposition \ref{norm on Qp}, \(\mathbb{Q}_{p}\) has a non-archimedian absolute value. It will generate a topology with some very weird properties so that the definition of the de Rham group as a quotient group of the cotangent space is infeasible. For instance, we have:
\begin{prop}\label{weirdtopology}
    The set \(\{x\in\mathbb{Q}_{p}:||x||\leq r\}\) is open for all \(r>0\). Furthermore, the topology of \(\mathbb{Q}_{p}\) is totally disconnected. That is, any subset in \(\mathbb{Q}_{p}\) consisting of more than one point is not connected.
\end{prop}
We seek now to define a de Rham cohomology algebraically. We begin by defining the module of derivations in the affine case. 
\begin{Def}\label{derivation}
    Let \(R\) be a commutative ring with identity, \(A\) an \(R\)-algebra and \(N\) an \(A\)-module. An \(\textbf{R-derivation}\) from \(A\) to \(N\) is an \(R\)-linear map \(d:A\to M\) satisfying the Leibniz rule
    \begin{equation}
        d(fg)=fd(g)+gd(f)\nonumber
    \end{equation}
    for \(f,g\in A\). The set of \(R\)-derivations from \(A\) to \(N\) has the structure of an \(A\)-module. We denote it by \(\mathrm{Der}_{R}(A,N)\).
\end{Def}
In differential geometry, we can consider \(A\) as the \(\mathbb{R}\)-algebra of smooth functions on \(\mathbb{R}^{n}\) and consider \(N\) as the module of differential \(1\)-forms. Then the \(d\) defined in differential geometry will satisfy the conditions defined in Definition \ref{derivation}.

Returning to algebraic geometry interpretations, we have the following lemma.
\begin{lem}\label{returntoalgebraic}
    Let \(A=R[T_{i};i\in I]\) be the polynomial ring in a family of variables \((T_{i})_{i\in I}\) over \(R\) and \(N\) an \(A\)-module. Then, for each family \((x_{i})_{i\in I}\) of elements in \(N\), there is a unique \(R\)-derivation \(d: A\to N\) such that \(d(T_{i})=x_{i}\) for all \(i\in I\).
\end{lem}
For the sketch of the proof, set \(d(T_{i})=x_{i}\) for \(i\in I\). Then for \(P\in A\), we set $$d(P)=\Sigma_{i\in I}\frac{\partial P}{\partial T_{i}}x_{i},$$ while \(\frac{\partial P}{\partial T_{i}}\) is the formal partial derivative built purely algebraically.

We can define now the \(\textbf{module of relative differential forms}\) as a universal object. We first give the universal property,
\begin{prop}\label{universalproperty}
    Let \(A\) be an \(R\)-algebra. Then there exists an \(A\)-module \(\Omega^{1}_{A/R}\) together with an \(R\)-derivation \(d_{A/R}:A\to \Omega^{1}_{A/R}\) such that for every \(R\)-derivation \(d:A\to M\) there exists a unique \(A\)-linear map \(\phi:\Omega^{1}_{A/R}\to M\) such that \(d=\phi\circ d_{A/R}\).
\end{prop}
\begin{Def}\label{module}
    \(\Omega^{1}_{A/R}\) in the previous proposition is called the \(\textbf{module of relative differential}\) \(\textbf{forms}\) of \(\textbf{degree}\) \(1\) of \(A\) over \(R\) and \(d_{A/R}\) the \(\textbf{exterior differential}\) of the \(R\)-algebra \(A\). Furthermore, the \(n\)-fold exterior power \(\Omega^{n}_{A/R}=\bigwedge\nolimits^{n}\Omega^{1}_{A/R}\) for \(n\in\mathbb{N}\) is called the \(\textbf{module of relative}\) \(\textbf{differential forms}\) of \(\textbf{degree}\) \(n\) of \(A\) over \(R\).
\end{Def}

For the sketch of the proof, we first consider the free polynomial ring \(A=R[T_{i};i\in I]\) where \((T_{i})_{i\in I}\) is a family of variables. Set \(\Omega^{1}_{A/R}=A^{(I)}\). From Lemma \ref{returntoalgebraic}, there is a unique \(R\)-derivation $$d_{A/R}:A\to\Omega^{1}_{A/R} \text{ such that } d_{A/R}(T_{i})=x_{i}.$$ Therefore, \((\Omega^{1}_{A/R},d_{A/R})\) will be a pair that satisfies the following property. In particular, for any \(R\)-derivation \(d:A\to M\), construct \(\phi\) by setting \(\phi(x_{i})=d(T_{i})\) for \(i\in I\).

To generalize the definition to all schemes, we need an equivalent definition. Starting from an \(R\)-algebra \(A\), consider the map \(\mu(x\otimes y)=xy\) from \(A\otimes_{R}A\) to \(A\). Set \(\mathcal{J}=\mathrm{ker}\ \mu\). If \(x\in A\), then \(1\otimes x-x\otimes 1\in\mathcal{J}\). Let the map \(\delta:A\to\mathcal{J}/{\mathcal{J}^2}\) be the one induced by the map \(x\to1\otimes x-x\otimes1\) from \(A\) to \(\mathcal{J}/{\mathcal{J}^2}\). It can be checked that \(\mathcal{J}/{\mathcal{J}^2}\) has a structure of an \(A\)-module. By some technical calculation, it is an \(R\)-derivation. We have the following important result.
\begin{thm}\label{equivalence2}
    \(\mathcal{J}/{\mathcal{J}^2}\) is isomorphic to \(\Omega^{1}_{A/R}\) as an \(A\)-module. Moreover, \(\delta\) is isomorphic to \(d_{A/R}\) as an \(R\)-derivation to the same \(A\)-module.
\end{thm}

In other words, it is equivalent to define the aforementioned pair as \((\mathcal{J}/{\mathcal{J}^2}, \delta)\). For a heuristic argument in differential geometry, if \(A\) is analogous to the sheaf of polynomial functions on, for instance, a manifold \(X\), \(\mathcal{J}\) is an analogue of the sheaf of polynomial functions that vanishes at the diagonal \(\Delta\) of \(X\times X\), which is canonically isomorphic to \(X\) itself. Then \(\mathcal{J}^2\) is analogous to the sheaf of polynomial functions that have zeros of degree at least \(2\) on the diagonal. \(\mathcal{J}/{\mathcal{J}^2}\) is then analogous to the dual space of the normal space of \(\Delta\) as a submanifold of \(X\times X\). And the normal space is canonically isomorphic to the tangent space through the map \(v\to(v,-v)\). It conforms to the definition of the space of differential \(1\)-forms as the dual space of the tangent space.

To generalize our previous constructions to schemes, consider a base scheme \(S\) and an \(S\)-scheme \(\sigma:X\to S\). Let \(\Delta:X\to X\times_{S}X\) be the diagonal embedding. It can be proved that there exists an open subscheme \(W\subset X\times_{S}X\) such that \(\Delta(X)\subset W\) and the induced morphism \(X\to W\) is a closed immersion. Let \(\mathcal{J}\subset \mathcal{O}_{W}\) be the associated quasicoherent ideal that defines \(\Delta(X)\), then the \(\mathcal{O}_{X}\)-module \(\Omega^1_{X/S}=\Delta^{*}(\mathcal{J}/{\mathcal{J}^2})\) is called the \(\textbf{sheaf of relative differential forms}\) of \(\textbf{degree}\) \(1\). Here \(\Delta^{*}\) is the inverse image sheaf, by whose construction \(\Omega^{1}_{X/S}\) is independent of the choice of \(W\). Note that \(\Omega^{1}_{A/R}\) that we define in the affine case is written as \(\Omega^{1}_{\mathrm{Spec}(A)/\mathrm{Spec}(R)}\) in this notation as a specific case. Then consider the map \(d_{X/S}:\mathcal{O}_{X}\to\Omega^{1}_{X/S}\) mapping a section \(f\) of \(\mathcal{O}_{X}\) to the residue class of \(p_{2}^{*}(f)-p_{1}^{*}(f)\) in \(\mathcal{J}/{\mathcal{J}^2}\) where \(p_{1},p_{2}:X\times_{S}X\to X\) are projections onto the factors \(X\). Then \(d_{X/S}\) is called the \(\textbf{exterior differential}\) of \(X\) over \(S\). We can define \(\textbf{sheaf of relative differential forms}\) of \(\textbf{degree}\) \(n\) as \(\Omega^{n}_{X/S}=\bigwedge^{n}\Omega^{1}_{X/S}\) as before. Note that by the construction of the wedge product, the derivative \(d_{X/S}\) can be extended to derivative from \(d:\Omega^{n}_{X/S}\to\Omega^{n+1}_{X/S}\). Then we have the \(\textbf{algebraic de Rham complex}\) \(\Omega^{\bullet}_{X/S}\) (with a little abuse of notation)
\begin{equation}\label{algebraiccomplex}
0\to\mathcal{O}_{X}\xrightarrow{d}\Omega^{1}_{X/S}\xrightarrow{d}\Omega^{2}_{X/S}\xrightarrow{d}\dots
\end{equation}
that satisfies \(d\circ d=0\). It can be proved that both \(\Omega^{n}_{X/S}\) for each \(n\) and \(\Omega^{\bullet}_{X/S}\) has the structure of an \(\mathcal{O}_{X}\)-module.

\subsection{Smooth and \'Etale Morphisms}\label{smooth and etale}
In differential geometry, the implicit function theorem gives the condition for the zero set of a smooth function to be a locally smooth manifold. Since it is a local theorem, without loss of generality, we consider the smooth map $f:\mathbb{R}^{m+n}\to\mathbb{R}^{n}$   such that $$f(x_{1},...,x_{m+n})=0 \text{ at }(x_{1},...,x_{m+n})\in\mathbb{R}^{m+n}.$$ Then the set \(f^{-1}(0)\) is locally a smooth manifold around \((x_{1},...,x_{m+n})\) if and only if the vector space spanned by the vectors $$\{(\frac{\partial f_{i}}{\partial x_{1}},\frac{\partial f_{i}}{\partial x_{2}},...,\frac{\partial f_{i}}{\partial x_{m+n}}):\text{  for } i\in\{1,2,...,m\}\}$$ has the full rank \(m\). We need to transform this insight into algebraic terms. Let us set \(\mathbb{A}^{n}_{S}\) to be the \(\textbf{affine n-space}\) \(\mathrm{Spec}\ \mathcal{O}_{S}[t_{1},...,t_{n}]\). There is also a canonical map from \(\mathbb{A}^{n}_{S}\to S\) defined by gluing the maps \(\mathbb{A}^{n}_{S'}\to S'\) generated by the immersion \(\mathcal{O}_{S'}\to\mathcal{O}_{S'}[t_{1},...,t_{n}]\) with \(S'\) being the affine components of \(S\). Define \(\textbf{S-immersion}\) to be an \(S\)-morphism that is an immersion.

\begin{Def}\label{smooth morphism}
    A morhpism of schemes \(f:X\to S\) is called \(\textbf{smooth}\) at a point \(x\in X\) of \(\textbf{relative dimension}\) \(r\) if there is an open neighborhood \(U\subset X\) of \(x\) together with a closed \(S\)-immersion \(j:U\to W\subset\mathbb{A}^{n}_{S}\) into an open subscheme \(W\) of some affine \(n\)-space \(\mathbb{A}^{n}_{S}\) such that:
    \begin{enumerate}
    \item If \(\mathcal{I}\subset\mathcal{O}_{W}\) is the sheaf of ideals corresponding to the closed immersion \(j\), there are \(n-r\) sections \(g_{r+1},...,g_{n}\) in \(\mathcal{I}\) that generate \(\mathcal{I}\) in a neighborhood of \(z=j(x)\); in particular, we assume \(r\leq n\).
    
    \item The residue classes \(dg_{r+1}(z),...,dg_{n}(z)\in\Omega^{1}_{\mathbb{A}_{S}^{n}/S}\otimes k(z)\) of the differential forms \(dg_{r+1},...,dg_{n}\) are linearly independent over \(k(z)\).

    We say that \(f\) is \(\textbf{smooth}\) if it is smooth at any point \(x\in X\).
    \end{enumerate}
\end{Def}
Here \(k(z)=\mathcal{O}_{\mathbb{A}^{n}_{S}/S,z}/\mathfrak{m}_{z}\) where \(\mathcal{O}_{\mathbb{A}^{n}_{S}/S,z}\) is the stalk of the sheaf \(\mathcal{O}_{\mathbb{A}^{n}_{S}/S}\) at \(z\), and \(\mathfrak{m}_{z}\) is  its maximal ideal. Note that the sheaf of ideals corresponding to \(j\) is the collection of sections that vanishes in the image of \(j\). That is analogous to a set of functions whose common zeros define the image of \(j\). Condition \(2\) is analogous to the linear independence condition given by the implicit function theorem. We continue to define a concept analogous to local diffeomorphism.

\begin{Def}\label{etale}
    A morphism of schemes \(f:X\to S\) is called \(\textbf{\'etale}\) at \(x\) (respectively, \(\textbf{\'etale}\) on \(X\)), if it is smooth of relative dimension \(0\) at \(x\) (respectively, smooth of relative dimension \(0\) at all points of \(X\)).
\end{Def}

\begin{ex}
    An open immersion is \'etale as an local diffeomorphism.
\end{ex}

\subsection{Grothendieck and \'Etale Topologies}\label{Grothendieck}
To further construct a cohomology theory on \(p\)-adic spaces, we use a weaker version of topology that gives a pseudo open covering.

\begin{Def}\label{Gr topology}
    A \(\textbf{Grothendieck topology}\) \(\Sigma\) consists of a category \(\mathrm{Cat}\ \Sigma\) and a set \(\mathrm{Cov}\ \Sigma\) of families \((U_{i}\to U)_{i\in I}\) of morphisms in \(\mathrm{Cat}\ \Sigma\), called \(\textbf{coverings}\), such that the following hold:
\begin{enumerate}
    \item If \(\Phi:U\to V\) is an isomorphism in \(\mathrm{Cat}\ \Sigma\), then \((\Phi)\in \mathrm{Cov}\ \Sigma\).

    \item If \((U_{i}\to U)_{i\in I}\) and \((V_{ij}\to U_{i})_{j\in J}\) for \(i\in I\) belong to \(\mathrm{Cov}\ \Sigma\), then the same is true for the composition \((V_{ij}\to U_{i}\to U)_{i\in I,j\in J_{i}}\).

    \item If \((U_{i}\to U)_{i\in I}\) is in \(\mathrm{Cov}\ \Sigma\) and if \(V\to U\) is a morphism in \(\mathrm{Cat}\ \Sigma\), then the fiber products \(U_{i}\times_{U}V\) exist in \(\mathrm{Cat}\ \Sigma\).
    \end{enumerate}
\end{Def}

\(\mathrm{Cat}\ \Sigma\) is analogous to the collection of open sets and the morphisms of \(\mathrm{Cat}\ \Sigma\) is analogous to the collection of inclusions of open sets. The first condition says that the isomorphism from an open set to itself is also an open covering. The second condition says that open coverings have a certain "transitivity". The third condition is analogous to that if \((U_{i}\to U)_{i\in I}\) is an open covering of \(U\) and \(V\subset U\) is open, then \((U_{i}\cap V\to V)_{i\in I}\) is a covering of \(V\). Note that there is nothing guaranteeing that the union of "open sets" is also "open". Especially, if the objects of the \(\mathrm{Cat}\ \Sigma\) are certain subsets of a given set \(X\) with inclusions as morphisms and if the elements of \(\mathrm{Cov}\ \Sigma\) \((U_{i}\to U)_{i\in I}\) are coverings of \(U\) in set-theoretic sense, then we call \(X\) a \textbf{G-topological space}. We further introduce a notion of presheaf and sheaf on Grothendieck topology.

\begin{Def}\label{sheaf}
    Let \(\Sigma\) be a Grothendieck topology and \(\mathfrak{C}\) a category admitting cartesian products. A \(\textbf{presheaf}\) on \(\Sigma\) with values in \(\mathfrak{C}\) is defined as a contravariant functor \(\mathfrak{F}:\mathrm{Cat}\ \Sigma\to\mathfrak{C}\). We call \(\mathfrak{F}\) a \(\textbf{sheaf}\) if the diagram
    \begin{equation}\label{exactsequence}
        \mathfrak{F}(U)\to\prod_{i\in I}\mathfrak{F}(U_{i})\to\prod_{i,j\in I}\mathfrak{F}(U_{i}\times_{U}U_{j})
    \end{equation}
    is exact for any covering \((U_{i}\to U)_{i\in I}\) in \(\mathrm{Cov}\ \Sigma\).
\end{Def}

Here, \(\mathrm{Cat}\ \Sigma\) is a category with morphisms as arrows. Since a morphism \(U\to V\) is analogous to an open immersion, \(\mathfrak{F}(U\to V)\) is analogous to the restriction map. And the first arrow is canonically defined. As for the second arrow, if \(\sigma_{i},\sigma_{j}\) are, respectively, the \(i^{th}\) and \(j^{th}\) components, then the \((i,j)^{th}\) component of the image will be \(\mathfrak{F}(U_{i}\times_{U}U_{j}\to U_{i})(\sigma_{i})-\mathfrak{F}(U_{i}\times_{U}U_{j}\to U_{j})(\sigma_{j})\). 

We further add the \'etale condition to construct \'etale topology.

\begin{Def}\label{etale topology}
    Let \(X\) be a scheme, then the \(\textbf{\'etale topology}\) \(\Sigma\) on \(X\) is a Grothendieck topology on \(X\) such that all morphisms in \(\mathrm{Cat}\ \Sigma\) are \'etale. A scheme with an \'etale topology is called an \textbf{\'etale site}.
\end{Def}

\subsection{Constructions of Cohomology Theories}\label{cohomology theories}
We can further construct some analogues of classical cohomology theories algebraically. We come back to \(X\) in the complex \ref{algebraiccomplex}. It might be tempting to simply define de Rham cohomology by the complex \ref{algebraiccomplex}. But such a definition does not work very well. For instance, for a compact and connected complex manifold \(X\) of dimension \(n\) over \(\mathbb{C}\), differential geometry tells us that the \(2n^{th}\) de Rham cohomology group will be one-dimensional. However, if calculated by the complex \ref{algebraiccomplex}, all \(i^{th}\)-cohomology groups over \(\mathbb{C}\) with \(i>n\) will vanish. To fix this problem, we introducte the concept of a hypercohomology.

For a left-bounded complex \(J^{\bullet}\) of the category \(\mathfrak{C}\) that has sufficient injectives. We can construct, by induction, a double complex \((I^{\bullet,\bullet},D_{1},D_{2})\) can be constructed such that each \(I^{p,q}\) is injective, each \((I^{p,\bullet},D_{2})\) is an injective resolution of \(J^{p}\) and the inclusion \((N,d_{m})\to(I^{\bullet,0},D_{1})\) given by these resolutions is a morphism of complexes. Then we take \((I^{\bullet},D)\) to be the simple complex associated to \((I^{\bullet,\bullet},D_{1},D_{2})\). Then \(I^{\bullet}\) is a left-bounded complex and each \(I^{p}\) is an injective object of \(\mathfrak{C}\). The injective morphism of complexes \(i:J^{\bullet}\to I^{^{\bullet,0}}\) generates a morphism \(\phi^{\bullet}:J^{\bullet}\to I^{\bullet}\). We need a concept to characterize \(\phi^{\bullet}\).

\begin{Def}\label{quasi-isomorphism}
    A morphism \(\phi^{\bullet}\) of complexes is called a \textbf{quasi-isomorphism} if the induced morphisms \(H^{i}(\phi)\) are isomorphisms for every \(i\).
\end{Def}

It can be proved that \(\phi:J^{\bullet}\to I^{\bullet}\) is a quasi-isomorphism such that \(\phi^{p}:J^{p}\to I^{p}\) is injective for each \(p\).

Let \(\mathfrak{D}\) be another abelian category and \(\mathfrak{F}:\mathfrak{C}\to\mathfrak{D}\) be a left-exact functor. Let \(\phi^{\bullet}:N^{\bullet}\to I^{\bullet}\) be a quasi-isomorphism (which exists as we have shown above). We define the derived object as $$R^{i}\mathfrak{F}(N^{\bullet})=H^{i}(\mathfrak{F}(I^{\bullet})).$$ It can be proved to be independent of the choice of the quasi-isomorphism.

\begin{rem}
    One may wonder why \(R^{i}\mathfrak{F}\) is a derived functor. But we can consider a complex \(C^{\bullet}\) that is \(C\in\mathrm{Ob}\ \mathfrak{C}\) at the entry \(0\) but \(0\) otherwise. Then if \((I^{\bullet},i)\) is one of the injective resolution of \(C\), then \(i\) generates naturally \(i^{\bullet}:C^{\bullet}\to I^{\bullet}\) which is a quasi-isomorphism. So then \(R^{i}\mathfrak{F}(C^{\bullet})=R^{i}\mathfrak{F}(C)\). So our terminology is justified.
\end{rem}

Since we have a notion of the derived functor on complexes, we can define algebraic de-Rham cohomology as a hypercohomology. We introduce the following definition

\begin{Def}\label{hyper-cohomology}
    Let \(\mathfrak{C}\) be of sheaves of abelian groups over a topological space \(X\) and let \(\Gamma:\mathfrak{C}\to\mathbf{Ab}\) be the global function functor. Let \(C^{\bullet}\) be a complex of \(\mathfrak{C}\) such that \(C^{i}=0\) for \(i<0\). Then we define the \(i^{th}\) \textbf{hypercohomology} \(\mathbb{H}^{i}(C^{\bullet})\) to be \(R^{i}\Gamma(C^{\bullet})\).
\end{Def}
It is well defined since the sheaves of abelian groups have enough injectives.
We consider next the complex \(\Omega^{\bullet}_{X/S}\) defined in \ref{algebraiccomplex}. 
\begin{Def}\label{algebraic deRham}
    Let \(X\) be an \(S\)-scheme and let \(\Omega^{\bullet}_{X/S}\) be its algebraic de Rham complex. Then the \(i^{th}\) \(\textbf{algebraic de Rham cohomology}\) \(H^{i}_{dR}\) is defined as
    \begin{equation}
        H^{i}_{dR}(X/S)=\mathbb{H}^{i}(\Omega^{\bullet}_{X/S}).\nonumber
    \end{equation}
\end{Def}

We also set \(H^{i}_{dR}(X/S,\Omega^{j}_{X/S})\) to be \(R^{k}\Gamma(\Omega^{j}_{X/S})\).

Finally, we define an analogue of singular cohomology on a scheme equipped with an \'etale cohomology (or an \(\textbf{\'etale site}\)) \(X\).
\begin{Def}\label{etalecohomology}
    Let \(X\) be an \'etale site and let \(N\) be an \(\mathcal{O}_{X}\)-module. Let \(\Gamma(N)\) be the global section functor mapping an \(\mathcal{O_{X}}\)-module to its global section on \(X\). Then the \(i^{th}\) \(\textbf{\'etale cohomology}\) \(H^{i}_{et}(X,N)\) is defined by
    \begin{equation}
        H^{i}_{et}(X,N)=R^{i}\Gamma(N).\nonumber
    \end{equation}
\end{Def}
We state without proof that both algebraic de Rham cohomology and \'etale cohomology satisfies Poincar\'e duality and K\"unneth formula.
\begin{ex}\label{example coh}
    Let \(X\) be a smooth variety over a field \(K\). Then the Zariski topology naturally generates an \'etale topology on \(X\). Let \(l\) be a prime integer and by an abuse of notation, let \(\mathbb{Z}/{l^{n}}\) be the constant \'etale sheaf on \(X\). Then define the \(i^{th}\) \(\textbf{l-adic cohomology}\) as
    \begin{equation}
        H^{i}(X,\mathbb{Z}_{l})=\underset{n}{\varprojlim}H^{i}(X,\mathbb{Z}/{l^{n}})\nonumber
    \end{equation}
    and
    \begin{equation}
        H^{i}(X,\mathbb{Q}_{l})=\underset{n}{\varprojlim}H^{i}(X,\mathbb{Z}/{l^{n}})\otimes\mathbb{Q}\nonumber
    \end{equation}
    while the inverse limit is generated by the canonical inverse system. Note that the in this case the cohomology does not in general commute with the inverse limit.
\end{ex}
\begin{rem}
    Note the \(E[l^{n}]\) and \(T_{l}(E)\) that we defined at the end of Section \(6\). Actually \(E[l^{n}]=H_{et}(E,\mathbb{Z}/l^{n})\) and \(T_{l}(E)=H_{et}(E,\mathbb{Z}_{l})\).
\end{rem}
\section{Witt Vectors and Crystalline Cohomology}\label{Witt vectors}
The major references are \cite{ser} for \ref{Valuation}, \ref{extension}, \ref{structure}, \ref{witt} and \ref{more} and \cite{berth} for \ref{crystalline}. The major parts of the theory is developed by Ernst Witt, Alexander Grothendieck and Pierre Deligne.

We know that in complex geometry, de-Rham Theorem couples a closed \(n\)-form \(\phi\) with an \(n\)-simplex \(\eta\) by the integration \(\int_{\eta}\phi\). Now in \(p\)-adic cohomology theory, we have the algebraic de-Rham cohomology as an analogue of geometrical de-Rham cohomology, and we have \(l\)-adic cohomology as an analogue of singular cohomology. Now, consider Example \ref{example coh}) and assume that \(X=\overline{\mathbb{Z}/p\mathbb{Z}}\) as a smooth variety over \(\overline{\mathbb{Z}/p\mathbb{Z}}\) itself. Then if \(p\neq l\) the comparison can go smoothly (the basic theories will be discussed in the next section). But if \(p=l\) the de-Rham cohomology group will not have enough information for the comparing.

Indeed, to have a more satisfying comparison theorem we need an approach to "thicken" the space \(X\) to a characteristic \(0\) space.

\subsection{Valuation and Residue Field}\label{Valuation}
Recall that in Definition \ref{valuation} we gives a \textit{valuation} of an element \(x\in\mathbb{Q}_{p}\). We can generalize the concept of valuation as followings.
\begin{Def}\label{valuation2}
    Let \(R\) be a commutative ring with an identity. 
    \begin{enumerate}
    \item A \textbf{valuation} on \(R\) is a map \(\mathrm{val}:R\to\mathbb{R}\cup\{\infty\}\) such that:
\begin{enumerate}
    \item \(\mathrm{val}(x)=\infty\) if and only if \(x=0\).

    \item \(\mathrm{val}(x+y)\geq\mathrm{min}(||x||,||y||)\).

    \item \(\mathrm{val}(xy)=x+y\).
    \end{enumerate}
    \item Especially, it is called a \textbf{discrete valuation} if \(\mathrm{val}(x)\in\mathbb{Z}\cup\{\infty\}\) for all \(x\in R\).
    \end{enumerate}
\end{Def}
We further need to define a specifically fruitful concept of a discrete valuation ring.
\begin{Def}\label{discrete valuation ring}
    Let \(R\) be an integral domain and let \(K\) be its field of fractions. \(R\) is called a \textbf{discrete valuation ring} if \(K\) has a discrete valuation \(\mathrm{val}\) such that \(\mathrm{val}(K)=\mathbb{Z}\cup\{\infty\}\) and \(R=\{x\in K\ |\ \mathrm{val}(x)\geq0\}\).
\end{Def}
Note that Definition \ref{valuation2} gives us a commutative ring with discrete valuation. It does not automatically gives a discrete valuation ring as defined in Definition \ref{discrete valuation ring}. For example, we can give \(\mathbb{Z}\) a discrete valuation \(\mathrm{val}\) by setting \(\mathrm{val}(x)=n\) if \(x/2^{n}\) is an odd integer. But since  \(\mathbb{Z}\neq\{x\in \mathbb{Q}\ |\ \mathrm{val}(x)\geq0\}\), its does not make \(\mathbb{Z}\) a discrete valuation ring.

If \(R\) is a discrete valuation ring, then \(\mathrm{val}(1)=0\) by Part \(1(c)\) in Definition \ref{valuation2}. If \(x\in R\) is invertible in \(R\), we have \(\mathrm{val}(x)=\mathrm{val}(x^{-1})=0\). Otherwise, if \(\mathrm{val}(y)>0\) for \(y\in R\), then \(\mathrm{val}(y^{-1})<0\), so \(y\notin R\) by Definition \ref{discrete valuation ring}. So \(x\) in invertible in \(R\) if and only if \(\mathrm{val}(x)=0\) and \(R\) is a local ring with a maximal ideal \(\mathfrak{m}=\{x\in R\ |\ \mathrm{val}(x)>0\}\). We then consider an \(\pi \in R\) such that \(\mathrm{val}(\pi)=1\), \(\pi\) is then a generator of the maximal ideal \(\mathfrak{m}\). We then call \(\pi\) a \textbf{uniformizer} of the discrete valuation ring \(R\).

Recall that Proposition \ref{norm on Qp} gives the definition of a non-archimedian absolute value. Then there is a one-to-one correspondence between a valuation and a non-archimedian absolute value on a field \(K\), given by \(\mathrm{val}(x)=-\mathrm{ln}(||x||)\) and \(||x||=\mathrm{exp}(-\mathrm{val}(x))\).

Then we define some related concepts.
\begin{Def}\label{valuation ring}
    Let \(K\) be a field and with a valuation \(\mathrm{val}\). Define the \textbf{valuation ring} of \(K\) to be \(\{x\in K\ |\ \mathrm{val}(x)\geq0\}\) with its inherited ring structure.
\end{Def}

If \(R\) is a discrete valuation ring with identity granted with a valuation, then its maximal ideal will be \(\{x\in R\ |\ \mathrm{val}(x)>0\}\). We have the following concept.

\begin{Def}\label{residue ring}
    Let \(R\) be a discrete valuation ring and \(\mathfrak{m}\) be its maximal ideal. Then its \textbf{residue field} is defined as \(R/\mathfrak{m}\).
\end{Def}
\begin{ex}\label{consider valuation}
    Consider \(\mathbb{Q}_{p}\) with valuation defined in Definition \ref{valuation}. It is a discrete valuation. The valuation ring of it is \(\mathbb{Z}_{p}\) and \(\mathbb{Z}/p\mathbb{Z}=\mathbb{Z}_{p}/p\mathbb{Z}_{p}\) is the residue field of its valuation ring.
\end{ex}
\begin{rem}
    By Example \ref{consider valuation}, we know that valuation can be interpreted as an extension of the concept of prime ideals.
\end{rem}
\subsection{Extensions}\label{extension}
Here we use \(A_{\mathfrak{p}}\) to denote the localization of \(A\) at a prime ideal \(\mathfrak{p}\) and use \(\mathfrak{a}A_{\mathfrak{p}}\) for the image of an ideal \(\mathfrak{a}\subseteq A\) in the natural map \(A\to A_{\mathfrak{p}}\). In this section we work on Dedekind domains. The definition is given below.
\begin{Def}\label{Dedekind domain}
    A Noetherian integral domain \(A\) is called a \textbf{Dedekind domain} if for every prime ideal \(\mathfrak{p}\neq 0\) of \(A\), \(A_{\mathfrak{p}}\) is a discrete valuation ring.
\end{Def}

We need another concept for our further developments.
\begin{Def}\label{fractional ideal}
    Let \(R\) be an integral domain and \(K\) its field of fractions, then \(K\) has an induced structure of an \(R\)-module. A \textbf{fractional ideal} of \(R\) is an \(R\)-submodule \(I\) of \(K\) such that there exists a non-zero \(r\in R\) so that \(rI\subseteq R\).
\end{Def}

Since a Dedekind domain \(A\) is Noetherian, there can be only finitely many prime ideals containing \(x\in A\) if \(x\neq 0\). We denote \(v_{\mathfrak{p}}\) be the valuation of \(K\) generated by \(A_{\mathfrak{p}}\). Then for every \(x\in K-\{0\}\) \(v_{\mathfrak{p}}(x)\) are zero except for a finite number. If \(\mathfrak{a}\) is a fractional ideal of \(A\), the image \(\mathfrak{a}_{\mathfrak{p}}\) of \(\mathfrak{a}\) in \(A_{\mathfrak{p}}\) has the form \((\mathfrak{p}A_{\mathfrak{p}})^{v_{\mathfrak{p}}(\mathfrak{a})}\) since \(A_{\mathfrak{p}}\) is a discrete valuation ring. Here \(v_{\mathfrak{p}}(\mathfrak{a})\) are zero except for a finite number. We call \(v_{\mathfrak{p}}(\mathfrak{a})\) the \textbf{valuation} of the ideal \(\mathfrak{a}\) with respect to \(\mathfrak{p}\).

\begin{prop}\label{fractional ideal}
    Let \(A\) be a Dedekind domain. Every fractional ideal \(\mathfrak{a}\) of \(A\) can be written uniquely in the form
    \begin{equation}\label{decomposition}
        \mathfrak{a}=\prod\mathfrak{p}^{v_{\mathfrak{p}}(\mathfrak{a})},
    \end{equation}
    where the \(v_{\mathfrak{p}}(\mathfrak{a})\) are zero except for a finite number.
\end{prop}
\begin{ex}
  We observe that  \(\mathbb{Z}\) is Dedekkind as \(\mathbb{Z}_{p}\) for each prime integer \(p\) has a nonnegative discrete valuation. Its field of fractions is \(\mathbb{Q}\) with \(q\mathbb{Q}\) for any \(q\in\mathbb{Q}\). It can be written in the form in \eqref{fractional ideal} by an extension of the prime decomposition.
\end{ex}
\begin{ex}
    Let \(X\) be a smooth algebraic curve over \(\mathbb{C}\), then the ring of rational function \(\mathcal{O}_{X}\) of \(X\) is a Dedekkind domain.
\end{ex}
Let \(K\) be a field and \(F\) a finite extension of \(K\). Assume that \(K\) is the field of fractions of the a Noetherian integrally closed domain \(A\). Let \(B\) be the integral closure of \(A\) in \(F\). We can prove that \(KB=F\) and \(F\) is the field of fractions of \(B\). Assume further that \(B\) is a finitely generated \(A\)-module and \(A\) is Dedekind.

If \(\mathfrak{B}\) is a non-zero prime ideal of \(B\) and if \(\mathfrak{p}=\mathfrak{B}\cap a\), we say that \(\mathfrak{B}\) \textbf{divides} \(\mathfrak{p}\) and we will write \(\mathfrak{B}|\mathfrak{p}\). Set \(e_{\mathfrak{B}}=v_{\mathfrak{B}}(\mathfrak{pB})\), we will have
\begin{equation}
    \mathfrak{pB}=\Pi_{\mathfrak{B}|\mathfrak{p}}\mathfrak{B}^{e_{\mathfrak{p}}}.\nonumber
\end{equation}
Then \(e_{\mathfrak{B}}\) is called the \textbf{ramification index} of \(\mathfrak{B}\) in the extension \(F/K\).

On the other hand, if \(\mathfrak{B}\ |\ \mathfrak{p}\), the field \(B/\mathfrak{B}\) is an extension of the field \(A/\mathfrak{p}\). As \(B\) is finitely generated over \(A\), \(B/\mathfrak{B}\) is a finite extension of \(A/\mathfrak{p}\). The degree of extension \(f_{\mathfrak{B}}=[B/\mathfrak{B}:A/\mathfrak{p}]\) is called the \textbf{residue degree} of \(\mathfrak{B}\) in the extension \(F/K\).

We say that \(F/K\) is \textbf{totally ramified} at \(\mathfrak{p}\) if there is only one prime ideal \(\mathfrak{B}\) which divides \(\mathfrak{p}\) and \(f_{\mathfrak{B}}=1\). We say that \(F/K\) is \textbf{unramified} at \(\mathfrak{B}\) if \(e_{\mathfrak{B}}=1\) and \(B/\mathfrak{B}\) is separable over \(A/\mathfrak{p}\). If \(F/K\) is unramified for all prime ideals \(\mathfrak{B}\) dividing \(\mathfrak{p}\), we say that \(F/K\) is \textbf{unramified at} \(\mathfrak{p}\). We say that \(F/K\) is \textbf{ramified at} \(\mathfrak{p}\) if it is not unramified at \(\mathfrak{p}\).

If \([F:K]=n\), it can be proved by some technical operation that \(B/\mathfrak{p}B\) is an algebra of dimension \(n\) over \(A/\mathfrak{p}\). It can also be proved that there is an isomorphism between \(B/\mathfrak{p}B\) and \(\Pi_{\mathfrak{B}|\mathfrak{p}}B/\mathfrak{B}^{e_{\mathfrak{B}}}\) by the Chinese Remainder Theorem. We therefore have $$[B/\mathfrak{p}B:A/\mathfrak{p}]=\sum_{\mathfrak{B}|\mathfrak{p}}[B/\mathfrak{B}^{e_{\mathfrak{B}}}:A/\mathfrak{p}]=\sum_{\mathfrak{B}|\mathfrak{p}}\sum_{0\leq i\leq e_{\mathfrak{B}}-1}[\mathfrak{B}^{i}/\mathfrak{B}^{i+1}:A/\mathfrak{p}]=e_{\mathfrak{B}}f_{\mathfrak{B}}.$$ 
We sketched the proof of the following proposition.
\begin{prop}\label{dimension}
    Let \(F\), \(K\), \(A\), and \(B\) be the same integral domains as were previously denoted in this subsection. Let \(\mathfrak{p}\) be a non-zero prime ideal of \(A\), then we have
    \begin{equation}\label{dimensionsummation}
        [F:K]=\sum_{\mathfrak{B}|\mathfrak{p}}e_{\mathfrak{B}}f_{\mathfrak{B}}.
    \end{equation}
\end{prop}
\begin{rem}
    In number theory, we extend the integer ring \(\mathbb{Z}\) in order to solve Diophantine equations. Indeed, the problem of solving Diophantine equations can be studied by studying the prime factorization of factional ideals of the ring of integers of some finite extension of \(\mathbb{Q}\), that is, some "extension" of \(\mathbb{Z}\). The simplest example of such a ring of integers is the ring of Gaussian integers \(\mathbb{Z}[i]\). For a prime integer \(p\neq2\), if we want to study the Diophantine equation \(a^2+b^2=p\), we pass to study the prime factorization of \(p\) in \(\mathbb{Z}[i]\), since the equation \(a^2+b^2=p\) can be expressed as \(p=(a+bi)(a-bi)\) in \(\mathbb{Z}[i]\).

    Note that the ring of integers \(\mathcal{O}_{F_{\mathbb{Q}}}\) of some finite extension \(F_{\mathbb{Q}}\) of \(\mathbb{Q}\) is the integral closure of the Noetherian integrally closed domain \(\mathbb{Z}\) in \(F_{\mathbb{Q}}\). We aim at studying the decomposition of a integer \(m\) in \(\mathcal{O}_{F_{\mathbb{Q}}}\). But the difficulty lies in the fact that \(\mathcal{O}_{F_{\mathbb{Q}}}\) is not necessarily a unique factorization domain, so the decomposition of \(m\), if studied merely as the decomposition of "number" in the ring of integers, is not necessarily unique. For instance, in 
   
    $$\mathcal{O}_{\mathbb{Q}[\sqrt{-5}]}\cong\mathbb{Z}+\mathbb{Z}[\sqrt{-5}] \text{ since  }
    21=3\cdot7=(1+2\sqrt{-5})(1-2\sqrt{-5}).$$

    Therefore the concept of ideal is used to generalize the concept of numbers to guarantee the uniqueness of factorization. Indeed, it can be proved that any such \(\mathcal{O}_{F_{\mathbb{Q}}}\) is a Dedekind domain and Proposition \ref{fractional ideal} can be applied to it.

    It is expected that Fermat's Last Theorem can be proved by extending the number field. But the project was obsoleted and it was proved by Andrew Wiles through elliptic curves (see Section \ref{elliptic curves} instead.
\end{rem}
\begin{ex}
    Let \(A=\mathbb{Z}\), then \(K=\mathbb{Q}\). Let \(F=\mathbb{Q}[i]\), then \(B=\mathbb{Z}[i]\). 
    \begin{enumerate}
   \item Let \(\mathfrak{p}=2\mathbb{Z}\). Then \((1+i)\) is a prime ideal in \(\mathbb{Z}[i]\) since
    
     $$2\mathbb{Z}[i]=((1+i)\mathbb{Z}[i])^{2} \text{ in } \mathbb{Z}[i], \text{ and } [\mathbb{Z}[i]/(1+i)\mathbb{Z}[i]:\mathbb{Z}/2\mathbb{Z}]=1$$
     
     We have that \(\mathbb{Q}[i]/\mathbb{Q}\) is totally ramified at \(2\mathbb{Z}[i]\).
\vskip 0.5pt
   \item Let now \(\mathfrak{p}=3\mathbb{Z}\). Both \((2+i)\mathbb{Z}[i]\) and \((2-i)\mathbb{Z}[i]\) are prime ideals in \(\mathbb{Z}[i]\). Moreover, 
   $$3\mathbb{Z}[i]=((2+i)\mathbb{Z}[i])((2-i)\mathbb{Z}[i] \text{ in } \mathbb{Z}[i].$$
   
   We can check routinely other conditions and see that \(\mathbb{Q}[i]/\mathbb{Q}\) is unramified at \(3\mathbb{Z}\).
   \end{enumerate}
\end{ex}
\begin{rem}
    Set \(X=\mathrm{Spec}\ B\) and \(Y=\mathrm{Spec}\ A\). The inclusion of of \(A\) in \(B\) naturally generates a map \(f:X\to Y\). The decomposition \(\mathfrak{pB}=\Pi_{\mathfrak{B}|\mathfrak{p}}\mathfrak{B}^{e_{\mathfrak{p}}}\) gives mappings from \(\mathfrak{B}\)s to \(\mathfrak{p}\) by \(f\). Then according to Proposition \ref{fractional ideal}, there are \([F:K]\) points above \(\mathfrak{p}\), except when \(F/K\) is ramified at \(\mathfrak{p}\).

    If \(F/K\) is Galois and unramified, that \(f\) is an \([F:K]\)-fold covering of \(Y\). The maximal such extension \(K'/K\) is called the \textbf{maximal unramified extension} of \(K\). Set \(Y'=\mathrm{Spec}\ K'\), then \(f:Y'\to Y\) gives the universal covering of \(Y\). Moreover, let the \textbf{fundamental group} of \(Y\) be the group of covering transformation of the universal covering of \(Y\), that is, the group \(\mathrm{Gal}(K'/K)\). Indeed, this is the fundamental group that is associated with \'etale cohomology (see Example \ref{let y and z}).
\end{rem}
\begin{rem}
    Let \(A\) be a Dedekkind domain and \(K\) its fractional field, then we can view elements \(f\in K\) as functions on \(\mathrm{Spec}\ A\). By Proposition \ref{fractional ideal}, the group of fractional ideals can be viewed as the group of divisors on \(\mathrm{Spec}\ A\) \(Div(A)\). Moreover, if \((f)=\prod\mathfrak{p}^{v_{\mathfrak{p}}((f))}\), we can interpret it as that the function \(f\) has a zero of degree \(v_{\mathfrak{p}}((f))\) at the point \(\mathfrak{p}\). In particular, write \(P(A)\) to be the group of fractional ideals of \(A\) that is generated by a single element and define the \textit{Picard group} as \(Pic(A)=Div(A)/P(A)\). Indeed, results like Riemann-Roch formula can be extended in this context.
\end{rem}
\begin{ex}\label{let y and z}
    Let \(Y\) and \(Z\) be smooth and connected algebraic curves over \(\mathbb{C}\) and let \(f:Y\to Z\) be a finite morphism. Note that \(f\) generates an inclusion \(f^{\#}:\mathcal{O}_{Y}\to\mathcal{O}_{X}\) of Dedekind domains of rational functions. Observe that \(\mathcal{O}_{X}\) is a finitely generated \(\mathcal{O}_{Y}\)-module and \(\mathcal{O}_{Y}\) is the integral closure of \(\mathcal{O}_{X}\) in \(\mathrm{Frac}\ \mathcal{O}_{Y}\) as a finite field extension over \(\mathrm{Frac}\ \mathcal{O}_{X}\).  Let \(p\in Y\) and \(q=f(p)\) and let \(\mathfrak{m}_{p}\) and \(\mathfrak{m}_{q}\) be the associated maximal ideal. The map \(f^{\#}\) maps \(\mathfrak{m}_{q}\) into \(\mathfrak{m}_{p}\). The \(e_{\mathfrak{m}_{p}}=1\) if and only if \(f\) is a local isomorophism around \(p\) and otherwise \(q\) will be a branch point of \(f\). We call \(f\) an unramified map at \(p\) if \(e_{\mathfrak{m}_{p}}=1\) and we call \(f\) a globally unramified map if \(f\) is unramified at any \(p\in Y\). The concept of a locally (respectively, globally) unramified map can be generalized to all maps between schemes. We can rewrite Definition \ref{etale} to say that \(f\) is locally (respectively, globally) \'etale if and only if it is locally (respectively, globally) smooth and unramified.
\end{ex}

\subsection{Structure of Complete Discrete Valuation Rings}\label{structure}
In this section, let \(A\) be a \textit{complete} discrete valuation ring (we recall that the valuation gives a non-archimedian absolute value, which further gives a topology). Let \(K\) be its field of fractions and \(K_{res}\) its residue field. Let \(S\) be a system of representatives of \(K_{res}\) in \(A\) and \(\pi\) a uniformizer of \(A\). We have the following result by inductions.
\begin{prop}\label{convergent series}
    Every element \(a\in A\) can be written uniquely as a convergent series
    \begin{equation}\label{breaking1}
        a=\sum_{n=0}^{\infty}s_{n}\pi^{n}
    \end{equation}
    with \(s_{n}\in S\).

    Similarly, every element \(x\in K\) can be written as
    \begin{equation}\label{breaking2}
        x=\sum_{n=m}^{\infty}s_{n}\pi^{n}
    \end{equation}
    with \(s_{n}\in S\) and \(m\in\mathbb{Z}\).
\end{prop}
\begin{ex}
    If \(A=\mathbb{Z}_{p}\), then the string representation as exemplified in Example \ref{string rep} constructs such a series. Here the residue field \(\mathbb{Z}/p\mathbb{Z}\) is non-negative integers less than \(p\).
\end{ex}
By Proposition \ref{convergent series}, we have that if \(S\) is a subfield of \(K\) (necessarily isomorphic to \(K_{res}\)), the ring \(A\) may be identified with the ring \(K_{res}[[T]]\). This is possible only if \(K\) and \(K_{res}\) have the same characteristic. Indeed, we can also prove that whenever \(A\) and \(K_{res}\) have the same characteristic, there is always a system of representatives of \(K\) that is a field isomorphic to \(K_{res}\). We have the following proposition.
\begin{prop}\label{characteristic}
    Let \(A\) be a complete discrete valuation ring with residue field \(K_{res}\). Then \(A\cong K_{res}[[T]]\) if and only if \(A\) and \(K_{res}\) have the same characteristic.
\end{prop}
We further suppose that the characteristics of \(A\) and \(K_{res}\) are different. (If a ring \(A\) has an ideal \(I\) such that \(A\) and \(A/I\) do not have the same character, then \(A\) is called a ring of \textbf{mixed characteristics}.) It can only the case that \(A\) has characteristic \(0\) and \(K_{res}\) has characteristic \(p\neq0\). We can then identify \(\mathbb{Z}\) with a subring of \(A\) and \(p\in\mathbb{Z}\) with an element of \(A\). Since \(p\) goes to \(0\) in \(K_{res}\), \(\mathrm{val}(p)>0\) in \(A\). The integer \(e=\mathrm{val}(p)\) is called the \textbf{absolute ramification index} of \(A\). Note that the injection \(\mathbb{Z}\to A\) extends by continuity to an injection of the ring \(\mathbb{Z}_{p}\) of \(p\)-adic integers into \(A\). When \(K_{res}\) is a finite field with \(q=p^{f}\) elements, by Proposition \ref{convergent series} \(A\) is a free \(\mathbb{Z}_{p}\)-module of rank \(n=ef\) and \(K\) is an extension of degree \(n\) of the \(p\)-adic field \(\mathbb{Q}_{p}\). The integer \(e\) can then be interpreted as the ramification index of the extension \(K/\mathbb{Q}_{p}\). We say that \(A\) is \textbf{absolutely unramified} if \(e=1\).

\subsection{Witt Vectors}\label{witt}

Now we consider a converse problem: if we have a field \(P\) of characteristic \(p\), can we find a complete discrete valuation ring of which \(P\) is a residue field? In the case of that \(P\) is a perfect field, the answer is yes. We construct a ring of \textit{Witt Vectors} to solve that.

For a prime integer \(p\), \((X_{0},..,X_{n},...)\) a sequence of indeterminates over \(R\) where \(R\) is a commutative ring with identity, we define \textbf{Witt Vectors} as the series of polynomials
\begin{equation}\label{WittVectors}
    W_{n}=\sum_{i=0}^{n}p^{i}X^{p^{n-i}}_{i}.
\end{equation}
We have \(W_{0}=X_{0}\), \(W_{1}=X^{p}_{0}+pX_{1}\) and others. Let \((Y_{0},...,Y_{n},...)\) be another sequence of indeterminates in \(R\). We try to find polynomials \(S_{0},...,S_{n},...\) and \(P_{0},...,P_{n},...\) in the ring 
\begin{equation}
\mathbb{Z}[X_{0},...,X_{n},...;Y_{0},...,Y_{n},...]\nonumber
\end{equation}
such that
\begin{equation}\label{WittSummation}
    W_{n}(S_{0},...,S_{n},...)=W_{n}(X_{0},...,X_{n},...)+W_{n}(Y_{0},...,Y_{n},...)
\end{equation}
and
\begin{equation}\label{WittMultiplication}
    W_{n}(P_{0},...,P_{n},...)=W_{n}(X_{0},...,X_{n},...)\cdot W_{n}(Y_{0},...,Y_{n},...)
\end{equation}
and then there will be a ring structure through Witt polynomials on \(R^{\mathbb{N}}\). We call a vector \((a_{0},...,a_{n},...)\) a \textbf{Witt vector} with coefficients in \(R\).

To find \(S_{n}\) and \(R_{n}\), we first note that if \(\mathbb{Z}'=\mathbb{Z}[\frac{1}{p}]\), then by induction we can give \(S_{n}\) and \(P_{n}\) for all \(n\in \mathbb{N}\) uniquely as polynomials in \(\mathbb{Z}'\). By some technical operations we can then prove that the coefficients are actually in \(\mathbb{Z}\). For instance, we have $$S_{0}=X_{0}+Y_{0}$$ and $$S_{1}=X_{1}+Y_{1}+\frac{X_{0}^{p}+Y_{0}^{p}-(X_{0}+Y_{0})^{p}}{p}.$$ Also we have $$P_{0}=X_{0}Y_{0}$$ and $$P_{1}=Y_{0}^{p}X_{1}+Y_{1}X_{0}^{p}+pX_{1}Y_{1}.$$ Observe that although there is a dividing over \(p\) in the expression of \(S_{1}\), $$X_{0}^{p}+Y_{0}^{p}-(X_{0}+Y_{0})^{p}$$ is a polynomial whose coefficients are divisibe by \(p\). So our expression of \(S_{1}\) gives a polynomial with integer coefficients. We conclude the following:
\begin{thm}\label{rewritingsequence}
    There exist unique sequences \((S_{0},...,S_{n},...)\) and \(P_{0},...,P_{n},...\) of elements of the ring \(\mathbb{Z}[X_{0},...,X_{n}...;Y_{0},...,Y_{n},...]\) such that Equation \eqref{WittSummation} and Equation \eqref{WittMultiplication} hold.
\end{thm}

The set of Witt vectors \(R^{\mathbb{N}}\) is granted a structure of ring by \eqref{WittSummation} and \eqref{WittMultiplication}. It is easy to check that it is commutative and has an identity element \((1,0,...,0,...)\). It is therefore a commutative ring with an identity. We call it the \textbf{ring of Witt vectors} with \textbf{coefficients} in \(R\) and denote it \(W(R)\).
\subsection{More on the Structure of Complete Discrete Valuation Rings}\label{more}
Now we come back to our case where \(P\) is a perfect field of characteristic \(p\). We have the ring of Witt vectors \(W(P)\). Now, it can be proved that a Witt vector \((x_{0},...,x_{n},...)\in W(P)\) is a unit if and only if \(x_{0}\neq0\) (the only if direction is fairly obvious), so it is a commutative ring with identity with a maximal ideal $$\mathfrak{m}=\{(x_{0},...,x_{n},...)\in W(P)\ |\ x_{0}=0\}. $$
Therefore it has a residue field isomorphic to \(P\). Furthermore, it can be proved that \(\mathfrak{m}\) is generated by the element \(\pi=(0,1,0,...,0,...)\), which can serve as a uniformizer of a discrete valuation ring. And \(W(P)\) is also proved to be complete and absolutely unramified.

We further consider a complete discrete valuation ring \(A\) which is absolutely unramified and has \(P\) as its residue field, then we can construct (not obviously) a unique system of representation \(f:P\to A\) such that \(f(x^{p})=f(x)^{p}\) known as the \textbf{Teichm\"uller representation}. It is also a ring homomorphism. We then give a map \(\theta:W(P)\to A\) such that
\begin{equation}\label{TeichmullerRepresentation}
    \theta((x_{0},...,x_{n},...))=\sum_{i=0}^{\infty}f(x_{i})^{p^{-i}}p^{i}
\end{equation}
which can be proved to be an isomorphism. We conclude the following.
\begin{thm}\label{uniqueness}
    For every perfect field \(P\) of characteristic \(p\), \(W(P)\) is the unique (up to isomorphism) complete discrete valuation ring which is absolutely unramified and has \(P\) as its residue field. 
\end{thm}
\begin{ex}
    Let \(P=\mathbb{Z}/p\mathbb{Z}\), since \(\mathbb{Z}_{p}\) is a complete discrete valuation ring that is absolutely unramified and has \(P\) as its residue field, Theorem \ref{uniqueness} gives us that \(W(P)=\mathbb{Z}_{p}\). We can calculate the Teichm\"uller representation explicitly noting that \(x^{p}=x\) automatically in \(\mathbb{Z}/p\mathbb{Z}\) so we only need to solve the equation \(y^{p}-y=0\) in \(\mathbb{Z}_{p}\) such that \(y\equiv a\ (\mathrm{mod}\ p)\) for each \(a\in\{0,1,...,p-1\}\). That can be done by Theorem \ref{Hensel}. Let \(f\) be the Teichm\"uller representation. We notice that here \(f(x)^{p^{-i}}=f(x)\) automatically holds for all \(x\in\mathbb{Z}_{p}\) and \(i\in \mathbb{N}\), so Equation \eqref{TeichmullerRepresentation} gives an alternative expression of \(\mathbb{Z}_{p}\) as \(\sum_{i=0}^{\infty}f(x_{i})p^{i}\) where \(x_{i}\in\mathbb{Z}/p\mathbb{Z}\). With an abuse of terminology, we may call this a Teichm\"uller representation of \(\mathbb{Z}_{p}\).

    Therefore we have an algorithm to calculate sum on \(\mathbb{Z}_{p}\). We write down polynomials \(S_{n}\). Then if we want to calculate \(z_{1}+z_{2}\), while \(z_{1},z_{2}\in\mathbb{Z}_{p}\), we first write down the Teichm\"muller representation for each \(x\in\mathbb{Z}/p\mathbb{Z}\). Then we can represent \(z_{1}\) and \(z_{2}\) as in Equation \eqref{TeichmullerRepresentation} and find out respective \(x_{i}\)'s for \(z_{1}\) and \(z_{2}\). We plug them into \(S_{n}\). Finally, we output (with a little abuse of terminology) \(\sum_{i=0}^{\infty}f(S_{i}(z_{1},z_{2}))p^{i}\), using Teichm\"uler representation. The product is carried out similarly.

    The technique of Witt vectors is primarily developed as a way to simplify summation and production on \(p\)-adic integers. But it ends up as a general approach to "deepen" a perfect field of non-zero characteristics.
\end{ex}
\begin{ex}
    Let \((X_{\alpha})\) be a family of indeterminates indexed by \(\alpha\) and let \(S=\mathbb{Z}[X_{\alpha}^{p^{-\infty}}]\), which is the direct limit of the rings \(\mathbb{Z}[X_{\alpha}^{p^{-n}}]\) for all \(n\). If we give \(S\) a \(p\)-adic filtration \(\{p^{n}S\}_{n\geq0}\) and complete it, we will have the completion \(\overline{S}=\mathbb{Z}_{p}[X^{p^{-\infty}}]\). The residue ring, $T=(\mathbb{Z}/p\mathbb{Z})[X_{\alpha}^{p^{-\infty}}]$, is  perfect of characteristic \(p\). We also have \(\overline{S}=W(T)\).
\end{ex}

To complete this subsection, for a perfect field with characteristic \(p\) \(P\) and its ring of Witt vectors \(W(P)\), let \(A\) be the free \(W(P)\)-module of rank \(e\), then \(A\) can be proved to have a commutative ring structure and be the unique complete discrete valuation ring whose residue field is \(P\) and has an absolute ramification index \(e\). In particular, if \(\pi\) is a uniformizer of \(A\), then \(p=\pi^{e}\) and \(\{1,\pi,...,\pi^{e-1}\}\) is a basis of \(A\) considered as a \(W(P)\)-module. So we have the following theorem.
\begin{thm}\label{existence}
    Let \(P\) be a complete discrete valuation ring of characteristic \(p\neq0\), then the free \(W(P)\)-module of rank \(e\) is the unique (up to isomorphism) complete discrete valuation ring such that its residue field being \(P\) and the absolute ramification index is \(e\).
\end{thm}

This theorem gives the structure of any complete discrete valuation ring of characteristic \(0\) with a perfect field of non-zero character as its residue field.
\subsection{Crystalline Cohomology}\label{crystalline}
To construct crystalline cohomology we need a structure of \textit{divided power algebra} on a ring of Witt Vectors. It is a generalization of the expression \(x^{n}/n!\) widely used in Taylor expansion.
\begin{Def}\label{divided power structure}
    Let \(R\) be a commutative ring with an identity and \(I\subseteq R\) an ideal. A \textbf{divided power structure} over \(I\) is a collection of maps \(\gamma_{n}:I\to R\) for \(n\in \mathbb{N}\) that satisfies the following conditions:
\begin{enumerate}
    \item \(\gamma_{0}(x)=1\) and \(\gamma_{1}(x)=x\) for every \(x\in R\).

    \item \(\gamma_{n}(x)\in I\) if \(n\geq 1\).

    \item \(\gamma_{n}(x+y)=\sum_{i+j=n}\gamma_{i}(x)\gamma_{j}(y)\).

    \item \(\gamma_{n}(\lambda x)=\lambda^{n}\gamma_{n}(x)\) for all \(\lambda\in R\).

    \item \(\gamma_{n}(x)\gamma_{m}(x)={n+m\choose n}\gamma_{m+n}(x)\).

    \item \(\gamma_{m}(\gamma_{n}(x))=\frac{(mn)!}{m!(n!)^{m}}\gamma_{mn}(x)\).
    \end{enumerate}

    For convenience, we denote \(\gamma_{n}(x)\) as \(x^{[n]}\). We then call \(R\) a \textbf{divided power ring} and \(I\) a \textbf{divided power ideal}. In general, we denote the divided power structure as \((R,I,\delta)\).
\end{Def}
\begin{rem}\label{polynomials}
    For any commutative ring with an identity \(R\), we can construct a divided power ring \(R<x_{1},x_{2},...,x_{n}>\) that is generated by \textbf{divided power monomials} of the form \(cx_{1}^{[i_{1}]}x_{2}^{[i_{2}]}...x_{n}^{[i_{n}]}\). Here the divided power ideal is \(R<x_{1},x_{2},...,x_{n}>-R\). Especially we have the \textbf{free} divided ring algebra over \(\mathbb{Z}\) with one generator 
    \begin{equation}
    \mathbb{Z}<x>=\mathbb{Z}[x,\frac{x^2}{2},...,\frac{x^{n}}{n!},...].\nonumber
    \end{equation}
\end{rem}
\begin{ex}\label{let K}
    Let \(K\) be a perfect field of characteristic \(p\). Then we can define a divided power structure on \(W(K)\) over its maximal ideal \((p)\) defined by \(\gamma_{n}(p)=p^{n}/n!\). It can be routinely checked that all conditions in Definition \ref{divided power structure} are satisfied. We set $$W_{n}(K)=W(K)/(p^{n})$$ and see that this divided power structure naturally generates a divided power structure on \(W_{n}(K)\).
\end{ex}
 Also note that if \(f:X_{1}\to X_{2}\) is an arbitrary map of schemes and then \(f\) will give \(\mathcal{O}_{X_{1}}\) a structure of an \(\mathcal{O}_{X_{2}}\)-algebra. We will define the compatibility of divided power structures in the following ways.
 \begin{Def}\label{compatibility}
     Let \(X_{1}\) and \(X_{2}\) be schemes and \(f:X_{1}\to X_{2}\) be a scheme morphism. Let \((\mathcal{O}_{X_{1}},I,\gamma)\) and \((\mathcal{O}_{X_{2}},J,\delta)\) be divided power structures. Then we say that \((\mathcal{O}_{X_{1}},I,\gamma)\) is \textbf{compatible} with \((\mathcal{O}_{X_{2}},J,\delta)\) if it satisfies the following conditions:
     \begin{enumerate}
     \item There exists a divided power structure \((\mathcal{O}_{X_{1}},J\mathcal{O}_{X_{1}},\bar{\delta})\) on \(\mathcal{O}_{X_{1}}\) such that \(\bar{\delta}_{n}(bx)=b^{n}\delta_{n}(x)\) for any \(x\in J\) and \(b\in\mathcal{O}_{X_{1}}\).

     \item \(\bar{\delta}|_{J\mathcal{O}_{X_{1}}\cap I}=\gamma|_{J\mathcal{O}_{X_{1}}\cap I}\).
     \end{enumerate}
 \end{Def}
Now we can define a crystalline site.
\begin{Def}\label{crystalline site}
    Let \(K\) be a perfect field of characteristic \(p\). Let \(X\) be a scheme over \(K\). Then a \textbf{crystalline site}  \(\mathrm{Crys}(X/W_{n})\) is a category with the following characterizations. 
    \begin{enumerate}
        \item  Its objects are closed immersions \(U\to V\) such that \(U\) is a Zariski open set in \(X\), \(V\) is a \(W_{n}(K)\) scheme and there is a divided power structure \((\mathcal{O}_{V},\mathrm{Ker}(\mathcal{O}_{V}\to\mathcal{O}_{U}),\delta)\) compatible with the divided power structure on \(W_{n}(K)\) defined in Example \ref{let K}. 
    
\item Its morphisms goes from \((U,V,\delta)\) to \((U',V',\delta')\) are commutative diagrams formed by an open immersion \(U\to U'\) and a morphism \(V\to V'\) such that \(\delta\) is compatible with \(\delta '\).
 \end{enumerate}
\end{Def}

We further want to define a Grothendieck topology (see Definition \ref{Gr topology} on \(\mathrm{Crys}(X/W_{n})\). The coverings are given by maps \((U_{i},V_{i},\delta_{i})\to(U,V,\delta)\) such that \((V_{i}\to V)\) is an open immersion and \(U=\cup_{i}U_{i}\), \(V=\cup_{i}V_{i}\). We can also construct a sheaf \(\mathcal{O}_{X/W_{n}}\) on \(\mathrm{Crys}(X/W_{n})\) through the structural sheaf by assigning each object \((U,V,\delta)\) the Abelian group \(\mathcal{O}_{V}\). We can define a cohomology by $$H^{i}_{crys}(X/W_{n})=H^{i}(\mathrm{Crys}(X/W_{n}),\mathcal{O}_{X/W_{n}}).$$

Since \(W(K)=\underset{n}{\varprojlim}W_{n}(K)\), we can give an inverse system on the groups \(H^{i}(\mathrm{Crys}(X/W_{n}),\mathcal{O}_{X/W_{n}})\),  we can define \textbf{crystalline cohomology} by
\begin{equation}\label{crystallinecohomology}
    H^{i}_{crys}(X/W)=\underset{n}{\varprojlim}H^{i}_{crys}(X/W_{n}).
\end{equation}
\begin{ex}
    Let \(X\) be an elliptic curve, then \(H^{n}_{crys}(X/W)\) are torsion-free \(W(K)\)-modules. Especially, \(H^{1}_{crys}(X/W)\) is free of rank \(2g\) and \(H^{n}_{crys}(X/W)\cong\wedge^{n}H^{1}_{crys}(X/W)\).
\end{ex}
As previously, we state without proof that the crystalline cohomology satisfies Poincar\'e's duality and K\"unneth formula.
\section{Hodge-Tate Decomposition and Comparison Theorems}\label{htd}
The major references in this section are \cite{Fal1}, \cite{Fal2}, \cite{fon1}, \cite{fon2}, \cite{bri} for \ref{hodgetate} and \ref{fontaniesperiod}, \cite{Liu} for \ref{models}, \cite{Fal1}, \cite{Fal2}, \cite{fon1}, \cite{fon2} and \cite{xar} for \ref{comparisontheorem}. The major part of this theory is developed by Gerald Faltings and Jean-Marc Fontaine.

We need to recall the results in Sections \ref{Hodge theory} and \ref{HdR spectral sequence}. De Rham Theorem and Dolbeault Theorem establish relationships between de Rham-like cohomology theories (defined by complexes connected by differential operators) and singular-like cohomology theories (defined by Cech complexes). Also, Theorem \ref{decomposition} gives us a decomposition of cohomology groups. We now need to construct similar parallelism between different cohomology theories in \(p\)-adic case and find out whether \(p\)-adic cohomology groups have a similar decomposition. In the introductory part of Chapter \ref{Witt vectors} we mentioned the coupling between a closed \(n\)-form with an \(n\)-simplex by integration. In complex geometry, we call a \textit{period} all complex numbers that can be written as the result of a coupling integration. In complex cases we have sufficient periods so that every coupling of non-zero elements yield a non-zero result. But it is not necessarily the case in \(p\)-adic geometry. We need to build a ring with enough periods in order to state comparison theorems. We need to state Faltings original result in Hodge-Tate decomposition and then introduce Fontaine's period rings. Besides the geometric part, Galois representation is involved in number-theoretic part of the theory.
\subsection{Hodge-Tate Decomposition}\label{hodgetate}
In this section we will abuse our notations by using \(K\) to denote both a field and the scheme \((\mathrm{Spec}\ K,K)\) generated by it.

Recall our introduction of \(l\)-adic cyclotomic character at the end of Chapter \ref{Galois representation}. It is a Galois representation on \(\mathbb{Z}_{l}\) as a \(1\)-dimensional free \(\mathbb{Z}_{l}\)-module. We write it as the \(1\)-fold \(\textbf{Tate Twist}\) \(\mathbb{Z}_{l}(1)\). For \(n\in\mathbb{N}\), we define the \(n\)-fold Tate Twist \(\mathbb{Z}_{l}(n)\) to be \(\mathbb{Z}_{l}(1)^{\otimes n}\) and \(\mathbb{Z}_{l}(-n)=\mathrm{Hom}(\mathbb{Z}_{l}(n),\mathbb{Z}_{l})\). Note that a Galois representation over a free \(\mathbb{Z}_{l}\)-module of rank \(1\) naturally gives a Galois representation over \(\mathbb{Q}_{l}\). We can analogously define \(\mathbb{Q}_{l}(n)\) for all \(n\in\mathbb{Z}\).

Let \(A\) be a complete discrete valuation ring of characteristic \(0\) with a perfect residue field of characteristic \(l\). Let \(K\) be its field of fractions. It can be proved that its algebraic closure, \(\overline{K}\), is equipped with a non-archimedian absolute value (and therefore a valuation) that is an extension of that of \(K\). Let \(\mathbb{C}_{K}\) be the completion of \(\overline{K}\) under the topology defined by the absolute value. It can be proved that \(\mathbb{C}_{K}\) is an algebraically closed field itself. In particular, we set \(\mathbb{C}_{\mathbb 
Q_{l}}=\mathbb{C}_{l}\) (however, the algebraic completion of \(\mathbb{Q}_{l}\) is not itself complete). We set \(\mathbb{C}_{K}(n)=\mathbb{C}_{K}\otimes_{\mathbb{Q}_{l}}\mathbb{Q}_{l}(n)\) for all \(n\in\mathbb{Z}\).

To proceed further we need to define several concepts concerning \(K\)-schemes.
\begin{Def}\label{smoothandproper}
    Let \(K\) be a field and let \(X\) be a \(K\)-scheme with the canonical morphism of schemes \(p:X\to K\). Then \(X\) is called a \textbf{proper} \(K\)-scheme if \(p\) is separated, of finite type and for every \(K\)-scheme \(Y\), the projection from the fiber product \(X\times_{K}Y\to Y\) is a closed map of the underlying topological spaces. \(X\) is called a \textbf{smooth} \(K\)-scheme if \(p\) is smooth.
\end{Def}
\begin{ex}
    A projective variety over \(K\) is a proper \(K\)-scheme.
\end{ex}
\begin{rem}
    Properness previously defined is an analogue of the concept of a proper map in complex geometry (the inverse image of a compact set is a compact set). Affine varieties of positive dimension over \(K\) are never proper over \(K\). For instance, the affine line \(\mathbb{A}^{1}\) over \(K\) gives a projective morphism \(\mathbb{A}^{1}\times_{\mathrm{Spec}\ K}\mathbb{A}^{1}\to\mathbb{A}^{1}\) that is not closed since the closed set \(xy=1\) in \(\mathbb{A}^{1}\times\mathbb{A}^{1}\) is \(\mathbb{A}^{1}-\{0\}\), which is not closed in \(\mathbb{A}^{1}\).
\end{rem}
\begin{rem}
    Consider the algebraic de Rham cohomology constructed in Section \ref{witt}. If \(X\) is a smooth, projective \(\mathbb{C}\)-scheme, then \(X\) will also have an analytic structure as a complex manifold \(X^{an}\). We write the analytic de Rham cohomology groups over \(\mathbb{C}\) as \(H^{i}_{dR}(X^{an},\mathbb{C})\). Then we state without proof that there is a canonical isomorphism from \(H^{i}_{dR}(X/\mathbb{C})\) to \(H^{i}_{dR}(X^{an},\mathbb{C})\) such that \(H^{i}_{dR}(X/\mathbb{C})\cong H^{i}_{dR}(X^{an},\mathbb{C})\) for all \(i\in\mathbb{N}\).
\end{rem}
Fix \(X\) to be a proper and smooth scheme over \(K\). The absolute Galois group over \(K\) (recall that our \(K\) has characteristic \(0\) and is thus perfect, so \(\overline{K}^{sep}=\overline{K}\)) generate an \(l\)-adic representation over \(H^{i}_{et}(X_{\overline{K}},\mathbb{Q}_{l})\) for \(i\geq0\) as a vector space over \(\mathbb{Q}_{l}\) (here \(X_{\overline{K}}\) is the scheme \(X\times_{K}\overline{K}\)). Note that here we assume an important result (we will not prove it in this note).

\begin{thm}\label{finiteness}
    Let \(K\) be a field and let \(X\) be a proper and smooth \(K\)-scheme, then \(H^{i}_{et}(X_{\overline{K}},\mathbb{Q}_{l})\) is a finite dimensional \(\mathbb{Q}_{l}\) vector space.
\end{thm}

Also, by isometries (given by the absolute value) the action of the absolute Galois group can be uniquely extended to \(\mathbb{C}_{K}\), so there is an continuous action \(\mathrm{Gal}(\overline{K}/K)\times \mathbb{C}_{K}\to \mathbb{C}_{K}\) generated. We further consider the tensor product of vector spaces \(\mathbb{C}_{K}\otimes_{\mathbb{Q}_{l}}H^{n}_{et}(X_{\overline{K}},\mathbb{Q}_{l})\) and note that each side of the \(\otimes\) can be acted by \(\mathrm{Gal}(\overline{K}/K)\). We have a characterization of the tensor product as a \(\mathbb{C}_{K}\) representation.
\begin{Def}\label{CkRepresentation}
    A \(\mathbf{\mathbb{C}_{K}}\)\textbf{-representation} of the absolute Galois group \(\mathrm{Gal}(\overline{K}/K)\) is a finite dimensional \(\mathbb{C}_{K}\) vector space \(U\) equipped with a continuous map that is \(\mathrm{Gal}(\overline{K}/K)\times U\to U\) such that \(g(cw)=g(c)g(w)\) for any \(g\in\mathrm{Gal}(\overline{K}/K)\), \(c\in \mathbb{C}_{K}\) and \(u\in U\). The \textbf{morphism} between \(\mathbb{C}_{K}\)-representations is characterized by \(\mathbb{C}_{K}\)-linearity and \(\mathrm{Gal}(\overline{K}/K)\)-equivalence.
\end{Def}

Note that a \(\mathbb{C}_{K}\)-representation of the absolute Galois group is generally not equivalent to a linear representation since the condition that \(g(cu)=g(c)g(u)\) is different from that \(g(cu)=cg(u)\). Indeed, we can consider \(\mathbb{C}_{K}\) as an analogue of \(\mathbb{C}=\mathbb{R}[i]\). We have \(\mathrm{Gal}(\mathbb{C}/\mathbb{R}))=\mathbb{Z}/2\mathbb{Z}\) where the non-trivial element is taking the conjugacy (let us denote it by \(h\)). The vector space \(U\) is an analogue of \(\mathbb{C}^{n}\) for some \(n\in\mathbb{N}\). Then if \(h(cu)=h(c)h(u)=\overline{c}h(u)\), the operation gives an anti-linear map on \(\mathbb{C}^{n}\). Therefore a \(\mathbb{C}_{K}\) representation can be regarded as an analogue of some "anti-linear representation".

\begin{ex}
    If \(V\) is an \(l\)-adic Galois representation, then it generates a \(\mathbb{C}_{K}\)-representation structure on \(W=\mathbb{C}_{K}\otimes_{\mathbb{Q}_{l}}V\). In particular, there is a \(\mathbb{C}_{K}\)-representation structure on \(\mathbb{C}_{K}\otimes_{\mathbb{Q}_{l}}H^{n}_{et}(X_{\overline{K}},\mathbb{Q}_{l})\) and on \(\mathbb{C}_{K}(n)\).
\end{ex}

Now we can state the main result proved by Faltings (see \cite{Fal2}).

\begin{thm}\label{FaltingsDecompositionThm}
    Let \(A\) be a complete discrete valuation ring of characteristic \(0\) with a perfect residue field and let \(K\) be its field of fractions. Let \(X\) be a smooth and proper \(K\)-scheme, then there is a canonical isomorphism
    \begin{equation}\label{FaltingsDecomposition}
        \mathbb{C}_{K}\otimes_{\mathbb{Q}_{l}}H^{n}_{et}(X_{\overline{K}},\mathbb{Q}_{l})\cong\oplus_{p+q=n}(\mathbb{C}_{K}(-q)\otimes_{K}H^{p}_{dR}(X,\Omega^{q}_{X/K}))
    \end{equation}
    as \(\mathbb{C}_{K}\)-representations.
\end{thm}

Equation \eqref{FaltingsDecomposition} is called \textbf{Hodge-Tate Decomposition}. We should note the analogue between this theorem and Theorem \ref{fractional ideal} (Hodge de-Rham Decomposition). Up to this point we have not yet introduced the comparison theorem, although \eqref{FaltingsDecomposition} does illuminate some relation between de Rham and \'etale cohomology theories. To illustrate the analogue more clearly, we define the following concept.
\begin{Def}\label{HodgeTateRing}
    Let \(A\) be a complete discrete valuation ring of characteristic \(0\) with a perfect residue field and let \(K\) be its field of fractions. The \textbf{Hodge-Tate ring} of \(K\) is the \(\mathbb{C}_{K}\)-algebra \(B_{HT}=\oplus_{n\in\mathbb{Z}}\mathbb{C}_{K}(n)\) in which multiplication is defined via the natural maps \(\mathbb{C}_{K}(m)\otimes_{\mathbb{C}_{K}}\mathbb{C}_{K}(n)\to \mathbb{C}_{K}(m+n)\).
\end{Def}

Indeed it can be proved technically that Equation \eqref{FaltingsDecomposition} can be written therefore as
\begin{equation}\label{RewritingFaltings}
    B_{HT}\otimes_{K}(\oplus_{p+q=n}H^{p}_{dR}(X,\Omega^{q}_{X/K}))\cong B_{HT}\otimes_{\mathbb{Q}_{l}}H^{n}_{et}(X_{\overline{K}},\mathbb{Q}_{l})
\end{equation}
as \(B_{HT}\)-linear \(\mathrm{Gal}(\overline{K}/K)\)-equivalent isomorphisms. Here \(B_{HT}\) play the role of the coefficient ring in which the periods will be defined. We will further talk about such rings as period rings in the following subsection.
\subsection{Fontaine's Period Rings}\label{fontaniesperiod}
Let \(X\) be a \(K\)-scheme for some field \(K\). Recall our Definition of algebraic de Rham cohomology in Definition \ref{algebraic deRham} and subsequent definition of \(H^{q}_{dR}(X/K,\Omega^{p}_{X/K})\). Also recall the definition of a filtration in Chapter \ref{HdR spectral sequence}. The following result is stated without proofs.

\begin{thm}\label{filtration}
    Let \(K\) be a characteristic field \(0\) and let \(X\) be a \(K\)-scheme. There is a filtration of abelian groups on \(H^{n}_{dR}(X/K)\)
    \begin{equation}
        \cdots F^{p}H^{n}_{dR}(X/K)\xhookrightarrow{}\cdots\xhookrightarrow{}F^{0}H^{n}_{dR}(X/K)=H^{n}_{dR}(X/K)\nonumber
    \end{equation}
    such that \(F^pH^{n}_{dR}(X/K)/F^{p+1}H^{n}_{dR}(X/K)=H^{n-p}_{dR}(X,\Omega^{p}_{X/K})\).
\end{thm}

We further impose the condition that \(K\) is the field of fractions of a complete discrete valuation ring \(A\) with a perfect residue field of characteristic \(l\neq 0\) and that \(X\) is a smooth and proper \(K\)-scheme, as in the previous section. Consider Equation \eqref{RewritingFaltings} again and note that through \(B_{HT}\) we establish some relationship between \(\oplus_{p+q=n}H^{p}_{dR}(X,\Omega^{q}_{X/K})\) on the de Rham side and \(H^{n}_{et}(X_{\overline{K}},\mathbb{Q}_{l})\) on the \'etale side. Theorem \ref{FaltingsDecompositionThm}, however, does not establish that
$H^{n}_{dR}(X/K)\cong\oplus_{p+q=n}H^{p}_{dR}(X/K,\Omega^{q}_{X/K})$, as we may expect from the analogue in the classical Hodge theory for K\"ahelrian manifolds. A finer period ring is required to establish a more satisfying comparison theorem.

Recall that the \textit{ring of integers} of a field is the ring composed of all zero points of monomials with integer coefficients in that field. Let \(O_{\mathbb{C}_{K}}\) be the ring of integers on \(\mathbb{C}_{K}\). We state without proof that \(\mathbb{C}_{K}=O_{\mathbb{C}_{K}}[\frac{1}{l}]\). We define that
\begin{equation}
    R=\underset{x\to x^{l}}{\varprojlim}O_{\mathbb{C}_{K}}/(l)\nonumber
\end{equation}
with the ring structure. That is, it is the set of all series of the form \((x_{0},x_{1},...)\) such that for each \(i\), \(x_{i}\in O_{\mathbb{C}_{K}}/(l)\) and \(x_{i+1}^{l}=x_{i}\). There is a natural ring structure generated by that on \(\mathbb{C}_{K}\). Note that there is a natural ring homomorphism from \(\mathbb{Z}/l\mathbb{Z}\) to \(O_{\mathbb{C}_{K}}/(l)\) that gives the latter a \(\mathbb{Z}/l\mathbb{Z}\)-algebra structure, which extends to a \(\mathbb{Z}/l\mathbb{Z}\)-algebra structure on \(R\). It is obvious that \(R\) has a characteristic \(l\). By construction, the Frobenius map on \(R\) is surjective, and note that it is also injective since if \((x_{0},x_{1},...)^{l}=0\), then \(x_{i-1}=x_{i}^{l}=0\) for all \(i\geq 1\) and therefore \((x_{0},x_{1},...)=0\). Therefore \(R\) is perfect. There is also a natural action of the absolute Galois group \(\mathrm{Gal}(\overline{K}/K)\) on \(R\).

Indeed, our constructions of \(W\) and \(W_{n}\) in \ref{witt} and \ref{more} can be extended if \(P\) is any perfect commutative ring with an identity with characteristic \(l\). In particular, we can construct \(W(R)\) and \(W_{n}(R)=W(R)/(l^{n})\). We will also have that $$W_{0}(R)=W(R)/(lW(R))=R.$$ Besides, we can also construct Teichm\"uller representation \(f:R\to W(R)\) as in \ref{more}. Indeed every element of \(W(R)\) can be represented as \(\sum_{i=0}^{\infty}f(x_{i})l^{i}\) such that \(x_{i}\in R\) for all \(i\in\mathbb{N}\). Then we can define \(\gamma:W(R)\to O_{\mathbb{C}_{K}}\) by
\begin{equation}
    \gamma(\sum_{i=0}^{\infty}f(c_{i})l^{i})=\sum_{i=0}^{\infty}c_{i}^{(0)}l^{i}\nonumber
\end{equation}
while \(c_{i}^{(0)}\) refers to the \(0\)-th component of \(c_{i}\) as an element of \(R\) as an inverse limit. We can interpret \(\gamma\) as a lifting of the projection map \(\pi_{0}:R\to O_{\mathbb{C}_{K}}\) (to the \(0\)-th component of \(R\)) to \(W(R)\). The map is also compatible with the action of the absolute Galois group.

Our \(\gamma\) can be canonically extended to a surjective ring homomorphism
\begin{equation}
    \gamma_{\mathbb{Q}}:W(R)[\frac{1}{l}]\to O_{\mathbb{C}_{K}}[\frac{1}{l}]=\mathbb{C}_{K}\nonumber
\end{equation}
that is compatible with the absolute Galois group action. The following result is stated without proof.
\begin{prop}\label{kernel}
    Let \(R\), \(\gamma\) and \(\gamma_{\mathbb{Q}}\) be defined as above. Choose \(r\in R\) such that \(r^{(0)}=l\). Set \(\zeta=r-l\in W(R)\). Then the ideal \(\mathrm{ker}\ \zeta\subseteq W(R)\) is the principal ideal generated by \(\zeta\). Moreover, for all \(j\geq 1\),
    \begin{equation}
        W(R)\cap(\mathrm{ker}\ \gamma_{\mathbb{Q}})^{j}=(\mathrm{ker\ \gamma})^{j}\nonumber
    \end{equation}
    and
    \begin{equation}
        \cap_{j=0}^{\infty}(\mathrm{ker}\ \gamma)^{j}=0.\nonumber
    \end{equation}
\end{prop}
Note that here there is an abuse of notations: by \(r^{(0)}=l\) we mean that \(r\) is lifted to \((l,l^{\frac{1}{l}},...)\) in \(\underset{x\to x^{l}}{\varprojlim}\ O_{\mathbb{C}_{K}}\). We can then conclude that \(W(R)[\frac{1}{l}]\) injects into the inverse limit
\begin{equation}
    B^{+}_{dR}=\underset{j}{\varprojlim}\ W(R)[\frac{1}{l}]/(\mathrm{ker}\ \gamma_{\mathbb{Q}})^{j}\nonumber
\end{equation}
with a naturation absolute Galois group action that is compatible with its action on \(W(R)[\frac{1}{l}]\). It is the \(l\)-adic completion of \(W(R)[\frac{1}{l}]\). The following is stated without proof.
\begin{prop}\label{IsACompleteDiscreteRing}
    The ring \(B^{+}_{dR}\) as defined above is a complete discrete valuation ring with residue field \(\mathbb{C}_{K}\) and any generator of \(\mathrm{ker}\ \gamma_{\mathbb{Q}}\) in \(W(R)[\frac{1}{l}]\) is a uniformizer of \(B^{+}_{dR}\).
\end{prop}

We can finally define the de Rham period ring.
\begin{Def}\label{PeriodRing}
    The \textbf{de Rham period ring} is \(B_{dR}=\mathrm{Frac}(B^{+}_{dR})\) equipped with its natural \(\mathrm{Gal}(\overline{K}/K)\)-action and \(\mathrm{Gal}(\overline{K}/K)\)-stable filtration via the \(\mathbb{Z}\)-powers of the maximal ideal of \(B_{dR}^{+}\).
\end{Def}

Note that here the Frobenius automorphism of \(W(R)[\frac{1}{l}]\) does not naturally extend to \(B_{dR}^{+}\) because it does not preserve \(\mathrm{ker}\ \gamma_{\mathbb{Q}}\). For instance, let \(\zeta\) be given as in Proposition \ref{kernel}, then \(\zeta^{l}=[r^{l}]-l\notin\mathrm{ker}\ \gamma_{\mathbb{Q}}\). We will further define crystalline period rings to remedy it. From then on we fix \(\zeta\) as the symbol of a generator of the ideal \(\mathrm{ker}\ \gamma\).

Recall our definition of a divided power structure in \(7.6\). We now set \(A^{0}_{crys}\) to be the divided power structure of \(W(R)\) over \(\mathrm{ker}\ \gamma\) by defining $$A^{0}_{crys}=W(R)<\gamma>$$ using the notation in Remark \ref{polynomials}. Moreover, define \(A^{0}_{crys}\) to be the \(l\)-adic completion of \(A^{0}_{crys}\), 
$$A_{crys}=\underset{n}{\varprojlim}\ A^{0}_{crys}/(p^{n}A^{0}_{crys}).$$ We further define \(B^{+}_{crys}=A_{crys}[\frac{1}{l}]\). We further define an essential concept.

\begin{Def}\label{cry period ring}
    Let \(K\), \(R\) be defined as above. The \textbf{crystalline period ring} \(B_{crys}\) is the \(W(R)[\frac{1}{l}]\)-algebra \(B^{+}_{crys}[\frac{1}{t}]\) while $$t=\sum_{n=1}^{\infty}(-1)^{n+1}([\epsilon]-1)^{n}/n\in B^{+}_{dR}$$ and \(\epsilon\) is an element in \(R\) satisfying that \(\epsilon^{(0)}=1\) and \(\epsilon^{(1)}\neq1\).
\end{Def}
\begin{rem}\label{t is uni}
    In the definition \(\epsilon\) is constructed such that \(\epsilon^{(n)}\) is the \(l^{n}\)-th root of unity, which exists because of the algebraic closedness of \(\mathbb{C}_{K}\). The series \(\sum_{n=1}^{\infty}(-1)^{n+1}([\epsilon]-1)^{n}/n\) can be written as \(\mathrm{log}(\epsilon)\). Indeed \(t\) is a uniformizer of \(B_{dR}^{+}\).
\end{rem}
Note that we can also give \(B_{crys}\) a natural topology and \(\mathrm{Gal}(\overline{K}/K)\)-action by the constructions. We can state the following without proofs.
\begin{prop}\label{Bcrys} The crystalline period ring \(B_{crys}\) has the following properties
\begin{enumerate}
    \item \(B_{crys}\) is a subring of \(B_{dR}\).

    \item The Frobenius endomorphism on \(B_{crys}\) is injective.
    \end{enumerate}
\end{prop}
\subsection{Models and Reductions}\label{models}
Before stating the comparison theorems, we need a conceptual tool to describe, for instance, a \(\mathbb{Z}\)-scheme that preserves certain properties when transformed into the corresponding \(\mathbb{Z}/l\mathbb{Z}\)-scheme. Such a transformation is called a \textit{reduction}.

Recall our definition of Dedekind domains in Definition \ref{Dedekind domain}. The following proposition gives an alternative definition.
\begin{prop}\label{Dedekind}
    A Noetherian integral domain \(A\) is a Dedekind domain if and only if \(A\) is integrally closed of dimension \(0\) or \(1\).
\end{prop}

Then the concept can be extended to schemes in general.

\begin{Def}\label{DedekindScheme}
    A \textbf{Dedekind scheme} \(X\) is an locally Noetherian irreducible scheme of dimension \(0\) or \(1\) such that for each \(x\in X\), \(\mathcal{O}_{X,x}\) is integrally closed.
\end{Def}

One obtains the following proposition.
\begin{prop}\label{Schemeanddomain}
    Let \(X\) be a Noetherian integral scheme. Then \(X\) is a Dedekind scheme if and only if \(\mathcal{O}_{X}(U)\) is a Dedekind domain for every open subset \(U\) of \(X\). In particular, the spectrum of a Dedekind domain is a Dedekind scheme.
\end{prop}

We further give the following definition.
\begin{Def}\label{fiberedvariety}
    Let \(S\) be a Dedekind scheme. We call a flat \(S\)-scheme of finite type \(\pi:X\to S\) a \textbf{fibered variety} over \(S\). A \textbf{morphism} of fibered varieties is a morphism that is compatible with the structure of \(S\)-schemes.
\end{Def}
Therefore a model can be defined.
\begin{Def}\label{model}
    Let \(S\) be a Dedekind scheme of dimension \(1\) and let \(K\) be its function field. Let \(X\) be a \(K\)-scheme of finite type. We call a fibered variety \(\mathcal{X}\to S\) together with an isomorphism \(f:\mathcal{X}_{\eta}\to X\) as \(K\)-schemes a \textbf{model} of \(X\). In particular, we say that model is \textbf{proper} (or \textbf{smooth}) if the fibered variety is proper (or smooth). A \textbf{morphism} between two models \(\mathcal{X}\) and \(\mathcal{X}'\) is an \(S\)-scheme morphism from \(\mathcal{X}\) to \(\mathcal{X}'\) that is compatible with the isomorphisms \(\mathcal{X}\cong X\) and \(\mathcal{X}'\cong X\).
\end{Def}
Note that here \(\eta\) is a generic point on \(\mathcal{X}\) (that is, \(\eta\in\mathcal{X}\) such that the closure of \(\{\eta\}\) is \(\mathcal{X}\) itself) and \(\mathcal{X}_{\eta}\) is the generic fiber $$\mathcal{X}\times_{S}\mathrm{Spec}\ k(\eta)$$ (here \(k(\eta)\) is the residue field of \(\mathcal{O}_{\mathcal{X},\eta}\)). For instance, if \(n\) is a non-zero integer and let $$f:\mathrm{Spec}\ \mathbb{Z}[T_{1},T_{2}]/(T_{1}T_{2}^{2}-n)\to\mathbb{Z}$$ be the canonical morphism. Then \(0\) is a generic point and the generic fiber is $$\mathrm{Spec}\ \mathbb{Q}[T_{1},T_{2}]/(T_{1}T_{2}^{2}-n).$$ By defining a model, we give a scheme another scheme of which our original scheme is a "localized" scheme at the generic point.

Finally the definition of reduction can be established.
\begin{Def}\label{reduction}
    Let \(S\) be a Dedekind scheme of dimension \(1\) and let \(K\) be its field of rational functions. Let \(X\) be a \(K\)-scheme of finite type. For each closed point \(s\in S\), we call the fiber \(\mathcal{X}_{s}\) (that is, the product \(\mathcal{X}\times_{S}\mathrm{Spec}\ k(s)\) while \(k(s)\) is the fraction field of \(\mathcal{O}_{S,s}\)) of a model \(\mathcal{X}\) of \(X\) a \textbf{reduction} of \(X\) at \(s\). Moreover, if \(S=\mathrm{Spec}\ A\) for some Dedekind domain \(A\) and \(\mathfrak{p}\) is the maximal ideal of \(A\) corresponding to \(s\), we call \(\mathcal{X}_{s}\) a \textbf{reduction} of \(X\) \textbf{modulo} \(\mathfrak{p}\).
\end{Def}

We consider the reduction as taking the local fiber. In particular, we say that \(X\) has \textbf{good reduction} at \(s\in S\) if it admits a smooth and proper model over \(\mathrm{Spec}\ \mathcal{O}_{S,s}\).

\begin{ex}
    If \(A\) is a discrete valuation ring and \(K\) is its field of fractions. Let \(E\) be an elliptic curve over \(K\) given by a Weierstrass equation. Let \(\mathfrak{m}\subseteq A\) be the maximal ideal. We can take the Weierstrass equation modulo \(\mathfrak{m}\) and give an algebraic curve over \(A/\mathfrak{m}\). It is then a reduction of \(E\) modulo \(\mathfrak{m}\). It is a good reduction if it has no singular points.
\end{ex}
\subsection{Comparison Theorems}\label{comparisontheorem}
Recall the definition of a decreasing filtration in Chapter \ref{HdR spectral sequence}. We should now extend the possible index of a filtration from \(\mathbb{N}\) to \(\mathbb{Z}\). Then we have a sequence $$\cdots F^{1}I\xhookrightarrow{}F^{0}I\xhookrightarrow{}F^{-1}I\cdots$$ for an object \(I\) such that \(I=\cup_{p=-\infty}^{\infty}F^{p}(I)\). We will now consider Equation \eqref{RewritingFaltings} that relates de Rham cohomology with \'etale cohomology. Recall that the filtration in Theorem \ref{filtration} does not necessarily split, but if we tensor product the de Rham cohomology group \(H^{n}_{dR}\) with some ring \(A\) with a decreasing filtration \(F\) such that $$B_{HT}=\oplus_{p=-\infty}^{\infty}F^{p}A/F^{p+1}A,$$ then a combination with Theorem \ref{filtration} and Equation \eqref{RewritingFaltings} may suggest some better comparison theorems. Indeed, \(B_{dR}\) is such a ring. Recall that in Proposition \ref{IsACompleteDiscreteRing} we established that \(B_{dR}^{+}\) is a discrete valuation ring. Recall that \(t\) as defined in \ref{cry period ring} is a uniformizer according to Remark \ref{t is uni}. The following is stated without proofs.
\begin{prop}\label{ExistenceFibration}
    Let \(F\) be the filtration such that \(F^{p}B_{dR}=t^{p}B_{dR}\) for \(p\in\mathbb{Z}\). Then 
    \begin{equation}
        F^{p}B_{dR}/F^{p+1}B_{dR}=\mathbb{C}(-p)\nonumber
    \end{equation}
    and the isomorphism is compatible with the \(\mathrm{Gal}(\overline{K}/K)\)-action structure. In particular, \\
    \(B_{HT}=\oplus_{p=-\infty}^{\infty}F^{p}A/F^{p+1}A\).
\end{prop}
Through Proposition \ref{ExistenceFibration} de Rham Comparison Theorem can be established (it is originally established by Fontaine in \cite{fon1}).
\begin{thm}\label{Comparison1Thm}
    Let \(K\) and \(X\) be given as in Theorem \ref{FaltingsDecompositionThm}. Then 
    \begin{equation}\label{Comparison1}
    B_{dR}\otimes_{\mathbb{Q}_{l}}H^{n}_{et}(X_{\overline{K}})\cong B_{dR}\otimes_{K}H^{n}_{dR}(X/K)
    \end{equation}
    for all \(n\in\mathbb{N}\) and the isomorphism is compatible with \(\mathrm{Gal}(\overline{K}/K)\)-action, Poincar\'e duality and K\"unneth formula. Moreover, if the left hand side is equipped with the filtration generated by that on \(B_{dR}\) and the right hand side is equipped with one generated by \(F^{p}=\sum_{a+b=p}F^{a}\otimes F^{b}\), then the isomorphism is compatible with the filtration.
\end{thm}
Note that Theorem \ref{Comparison1Thm} actually implies Theorem \ref{FaltingsDecompositionThm}.

If \(K\) and \(X\) are given as in Theorem \ref{FaltingsDecompositionThm}, then \(\mathrm{Spec}\ A\) is a Dedekind scheme of dimension \(1\) and \(K\) is its field of rational functions. Let us assume that \(X\) has a good reduction at the maximal ideal \(\mathfrak{m}\subseteq A\). One has the Crystalline Comparison Theorem (established by Faltings in \cite{Fal1}).

\begin{thm}\label{Comparison2}
    If \(K\) and \(X\) are given as in Theorem \ref{FaltingsDecompositionThm} and if \(X\) has a good reduction modulo \(\mathfrak{m}\) \(\mathcal{X}\). Let \(K_{0}\) be the subfield of \(B_{crys}\) fixed by the action of the absolute Galois group, then
    \begin{equation}\label{Comparison2}
        B_{crys}\otimes H^{i}_{et}(X_{\overline{K}},\mathbb{Q}_{l})\otimes_{\mathbb{Q}_{l}}\cong B_{crys}\otimes_{K_{0}} H^{i}_{crys}(\mathcal{X}/K_{0})
    \end{equation}
    for all \(n\in\mathbb{N}\). And the isomorphism is compatible with the action \(\mathrm{Gal}(\overline{K}/K)\), the Frobenius map, the Poincar 'e duality and the K 'unneth formula.
\end{thm}

\section{Rigid Analytic Space}
The major references are \cite{Ta} and \cite{Bo2} for \ref{definitionand}, \ref{gagafunctor} and \ref{cohomologytheories} and \cite{gro} for \ref{daggerspaces}. The theory is mainly developed by John Tate in \cite{Ta}.

Chapter \ref{htd} is concerned with smooth and proper \(K\) schemes while \(K\) is a field with a non-archimedian norm. We will further consider the analytic side of the theory by defining \textit{rigid analytic spaces} over \(K\). Especially, we will have a GAGA-style theorem claiming that any smooth and proper \(K\)-scheme is a rigid analytic space. The inverse, however, is not true. So the theory in this chapter is a generalization of the theory of the previous chapters.
\subsection{Definition and Topology}\label{definitionand}
Again consider a complete discrete valuation ring \(A\) with fraction field \(K\). We again denote \(\overline{K}\) as the algebraic closure of \(K\). We mentioned previously that \(\overline{K}\) extends the non-archimedian absolute value on \(K\) and that it might not be complete, but now we state without proof that for any \textit{finite} subextension of \(\overline{K}/K\), the absolute value is nevertheless complete. For \(n\geq 1\), set \(\mathbb{B}^{n}(\overline{K})\) to be the unit ball in \(\overline{K}^{n}\). It will be the building block of the theory in this section. There is the following lemma.
\begin{lem}\label{FormalPowerSeries}
    A formal power series \(f=\sum_{\nu\in\mathbb{N}^{n}}c_{\nu}\zeta^{\nu}\in K[[\zeta_{1},...,\zeta_{n}]]\) converges on \(\mathbb{B}^{n}(\overline{K})\) if and only if \(\mathrm{lim}_{|\nu|\to\infty}||c_{\nu}||=0\).
\end{lem}
Then the analogue of "functions" on the unit ball can be defined.
\begin{Def}\label{TateAlgebra}
    The \(K\)-algebra \(T_{n}=K(\zeta_{1},...,\zeta_{n})\) of all formal power series that converges on \(\mathbb{B}^{n}(\overline{K})\) is called the \textbf{Tate algebra} of \textbf{strictly convergent} power series.
    
    By convention we write \(T_{0}=K\).
\end{Def}
Note that Lemma \ref{FormalPowerSeries} gives the description of Tate algebras. The \textbf{Gauss norm} on \(T_{n}\) is defined by setting $$||f||=\mathrm{max}||c_{\nu}|| \text{ for } f=\sum_{\nu\in\mathbb{N}^{n}}c_{\nu}\zeta^{\nu}.$$ By \ref{FormalPowerSeries} the Gauss norm on \(f\) is finite. Indeed, the following properties of Gauss norm can be verified.
\begin{prop}\label{GaussNorm}
     For each \(c\in K\) and \(f,g\in T_{n}\), the Gauss norm \(||\ ||\) satisfies the following properties:
    \begin{enumerate}
    \item \(||f||\geq0\) and \(||f||=0\) if and only if \(f=0\),
    
    \item \(||cf||=||c||\cdot ||f||\),

    \item \(||fg||=||f||\cdot||g||\),

    \item \(||f+g||\leq\mathrm{max}\{||f||,||g||\}\).
    \end{enumerate}
\end{prop}

We see that Gauss norm is actually a non-archimedian norm on \(T_{n}\). From \(3\) of the previous proposition it follows that \(T_{n}\) is an integral domain. In particular, there is the following proposition.

\begin{prop}\label{completeness}
    \(T_{n}\) is complete with respect to the Gauss norm.
\end{prop}
Now if \(\mathfrak{a}\subseteq T_{n}\) is an ideal and \(V(\mathfrak{a})\) is its zero set, then we can consider \(A=T_{n}/\mathfrak{a}\) is an algebra of "functions" on \(V(\mathfrak{a})\). We give therefore the definition of an affinoid \(K\)-algebra.
\begin{Def}\label{affinoidalgebra}
    A \(K\)-algebra \(A\) is called an \textbf{affinoid} \(\mathbf{K}\)\textbf{-algebra} if there is an epimorphism of \(K\)-algebra from \(T_{n}\to A\) for some \(n\in\mathbb{N}\). A \textbf{morphism} between two affinoid \(K\)-algebras is the homomorhpism between them as \(K\)-algebras.
\end{Def}
Note that for any epimorphism \(\alpha:T_{n}\to A\), the Gauss norm of \(||\ ||\) of \(T_{n}\) induces a norm on \(A\) called the \textbf{residue norm}, given by
$$||\alpha(f)||=\mathrm{inf}_{\alpha\in\mathrm{ker}\ \alpha}||f-\alpha||.$$

\begin{rem}
    If we view \(f\in T_{n}/\mathfrak{a}\) as \(\overline{K}\)-valued functions on the zero set \(V(\mathfrak{a})\), then we can introduce the supremum \(||f||_{sup}\) of all values that are assumed by \(f\) as the \textbf{supremum norm}. It can be proved to be an equivalent to the residue norm.
\end{rem}
We therefore have a topology on an affinoid \(K\)-algebra. We state a few properties of \(K\)-algebras without proofs.
\begin{prop}\label{affinoid characteristics}
    Let \(A=T_{n}/\mathfrak{a}\) be an affinoid \(K\)-algebra with the residue norm \(||\ ||\). The following hold:
    \begin{enumerate}
    \item The canonical map \(\alpha:T_{n}\to A\) is continuous and open,

    \item \(A\) is complete under \(||\ ||\).

    \item For any \(\overline{f}\in A\) there is an inverse image \(f\in T_{n}\) such that \(||\overline{f}||=||f||\).

    \item For any \(\overline{f}\in A\), \(||\overline{f}||=0\) if and only if \(||\bar{f}||\) is nilpotent.
    \end{enumerate}
\end{prop}

Note that if \(f\in A\) and \(x\) is a maximal ideal of an affinoid \(K\)-algebra \(A\), then its residue class \(f(x)\in A/x\). We can embed \(A/x\) into an algebraic closure \(\overline{K}\) of \(K\) and the value \(f(x)\in\overline{K}\) is defined up to an conjugation over \(K\). But its absolute value is well-defined since it does not depend on the chosen embedding. So we can treat \(A\) as "functions" on the space of maximal ideals. We write \(\mathrm{Sp}\ A\) for the set of maximal ideals of \(A\) (we write it as \(\mathrm{Max}\ A\)) with its \(K\)-algebra of "functions" \(A\) and call it the \textbf{affinoid} \(\mathbf{K}\)\textbf{-space associated to} \(A\). The Zariski topology on an affinoid \(K\)-space \(\mathrm{Sp}\ A\) is defined as usual by defining zero sets $$V(\mathfrak{a})=\{x\in\mathrm{Sp}\ A:\mathfrak{a}\subseteq x\}$$ as Zariski closed sets. The following sets form a basis of the Zariski open subsets of \(\mathrm{Sp}\ A\), $$D_{f}=\{x\in\mathrm{Sp}\ A:f(x)\neq0\} \text{ with } f\in A.$$ 

The Zariski topology, however, is too coarse for our purposes. We now introduce a finer topology. We set 

$$X(f)=\{x\in X:||f(x)||\leq 1\} \text{ for } f\in A,$$ and we define the topology generated by all sets of type \(X(f)\) the \textbf{canonical topology} of \(X\). We write $$X(f_{1},...,f_{r})=X(f_{1})\cap...\cap X(f_{r}) \text{ for } f_{1},...,f_{r}\in A.$$ We state the following without proofs.
\begin{prop}\label{openness}
    Let \(\mathrm{Sp}\ A\) be an affinoid \(K\)-space. Then for \(f\in A\), the set $$\{x\in\mathrm{Sp}\ A:f(x)\neq0\}$$ is open with respect to the canonical topology.
\end{prop}
\begin{rem}
    We see then that the canonical topology is a finer topology than Zariski topology. Moreover, recall the weird behavior of the topology generated by non-archimedian absolute values in Proposition \ref{weirdtopology} and we can understand why \(\leq\) can be put into the definition of an open set. Roughly speaking, we can consider the canonical topology as the restriction of the topology on \(K^{n}\) to \(\mathrm{Sp}\ A\).
\end{rem}
We will then define a form of subdomains.
\begin{Def}\label{AffinoidSubdomain}
    Let \(X=\mathrm{Sp}\ A\) be an affinoid \(K\)-space. A subset \(U\subseteq X\) is called an \textbf{affinoid subdomain} of \(X\) if there exists a morphism of affinoid \(K\)-spaces \(i:X'\to X\) such that \(i(X')\subseteq U\) and for any morphism of affinoid \(K\)-spaces \(\phi:Y\to X\) satisfying \(\phi(Y)\subseteq U\) admits a unique factorization through \(i:X'\to X\) via a morphism of affinoid \(K\)-spaces \(\phi':Y\to X'\).
\end{Def}
The following proposition characterizes affinoid subdomains.
\begin{prop}\label{CharacterSubdomain}
    In the situation of the previous definition, set \(X=\mathrm{Sp}\ A\) and \(X'=\mathrm{Sp}\ A'\) and let \(i^{*}:A\to A'\) be the \(K\)-morphism corresponding to \(i\). Then we have the followings:
\begin{enumerate}
    \item \(i\) is injective and satisfies \(i(X')=U\). Hence, it induces a bijection of sets between \(X'\) and \(U\).

    \item For any \(x\in X'\) and \(n\in\mathbb{N}\), the map \(i^{*}\) induces an isomorphism of affinoid \(K\)-algebras $$A/\mathfrak{m}^{n}_{i(x)}\to A'/\mathfrak{m}^{n}_{x}.$$

    \item For \(x\in X'\) we have \(\mathfrak{m}_{x}=\mathfrak{m}_{i(x)}A'\).
    \end{enumerate}
\end{prop}
We can therefore identify affinoid subdomains with the subset \(U\subseteq X\). We state without proof that \(U\) is actually an open set in canonical topology.
\begin{prop}\label{openness2}
    Let \(U\to X\) be a morphism of rigid \(K\)-spaces defining \(U\) as an affinoid subdomain of \(X\). Then \(U\) is open in \(X\), and the canonical topology of \(X\) restricts to one of \(U\).
\end{prop}

We further give a more generalized concept.
\begin{Def}\label{AdmissibleOpen}
    Let \(X\) be an affinoid \(K\)-space. A subset \(U\subseteq X\) is called \textbf{admissible open} if there is a covering \(U=\cup_{i\in I}U_{i}\) of \(U\) by affinoid subdomains \(U_{i}\subseteq X\) such that for any morphism of affinoid \(K\)-spaces \(\phi:Z\to X\) satisfying \(\phi(Z)\subseteq U\) the covering \((\phi^{-1}(U_{i}))_{i\in I}\) of \(Z\) admits a refinement that is a finite covering of \(Z\) by affinoid subdomains.
\end{Def}

Recall Definition \ref{Gr topology} of Grothendieck topology and our following definition of a G-topological space. We call a covering \(V=\cup_{j\in J}V_{j}\) of some admissible open set \(V\subseteq x\) by means of admissible open sets \(V_{j}\) \textbf{admissible} if for each morphism of affinoid \(K\)-spaces \(\phi:Z\to X\) satisfying \(\phi(Z)\subseteq V\) covering \((\phi^{-1}(V_{j}))_{j\in J}\) of \(Z\) admits a refinement that is finite covering of \(Z\) by affinoid subdomains. Let \(X\) be an affinoid \(K\)-space, then the admissible open sets and admissible coverings give a structure of a G-topological space. We call it the \textbf{strong Grothendieck topology}. We state without proof the following proposition.

\begin{prop}\label{topology on affinoids}
    Let \(X\) be an affinoid \(K\)-space. The strong Grothendieck topology satisfies the following properties:
\begin{enumerate}
    \item \(\varnothing\) and \(X\) are admissible open.

    \item Let \((U_{i})_{i\in I}\) be an admissible covering of an admissible open subset \(U\subseteq X\) and let \(V\subseteq U\) be a subset such that \(V\cap U_{i}\) is admissible open for all \(i\in I\), then \(V\) is admissible open in \(X\).

    \item Let \((U_{i})_{i\in I}\) be a covering of an admissible open set \(U\subseteq X\) by admissible open subsets \(U_{i}\subseteq X\) such that \((U_{i})_{i\in I}\) admits an admissible covering of \(U\) as refinement. Then \((U_{i})_{i\in I}\) itself is admissible.
    \end{enumerate}
\end{prop}

Recall that Grothendieck topology is given as a weaker condition than general topology. Indeed, Zariski topology and canonical topology that we defined previously has structures as Grothendieck topologies and we can compare strong Grothendieck topology on an affinoid \(K\)-space with them. By Proposition \ref{openness2} and our definition, we know that the strong Grothendieck topology is coarser than the canonical topology. The following proposition is stated without proof.
\begin{prop}\label{finer}
    Let \(X\) be an affinoid \(K\)-sapce. The strong Grothendieck topology on \(X\) is finer than the Zariski topology. That is, every Zariski open subset \(U\subseteq X\) is admissible open and every Zariski covering is admissible.
\end{prop}
We further define a sheaf on our G-topological space. Let \(X\) be an affinoid \(K\)-space. For any affinoid subdomain \(U\subseteq X\) we denote by \(\mathcal{O}_{X}(U)\) the affinoid \(K\)-algebra corresponding to \(U\). It is easy to check that it is a presheaf over the strong Grothendieck topology. We call it the presheaf of \textbf{affinoid functions} on \(X\). Like in canonical algebraic geometry we define the \textbf{stalk} as the direct limit $$\mathcal{O}_{X,x}=\underset{x\in U}\varinjlim\ \mathcal{O}_{X}(U).$$ We will do a sheafification of \(\mathcal{O}_{X}\) with respect to the strong Grothendieck topology in the same way as we do in classical algebraic geometry and with an abuse of notations, denote the resultant sheaf also as \(\mathcal{O}_{X}\). One can now define the analogue of ringed space in Grothendieck topology.

\begin{Def}\label{Gringed} We consider the following concepts:
\begin{enumerate}
    \item A \textbf{G-ringed} \(\mathbf{K}\)\textbf{-space} is a pair \((X,\mathcal{O}_{X})\) where \(X\) is a G-topological space and \(\mathcal{O}_{X}\) a sheaf of \(K\)-algebras on it. \((X,\mathcal{O}_{X})\) is called a \textbf{locally G-ringed} \(\mathbf{K}\)\textbf{-space} if, in addition, all stalks \(\mathcal{O}_{X,x}=\underset{x\in U}\varinjlim\ \mathcal{O}_{X}(U)\in X\), are local rings.
    
      \item A \textbf{morphism} of G-ringed
    \(K\)-spaces \((X,\mathcal{O}_{X})\to(Y,\mathcal{O}_{Y})\) is a pair \((\phi,\phi^{*})\) where \(\phi:X\to Y\) is a map, continuous with respect to Grothendieck topologies, and where \(\phi^{*}\) is a system of \(K\)-homomorphisms $$\phi^{*}_{V}:\mathcal{O}_{Y}(V)\to\mathcal{O}_{X}(\phi^{-1}(V))$$ with \(V\) varying over the admissible open subsets of \(Y\) that is compatible with restriction homomorphisms. 
    
    \item Furthermore, assuming that \((X,\mathcal{O}_{X})\) and \((Y,\mathcal{O}_{Y})\) are locally G-ringed \(K\)-spaces, a morphism $$(\phi,\phi^{*}):(X,\mathcal{O}_{X})\to(Y,\mathcal{O}_{Y})$$ is called a \textbf{morphism} of locally G-ringed \(K\)-spaces if the ring homomorphisms $$\phi^{*}_{x}:\mathcal{O}_{Y,\phi(x)}\to\mathcal{O}_{X,x}$$ for \(x\in X\) induced from \(\phi^{*}_{V}\) are local in the sense that the maximal ideal of \(\mathcal{O}_{Y,\phi(x)}\) is mapped into the maximal ideal of \(\mathcal{O}_{X,x}\).
    \end{enumerate}
\end{Def}

We can finally define our major concept.

\begin{Def}\label{RigidAnalyticSpace}
    A \textbf{rigid analytic} \(\mathbf{K}\)\textbf{-space} is a locally G-ringed \(K\)-space \((X,\mathcal{O}_{X})\) such that the G-topology of \(X\) satisfies the propositions \(1\), \(2\) and \(3\) in Proposition \ref{topology on affinoids} and such that \(X\) admits an admissible covering \((X_{i})_{i\in I}\) where \((X_{i},\mathcal{O}_{X|X_{i}})\) is an affinoid \(K\)-space for all \(i\in I\). A \textbf{morphism} of rigid \(K\)-spaces \((X,\mathcal{O}_{X})\to(Y,\mathcal{O}_{Y})\) is a morphism in the sense of locally G-ringed \(K\)-spaces.
\end{Def}
Here \(\mathcal{O}_{X|X_{i}}\) refers to the restricted sheaf. We can, indeed, glue local rigid analytic \(K\)-spaces to a global rigid analytic \(K\)-space exactly as we glue local schemes to a global scheme. For two rigid analytic \(K\)-spaces over a third one we can construct the fiber product in the same way as in classical scheme theory.
\subsection{GAGA Functor}\label{gagafunctor}
With a rigid analytic \(K\)-space \(X\), we can define concepts like \(\mathbf{\mathcal{O}_{X}}\)\textbf{-module} (including modules of \textbf{finite type}, of \textbf{finite presentation} and \textbf{coherent} modules), \textbf{closed immersion}, \textbf{quasi-compactness} and \textbf{separatedness} analogous to the way we define them in classical algebraic geometry, as long as we take affinoid \(K\)-spaces in place of affine spaces and admissible affinoid coverings in place of open coverings. To further bridging the relationship between \(K\)-schemes of locally finite type and rigid analytic \(K\)-spaces, we give the following proposition.

\begin{prop}\label{BridgingRelationship}
    Let \((Z,\mathcal{O}_{Z})\) be a \(K\)-scheme of locally finite type. Then there exists a rigid analytic \(K\)-space \((Z^{rig},\mathcal{O}_{Z^{rig}})\) together with a morphism of locally \(G\)-ringed \(K\)-spaces $$(i,i^{*}):(Z^{rig},\mathcal{O}_{Z^{rig}})\to(Z,\mathcal{O}_{Z})$$ such that for any rigid analytic \(K\)-space \((Y,\mathcal{O}_{Y})\) and a morphism of locally \(G\)-ringed \(K\)-spaces \((Y,\mathcal{O}_{Y})\to(Z,\mathcal{O}_{Z})\), the latter factors through \((i,i^{*})\) via a unique morphism of rigid analytic \(K\)-spaces \((Y,\mathcal{O}_{Y})\to(Z^{rig},\mathcal{O}_{Z^{rig}})\).
\end{prop}

We call such an \((Z^{rig},\mathcal{O}_{Z^{rig}})\) a \textbf{rigid analytification} of \((Z,\mathcal{O}_{Z})\).

To give an example and to construct the rigid analytification, we start with the affine \(n\)-space \(\mathbb{A}^{n}_{K}\). We set \(T_{n}(r)\) for \(r>0\) to be the \(K\)-algebra of all power series \(\sum_{\nu\in\mathbb{N}^{n}}a_{\nu}\zeta^{\nu}\) with \(a_{\nu}\in K\) for each \(\nu\) satisfying $$\mathrm{lim}_{||\nu||\to\infty}a_{\nu}r^{\nu}=0.$$ Therefore \(T_{n}(r)\) consists of all power series converging on a closed \(n\)-dimensional ball of radius \(r\). Choose \(c\in K\) such that \(||c||>1\) and we may identify \(T_{n}^{(i)}=T_{n}(||c||^{i})\) with the Tate algebra $$K<c^{-i}\zeta_{1},...,c^{-i}\zeta_{n}>.$$ The relation
\begin{equation}
    K[\zeta]\xhookrightarrow{}...\xhookrightarrow{}T_{n}^{(1)}\xhookrightarrow{}T_{n}^{(0)}=T_{n}\nonumber
\end{equation}
give rise to inclusions of affinoid subdomains
\begin{equation}
    \mathbb{B}^{n}=\mathrm{Sp}\ T_{n}^{(0)}\xhookrightarrow{}\mathrm{Sp}\ T_{n}^{(1)}\xhookrightarrow{}...\nonumber
\end{equation}
where \(\mathrm{Sp}\ T_{n}^{(i)}\) can be interpreted as the \(n\)-dimensional ball of radius \(||c^{i}||\). Analogous to the way that we glue schemes, we can glue all \(\mathrm{Sp}\ T_{n}^{(i)}\) together as the "union" of these balls. And it is obvious that it is independent of \(c\). The resulting rigid analytic \(K\)-space \(\mathbb{A}^{n,rig}_{K}\) is the rigid analytification of \(\mathbb{A}^{n}_{K}\).

In general, for any affine \(K\)-scheme of finite type of the form \(\mathrm{Spec}\ K[\zeta_{1},...,\zeta_{n}]/\mathfrak{a}\) with an ideal \(\mathfrak{a}\subseteq K[\zeta_{1},...,\zeta_{n}]\). We can check the maps
\begin{equation}
    K[\zeta_{1},...,\zeta_{n}]\to...\to T^{(1)}_{n}/(\mathfrak{a})\to T^{(0)}/(\mathfrak{a})\nonumber
\end{equation}
and the associated sequence of inclusions
\begin{equation}
    \mathrm{Max}\ T^{(0)}_{n}/(\mathfrak{a})\xhookrightarrow{}...\xhookrightarrow{}\mathrm{Max}\ K[\zeta_{1},...,\zeta_{n}]/(\mathfrak{a})\nonumber
\end{equation}
while \(\mathrm{Max}\) is the set of maximal ideals. We state without proof that \(\mathrm{Max}\ K[\zeta_{1},...,\zeta_{n}]/\mathfrak{a}\) equals the union of all \(\mathrm{Max}\ T^{(i)}_{n}/(\mathfrak{a})\).  The union of all \(\mathrm{Sp}\ T_{n}^{(i)}/(\mathfrak{a})\) can be glued into a rigid analytic \(K\)-space can be proved to be the rigid analytification of \(\mathrm{Spec}\ K[\zeta_{1},...,\zeta_{n}]/\mathrm{a}\). We can glue the affine analytifications together for \(K\)-schemes of finite type in general. One has the following proposition to characterize analytification.

\begin{prop}\label{analytification}
    Every \(K\)-scheme \(Z\) of locally finite type admits an analytification \(Z^{rig}\to Z\). Furthermore, the underlying map of sets identifies the points of \(Z^{rig}\) with the closed points of \(Z\).
\end{prop}

With the universal property we know that morphisms between \(K\)-schemes of locally finite type admit analytifications as well. Therefore one has the following property.

\begin{prop}\label{functor}
    Rigid analytification defines a functor from the category of \(K\)-schemes of locally finite type to the category of rigid analytic \(K\)-spaces.
\end{prop}

We call it \textbf{GAGA functor} since it relates some analytic data with algebraic data.

\begin{rem}
    On the contrary not any rigid analytic \(K\)-space has a structure of \(K\)-scheme. We will talk about this on \ref{cohomologytheories}.
\end{rem}
\subsection{Dagger Spaces}\label{daggerspaces}
Before continuing to the cohomology theories, we define an antecedent concept of \textit{dagger spaces}. We recall our notation of \(T_{n}(r)\) in Section \ref{definitionand}. We set $$W_{n}=\cup_{r>1}T_{n}(r).$$ We define the \(\mathbf{K}\)\textbf{-dagger algebra} \(A\) as a quotient of some \(W_{n}\) and its \textbf{norm} to be the Gauss norm given by the epimorphism \(W_{n}\to A\). We can prove that \(T_{n}\) is the completion of \(W_{n}\) and any completion of a \(K\)-dagger algebra is a \(K\)-algebra. The \textbf{morphism} between \(K\)-dagger algebras is a \(K\)-algebra morphism between them.

For a \(K\)-dagger space \(A'\), we set \(\mathrm{Sp}\ A'\) to be the set of maximal ideals, as in \ref{definitionand}. It is an \textbf{affinoid} \(\mathbf{K}\)\textbf{-dagger algebra}. Let \(A\) be the completion of \(A'\). The natural map \(\mathrm{Sp}\ A\to\mathrm{Sp}\ A'\) can be proved to be bijective. Like in \ref{definitionand} we can glue affinoid spaces together into a \(\mathbf{K}\)\textbf{-dagger space}. Then there is a faithful functor from the category of \(K\)-dagger spaces to the category of \(K\) rigid analytic spaces, assigning to a \(K\)-dagger space \(X'\) a \(K\)-rigid analytic space \(X\) which is called the \textbf{associated rigid analytic space}. We state without proof that there exists a natural morphism of ringed spaces \(x:X\to X'\) which induces isomorphisms between the underlying \(G\) topological spaces and which is locally compatible to the map \(\mathrm{Sp}\ A\to \mathrm{Sp}\ A'\). It can inherit a strong \(G\) topology from the aforementioned bijection. We can also generate a definition of \textbf{admissible open sets} on the \(K\)-dagger space therefore.
\subsection{Cohomology Theories}\label{cohomologytheories}
Recall that we remarked at Section \ref{gagafunctor} that the concept of an \(\mathcal{O}_{X}\)-module is defined on a rigid analytic \(K\)-space \(X\). If \(\mathcal{F}\) is an \(\mathcal{O}_{X}\)-module, we define the \textbf{\v{C}ech cohomology} groups \(H^{q}(\mathcal{U},\mathcal{F})\) for any admissible covering \(\mathcal{U}\) as in classical algebraic geometry, and we define \textbf{\v{C}ech cohomology} groups in general as
\begin{equation}
    \check{H}^{q}(X,\mathcal{F})=\varinjlim_{\mathcal{U}}\ H^{q}(\mathcal{U},\mathcal{F})\nonumber
\end{equation}
where the limit runs over all admissible coverings of \(X\).

Again let \(\Gamma(X,\cdot)\) be the global section functor on \(X\). We define the \(q^{th}\) \textbf{Grothendieck cohomology group} as
\begin{equation}
    H^{q}(X,\mathcal{F})=R^{q}\Gamma(X,\mathcal{F}),\nonumber
\end{equation}
that is the \(q^{th}\) right-derived functor of the global section functor \(\Gamma(X,\cdot)\).

We state without proof that if \(X\) is a \(K\)-dagger space or a rigid analytic \(K\)-space, there exists a way analogous to that in \ref{differentials} to define a \textbf{module of relative differential forms} of \textbf{degree} \(n\) \(\Omega^{n}_{X/K}\). The concept of \textbf{smoothness} is defined analogously. We state without proof that \(X'\) is smooth if and only if \(X\) is smooth.

Now let \(X\) be a rigid analytic space. Let \(T'\) be a \(K\)-dagger space such that the rigid analytic space \(T\) is its completion. Let \(\phi:X\to T\) be a closed immersion. Set \(\Psi(\phi,X')\) to be the set of admissible open subsets \(U'\) of \(X'\) for which \(\phi\) factors as
\begin{equation}
    X\to U \to S\nonumber
\end{equation}
such that \(U\) is associated with \(U'\) and \(U\to S\) is the embedding of rigid analytic spaces associated with the embedding of \(K\)-dagger spaces \(U'\to S'\). Suppose that \(S\) is smooth. Let \(s:S'\to S\) be the natural morphism of ringed spaces, we set that 
\begin{equation}
    R\Gamma_{dR}(X/K)=R\Gamma_{dR}(X)=R\Gamma(X,(s\circ\phi)^{-1}\Omega^{\cdot}_{S/k})\nonumber
\end{equation}
and we define the \(q^{th}\) \textbf{de Rham cohomology group} as
\begin{equation}
    H^{q}_{dR}(X)=H^{q}(R\Gamma_{dR}(S)).\nonumber
\end{equation}
The \textbf{\'etale cohomology} is defined on \(X\) as an \'etale site.

We note that those cohomology groups also have a natural structure of a \(\mathbb{C}_{K}\)-representation.
\section{Proceedings and Problems}
We mainly refer to \cite{Sch1}, \cite{Sch2}, \cite{Sch3}, \cite{Sch6}, \cite{Sch4}, \cite{Sch5} and \cite{Has} for this Section.

Recall that in Proposition \ref{analytification} we give every \(K\)-scheme of locally finite type an analytification.

\begin{ex}
 Consider the \textbf{rigid analytic Hopf surface}, defined as
\begin{equation}
    X=(\mathbb{A}^{n,rig}_{K}-\{(0,0)\})/q^{\mathbb{Z}}\nonumber
\end{equation}
such that \(q\in K^{*}\) and \(||q||<1\). Here \(q^{n}\) refers to the map $$(\zeta_{1},\zeta_{2})\to(q^{n}\zeta_{1},q^{n}\zeta_{2})$$ locally on the affinoid \(K\)-algebra. It can be proved that this is not an analytification of any \(K\)-scheme of locally finite type. Although Proposition \ref{functor} establishes rigid analytification as a functor, the existence of a rigid analytic \(K\)-space that is not a rigid analytification of any \(K\)-scheme of locally finite type means that the theory of rigid analytic \(K\)-space is a non-trivial extension of the theory of \(K\)-scheme of locally finite type.
\end{ex}
We further consider Theorem \ref{FaltingsDecompositionThm} where \(K\) is a non-archimedian field and \(X\) is a smooth \(K\)-scheme of locally finite type. We ask the question whether Equation \eqref{FaltingsDecomposition} (which is equivalent to Equation \eqref{RewritingFaltings} still holds for rigid analytic \(K\)-spaces in general. The following result is proved by Peter Scholze in \(2013\) (see \cite{Sch2}).
\begin{thm}\label{scholze decomposition}
    For any proper smooth rigid analytic \(K\)-space, there is a canonical isomorphism
    \begin{equation}\label{analyticdecomposition}
        \mathbb{C}_{K}\otimes_{\mathbb{Q}_{l}}H^{n}_{et}(X_{\overline{K}},\mathbb{Q}_{l})\cong\oplus_{p+q=n}(\mathbb{C}_{K}(-q)\otimes_{K}H^{p}(X,\Omega^{q}_{X/K}))
    \end{equation}
    which an equivalence of structure as \(\mathbb{C}_{K}\) representatives. Moreover, we have
    \begin{equation}\label{dimensionsame}
        \mathrm{dim}\ H^{n}_{et}(X_{\overline{K}},\mathbb{Q}_{l})=\mathrm{dim}\ (B_{dR}\otimes_{\mathbb{Q}_{l}}V).
    \end{equation}
    And for \(i\geq 0\), we have
    \begin{equation}\label{KKQl}
        \sum_{j=0}^{i}\mathrm{dim}_{K} H^{i-j}(X,\Omega^{j}_{X})=\mathrm{dim}_{K}H^{i}_{dR}(X)=\mathrm{dim}_{\mathbb{Q}_{l}}H^{i}_{et}(X_{\overline{K}},\mathbb{Q}_{l}).
    \end{equation}
\end{thm}
\begin{rem}
    Note that on the right-hand-side of Equation \ref{analyticdecomposition}, we use Grothendieck cohomology instead of de Rham cohomology and de Rham cohomology begins to appear in Equation \ref{KKQl}. Unlike the complex case, moreover, \textit{no} K\"alheran condition is necessary. 
\end{rem}
From Corollary \ref{hodge numbers}, we conclude the analogous condition of that \(b_{r}=\sum_{p+q=r}h^{p,q}\), but we do not yet have the conclusion known as \textbf{Hodge Symmetry}, that is, \(h^{p,q}=h^{q,p}\). A partial progression on this issue is made by David Hansen and Shizhang Li in \cite{Has}. Before stating the result, several preparatory concepts are necessary.
\begin{Def}\label{FormalSpectrum}
    Let \(A\) be a complete discrete valuation ring with uniformiser \(\pi\), the \textbf{formal spectrum} \(\mathrm{Spf}(A)\) is defined as 
    \begin{equation}
        \mathrm{Spf}(A)=\varinjlim \mathrm{Spec}(A/(\pi^{n}))\nonumber
    \end{equation}
    while the inverse system is defined by the canonical immersion from \(\mathrm{Spec}(R/(\pi^{n}))\) to \(\mathrm{Spec}(R/(\pi^{n+1}))\). There is a topology on \(\mathrm{Spf}(A)\) generated naturally by the immersions.
\end{Def}
Recall Definition \ref{model} of models and we need to define the rigid analytic equivalence of models. We state without explicit definition that there is a well-defined concept of a \textbf{formal model} \(\mathfrak{X}\) of a rigid analytic space \(X\) over \(\mathrm{Spf}(\mathcal{O}_{K})\), that is, a flat scheme of finite type over \(\mathrm{Spf}(\mathcal{O}_{K})\) that has a \textbf{rigid analytic generic fiber} that is isomorphic to \(X\), something analogous to Definition \ref{model}. Also, we define the \textbf{specific fiber} of \(\mathfrak{X}\) as the scheme \(\mathfrak{X}_{P}\times_{\mathrm{Spf}(\mathcal{O}_{K})}\mathrm{Spec}(P)\) while \(P\) is the residue field of \(K\). We state the result by Hansen and Li (see \cite{Has}).
\begin{thm}\label{hodge symmetry}
    Let \(A\) be a complete discrete valuation ring and let \(X\) be a smooth proper rigid analytic space over \(K\). Assume that \(X\) has a formal model \(\mathfrak{X}\) over \(\mathrm{Spf}(\mathcal{O}_{K})\) whose speical fiber is projective. Then
    \begin{equation}\label{symmetry dimension}
        \mathrm{dim}_{K}H^{1}(X,\Omega^{0}_{X})=\mathrm{dim}_{K}H^{0}(X,\Omega^{1}_{X}).
    \end{equation}
\end{thm}
Note here the existence of a formal model is analogous to the K\"ahler condition in complex geometry.
\begin{rem}
    Given the conditions of Theorem \ref{hodge symmetry}  sometimes the higher Hodge symmetry still fails. Indeed, Alexander Petrov proved in \cite{Pet}  that for any pair of positive integers \(i\neq j\) with \(i+j>3\), there is a smooth proper rigid analtyic space over \(\mathbb{Q}_{p}\) admitting a smooth formal model \(\mathfrak{X}\) over \(\mathrm{Spf}\ \mathbb{Z}_{p}\) whose special fiber is projective that fails the Hodge symmetry.
\end{rem}
The question arises naturally: is it possible that, under a stricter condition, some Hodge symmetry of higher cohomology will hold?
\section{Perfectoid Space and Condensed Mathematics}\label{Perfectoids}
We mainly refer to \cite{Sch1}, \cite{Sch6}, \cite{Sch4} and \cite{Cla} for this Session.

Note that our previous works are based on \'etale cites. However, the proof of Theorem \ref{scholze decomposition}
involves lots of difficulties. In particular, the local structures of rigid analytic spaces are much more complicated under \'etale topology. In order to simplify the problem, Peter Scholze introduced the tool of \textit{pro-\'etale topology}, where inverse limits of \'etale morphism under certain conditions, instead of only \'etale morphisms, are allowed in constructing the site. Furthermore, Scholze showed that under pro-\'etale topology, a rigid analytic \(K\)-space is locally a \textit{perfectoid space}. The following are some definitions necessary for understanding this concept.
\begin{Def}\label{bounded}
    Let \(A\) be a topological ring. A subset \(T\subseteq A\) is called \textbf{bounded}, if for every neighbourhood \(U\) of zero, there exists an open neighbouhood \(V\) of zero such that 
    $$T\cdot V=\{t\cdot v\ |\ t\in T,v\in V\}\subseteq U.$$ An element \(a\in A\) is called \textbf{power-bounded} if the set \(\{a^{n}\ |\ n\in N\}\) is bounded.
\end{Def}
It is easy to see that for a discrete valuation ring, a power-boudned set forms a subring. Scholze defined perfectoid field (see \cite{Sch1}).
\begin{Def}\label{perfectoid field}
    Let \(K\) be a fractional field of some complete discrete valuation ring \(A\) and assume that \(p\neq 0\) is the characteristic of its residue field, \(K\) is a \textbf{perfectoid field} if the Frobenius endomorphism \(\Phi\) is surjective on \(K^{\circ}/p\) where \(K^{\circ}\) denotes the ring of power-bounded elements. 
\end{Def}
We define \textit{perfectoid algebras} and \textit{perfectoid spaces} over a perfect field in ways analogous to the way we define algebras and schemes in classical algebraic geometry. 

For any perfectoid field \(K\), there is a \textit{tilt} \(K^{\flat}\). The underlying set given a structure of a perfectoid field with characteristic \(p\),
$$\underset{x\to x^{p}}{\varprojlim} K.$$ If \(X\) is a perfecoid space over a perfectoid field \(K\), then we can form a perfectoid space \(X^{\flat}\) over \(K^{\flat}\). There is a theorem called \textit{tilting equivalence}, saying that the tilting functor \(X\to X^{\flat}\) induces an equivalence of categories. Note that the tilting equivalence associates the geometry over fields of mixed characteristics with the geometry over a field of finite characteristics. This largely simplifies the former.

The studies in perfectoid space lead to Dustin Clausen and Peter Scholze's development of \textit{condensed mathematics}, a framework that unites topology, complex analysis and functional analysis. The basic idea is: the pro-\'etale site is developed to domesticate some wild topologies in the rigid analytic space and the technique to domesticate misbehaved topologies is generalized into condensed mathematics. Instead of study topology spaces per se, we study \textit{condensed} objects.

\begin{Def}\label{profinite set}
    A \textbf{profinite set} is a compact, Hausdorff and totally disconnected topological space. The \textbf{site} of profinite sets is the set of all profinite sets with the Grothendieck topology given by finite, jointly surjective collections of maps.
\end{Def}
\begin{ex}\label{profinite set example}
    The set \(\mathbb{N}\cup\{+\infty\}\) under the topology inherited from \(\mathbb{R}\cup\{+\infty\}\) is a profinite set.
\end{ex}
We will now replace traditional topology on a space \(X\) with "continuous map from profinite sets to \(X\)". The following definition is worked out by Dustin Clausen and Peter Scholze (see \cite{Cla} and \cite{Sch6}).
\begin{Def}\label{condensed set}
    A \textbf{condensed set} (respectively, \textbf{condensed group, ring} etc) is a sheaf of sets (respectively, groups, rings) on the site of profinite sets.
\end{Def}
To any topological space (respectively, group or ring) \(X\) one can associate a condensed set (respectively, group or ring). Especially, the continuous maps from the set in Example \ref{profinite set example}  to a topological space \(X\) form a sheaf of group that codes the information that a series \(x_{n}\to x_{\infty}\).

There are notions of \textit{solid} and \textit{liquid} modules, and the former has more to do with the non-archimedian, algebraic side of the theory (since non-archimedian spaces are totally disconnected) and the latter has more relationship with the analytic, complex side of the theory (intuitively, analytic objects are more "fluid"). But the details of the theory is beyond the scope of this note.

\footnote{Email [1]: dolivia@unc.edu and Email [2]: ykx13579@gmail.com.}


\begin{bibdiv}
    \begin{biblist}
\bib{aty}{article}{
    AUTHOR = {Atiyah, Michael},
    AUTHOR = {Macdonald, Ian}
     TITLE = {Introduction to Commutative Algebra},
   JOURNAL = {Addison-Wesley Publishing Company},
    VOLUME = {},
      YEAR = {1969},
    NUMBER = {},
     PAGES = {},
      ISSN = {},
       DOI = {},
       URL = {},
}
\bib{berth}{article}{
    AUTHOR = {Berthelot, Pierre},
    AUTHOR = {Ogus, Arthur}
     TITLE = {Notes on Crystalline Cohomology},
   JOURNAL = {Princeton University Press and University of Tokyo Press},
    VOLUME = {},
      YEAR = {1978},
    NUMBER = {},
     PAGES = {},
      ISSN = {},
       DOI = {},
       URL = {},
}

    \bib{Bo1}{article}{
    AUTHOR = {Bosch, Siegfried},
     TITLE = {Algebraic Geometry and Commutative Algebra},
   JOURNAL = {Springer-Verlag London},
    VOLUME = {},
      YEAR = {2013},
    NUMBER = {},
     PAGES = {},
      ISSN = {},
       DOI = {},
       URL = {},
}


\bib{Bo2}{article}{
    AUTHOR = {Bosch, Siegfried},
     TITLE = {Lectures on Formal and Rigid Geometry},
   JOURNAL = {Springer International Publishing Switzerland},
    VOLUME = {},
      YEAR = {2014},
    NUMBER = {},
     PAGES = {},
      ISSN = {},
       DOI = {},
       URL = {},
}
\bib{bre}{article}{
    AUTHOR = {Breuil, Christophe},
    AUTHOR = {Conrad, Brian},
    AUTHOR = {Diamond, Fred},
    AUTHOR = {Taylor, Richard},
     TITLE = {On the Modularity of Elliptic Curves Over \(\mathbb{Q}\)},
   JOURNAL = {Journal of the American Mathematical Society},
    VOLUME = {14},
      YEAR = {2001},
    NUMBER = {},
     PAGES = {843-939},
      ISSN = {},
       DOI = {},
       URL = {},
}

\bib{bri}{article}{
    AUTHOR = {Brinon, Olivier},
    AUTHOR = {Conrad, Brian},
     TITLE = {CMI Summer Scholl Notes on \(p\)-adic Hodge Theory (Preliminary Version)},
   JOURNAL = {},
    VOLUME = {},
      YEAR = {2009},
    NUMBER = {},
     PAGES = {},
      ISSN = {},
       DOI = {},
       URL = {},
}
\bib{chern}{article}{
    AUTHOR = {Chern, S.S.},
    AUTHOR = {Chen, W.H.},
    AUTHOR = {Lam, K.S.},
     TITLE = {Lectures on Differential Geometry},
   JOURNAL = {World Scientific Publishing Co. Pte. Ltd.},
    VOLUME = {},
      YEAR = {2000},
    NUMBER = {},
     PAGES = {},
      ISSN = {},
       DOI = {},
       URL = {},
}
\bib{Cla}{article}{
    AUTHOR = {Clausen, Dustin},
    AUTHOR = {Scholze, Peter},
     TITLE = {Condensed Mathematics and Complex Geometry},
   JOURNAL = {},
    VOLUME = {},
      YEAR = {2026},
    NUMBER = {},
     PAGES = {},
      ISSN = {},
       DOI = {},
       URL = {},
}
\bib{Cox}{article}{
    AUTHOR = {Cox, David},
    AUTHOR = {Katz, Sheldon},
     TITLE = {Mirror Symmetry and Algebraic Geometry},
   JOURNAL = {American Mathematical Society},
    VOLUME = {},
      YEAR = {1999},
    NUMBER = {},
     PAGES = {},
      ISSN = {},
       DOI = {},
       URL = {},
}


\bib{Fal1}{article}{
    AUTHOR = {Faltings, Gerd},
     TITLE = {Crystalline Cohomology and \(p\)-adic Galois Representations},
   JOURNAL = {Algebraic Analysis, Geometry and Number Theory},
    VOLUME = {},
      YEAR = {1998},
    NUMBER = {1},
     PAGES = {25-80},
      ISSN = {},
       DOI = {},
       URL = {},
}
\bib{Fal2}{article}{
    AUTHOR = {Faltings, Gerd},
     TITLE = {\(p\)-adic Hodge Theory},
   JOURNAL = {Journal of the American Mathematical Society},
    VOLUME = {1},
      YEAR = {1998},
    NUMBER = {1},
     PAGES = {255-299},
      ISSN = {},
       DOI = {},
       URL = {},
}

\bib{fon1}{article}{
    AUTHOR = {Fontaine, Jean-Marc},
     TITLE = {Sur certains types de représentations \(p\)-adiques du groupe de Galois d’un corps local},
   JOURNAL = {Inventiones Mathematicae},
    VOLUME = {81},
      YEAR = {1985},
    NUMBER = {},
     PAGES = {515-538},
      ISSN = {},
       DOI = {},
       URL = {},
}
\bib{fon2}{article}{
    AUTHOR = {Fontaine, Jean-Marc},
     TITLE = {Le corps des périodes \(p\)-adiques},
   JOURNAL = {Astérisque},
    VOLUME = {223},
      YEAR = {1994},
    NUMBER = {},
     PAGES = {59-111},
      ISSN = {},
       DOI = {},
       URL = {},
}
\bib{gro}{article}{
    AUTHOR = {Große-Klönne, Elmar},
     TITLE = {De Rham Cohomology of Rigid Spaces},
   JOURNAL = {Math. Zeit.},
    VOLUME = {247},
      YEAR = {2004},
    NUMBER = {},
     PAGES = {223-240},
      ISSN = {},
       DOI = {},
       URL = {},
}
\bib{GrHa}{article}{
    AUTHOR = {Griffith, Phillip},
    AUTHOR = {Harris, Joseph},
     TITLE = {Priciples of Algebraic Geometry},
   JOURNAL = {John Wiley and Sons, Inc},
    VOLUME = {},
      YEAR = {1994},
    NUMBER = {5},
     PAGES = {},
      ISSN = {},
       DOI = {},
       URL = {},
}
\bib{FrSh}{article}{
    AUTHOR = {Freud, Diamond},
    AUTHOR = {Shurman, Jerry},
     TITLE = {A First Course in Modular Forms},
   JOURNAL = {Springer Science and Business Media, Inc},
    VOLUME = {},
      YEAR = {2005},
    NUMBER = {},
     PAGES = {},
      ISSN = {},
       DOI = {},
       URL = {},
}
\bib{Has}{article}{
    AUTHOR = {Hansen, David},
    AUTHOR = {Li, Shizhang}
     TITLE = {Line Bundles on Rigid Varieties and Hodge Symmetry},
   JOURNAL = {Math. Zeit.},
    VOLUME = {296},
      YEAR = {2017},
    NUMBER = {},
     PAGES = {1777-1786},
      ISSN = {},
       DOI = {},
       URL = {},
}

\bib{Hat}{article}{
    AUTHOR = {Hatcher, Allen},
     TITLE = {Algebraic Topology},
   JOURNAL = {Cambridge University Press},
    VOLUME = {},
      YEAR = {2001},
    NUMBER = {},
     PAGES = {},
      ISSN = {},
       DOI = {},
       URL = {},
}
\bib{Kont}{article}{
    AUTHOR = {Kontsevich, Maxim},
     TITLE = {Homological Algebra of Mirror Symmetry},
   JOURNAL = {Proceedings of the International Congress of Mathematicians},
    VOLUME = {},
      YEAR = {1994},
    NUMBER = {},
     PAGES = {120-139},
      ISSN = {},
       DOI = {},
       URL = {},
}
\bib{Lee}{article}{
    AUTHOR = {Lee, John M.},
     TITLE = {Introduction to Smooth Manifolds},
   JOURNAL = {Springer Science+Business Media New York},
    VOLUME = {},
      YEAR = {2003},
    NUMBER = {},
     PAGES = {},
      ISSN = {},
       DOI = {},
       URL = {},
}

\bib{Liu}{article}{
    AUTHOR = {Liu, Qing},
     TITLE = {Algebraic Geometry and Arithmetic Curves},
   JOURNAL = {Oxford University Press},
    VOLUME = {},
      YEAR = {2002},
    NUMBER = {},
     PAGES = {},
      ISSN = {},
       DOI = {},
       URL = {},
}

\bib{Mac}{article}{
    AUTHOR = {Mac Lane, Saunders},
     TITLE = {Categories for the Working Mathematician},
   JOURNAL = {Springer-Verlag New York Inc.},
    VOLUME = {},
      YEAR = {1971},
    NUMBER = {},
     PAGES = {},
      ISSN = {},
       DOI = {},
       URL = {},
}
\bib{mil}{article}{
    AUTHOR = {Milnor, John},
    AUTHOR = {Stasheff, James},
     TITLE = {Characterisitc Classes},
   JOURNAL = {Princeton University Press and University of Tokyo Press},
    VOLUME = {},
      YEAR = {1974},
    NUMBER = {},
     PAGES = {},
      ISSN = {},
       DOI = {},
       URL = {},
}
\bib{Neu}{article}{
    AUTHOR = {Neukirch, J\"urgen},
     TITLE = {Algebraic Number Theory},
   JOURNAL = {Springer-Verlag Berlin Heidelberg},
    VOLUME = {},
      YEAR = {1999},
    NUMBER = {},
     PAGES = {},
      ISSN = {},
       DOI = {},
       URL = {},
}
\bib{Pet}{article}{
    AUTHOR = {Petrov, Alexander},
     TITLE = {Rigid-Analytic Varieties with Projective Reduction Violating Hodge Symmetry},
   JOURNAL = {Compositio Mathematica},
    VOLUME = {157},
      YEAR = {2021},
    NUMBER = {},
     PAGES = {625-640},
      ISSN = {},
       DOI = {},
       URL = {https://math.mit.edu/~poonen/782/782notes.pdf},
}
\bib{Po}{article}{
    AUTHOR = {Poonen, Bjorn},
     TITLE = {Introduction to Arithmetic Geometry},
   JOURNAL = {Journal of the American Mathematical Society},
    VOLUME = {},
      YEAR = {},
    NUMBER = {},
     PAGES = {},
      ISSN = {},
       DOI = {},
       URL = {https://math.mit.edu/~poonen/782/782notes.pdf},
}

\bib{rom}{article}{
    AUTHOR = {Roman, Steven},
     TITLE = {Field Theory},
   JOURNAL = {Springer-Verlag},
    VOLUME = {},
      YEAR = {1995},
    NUMBER = {},
     PAGES = {},
      ISSN = {},
       DOI = {},
       URL = {},
}
\bib{Sch1}{article}{
    AUTHOR = {Scholze, Peter},
     TITLE = {Perfectoid spaces: a survey},
   JOURNAL = {Current Developments in Mathematics},
    VOLUME = {},
      YEAR = {2012},
    NUMBER = {},
     PAGES = {},
      ISSN = {},
       DOI = {},
       URL = {},
}

\bib{Sch2}{article}{
    AUTHOR = {Scholze, Peter},
     TITLE = {\(p\)-adic Hodge theory for rigid-analytic spaces},
   JOURNAL = {Forum of Mathematics, Pi},
    VOLUME = {1},
      YEAR = {2013},
    NUMBER = {},
     PAGES = {},
      ISSN = {},
       DOI = {},
       URL = {},
}

\bib{Sch3}{article}{
    AUTHOR = {Scholze, Peter},
     TITLE = {\(p\)-adic Geometry},
   JOURNAL = {Proceedings of the ICM},
    VOLUME = {1},
      YEAR = {2018},
    NUMBER = {},
     PAGES = {},
      ISSN = {},
       DOI = {},
       URL = {},
}

\bib{Sch6}{article}{
    AUTHOR = {Scholze, Peter},
     TITLE = {Lectures on Condensed Mathematics},
   JOURNAL = {},
    VOLUME = {},
      YEAR = {2019},
    NUMBER = {},
     PAGES = {},
      ISSN = {},
       DOI = {},
       URL = {},
}
\bib{Sch4}{article}{
    AUTHOR = {Scholze, Peter},
    AUTHOR={Weinstein, Jared}
     TITLE = {Berkeley Lectures on \(p\)-adic Geometry},
   JOURNAL = {Princeton University Press},
    VOLUME = {1},
      YEAR = {2020},
    NUMBER = {},
     PAGES = {},
      ISSN = {},
       DOI = {},
       URL = {},
}
\bib{Sch5}{article}{
    AUTHOR = {Scholze, Peter},
    AUTHOR={Bhatt, Bhargav},
    AUTHOR={Morrow, Matthew}
     TITLE = {Integral \(p\)-aidc Hodge Theory},
   JOURNAL = {Publ. Math. de l'IHÉS},
    VOLUME = {128},
      YEAR = {2018},
    NUMBER = {},
     PAGES = {},
      ISSN = {},
       DOI = {},
       URL = {},
}
\bib{ser}{article}{
    AUTHOR = {Serre, Jean-Pierre},
     TITLE = {Local Fields},
   JOURNAL = {Springer Science+Business Media New York},
    VOLUME = {},
      YEAR = {1979},
    NUMBER = {},
     PAGES = {},
      ISSN = {},
       DOI = {},
       URL = {},
}
\bib{sil}{article}{
    AUTHOR = {Silverman, Joseph},
     TITLE = {The Arithmetic of Elliptic Curves},
   JOURNAL = {Springer Dordrecht Heidelberg London New York},
    VOLUME = {},
      YEAR = {1986},
    NUMBER = {},
     PAGES = {},
      ISSN = {},
       DOI = {},
       URL = {},
}

\bib{Ta}{article}{
    AUTHOR = {Tate, John},
     TITLE = {Rigid Analytic Spaces},
   JOURNAL = {Inventiones Math.},
    VOLUME = {12},
      YEAR = {1971},
    NUMBER = {},
     PAGES = {257-289},
      ISSN = {},
       DOI = {},
       URL = {},
}


\bib{Vos1}{article}{
    AUTHOR = {Vosin, Claire},
     TITLE = {Hodge Theory and Complex Algebraic Geometry I},
   JOURNAL = {Cambridge University Press (translated by Leila Schneps)},
    VOLUME = {},
      YEAR = {2002},
    NUMBER = {},
     PAGES = {},
      ISSN = {},
       DOI = {},
       URL = {},
}


\bib{Vos2}{article}{
    AUTHOR = {Vosin, Claire},
     TITLE = {Hodge Theory and Complex Algebraic Geometry II},
   JOURNAL = {Cambridge University Press (translated by Leila Schneps)},
    VOLUME = {},
      YEAR = {2003},
    NUMBER = {},
     PAGES = {},
      ISSN = {},
       DOI = {},
       URL = {},
}
\bib{wei}{article}{
    AUTHOR = {Weibel, Charles},
     TITLE = {An Introduction to Homological Algebra},
   JOURNAL = {Cambridge University Press},
    VOLUME = {},
      YEAR = {1997},
    NUMBER = {},
     PAGES = {},
      ISSN = {},
       DOI = {},
       URL = {},
}
\bib{wil}{article}{
    AUTHOR = {Wiles, Andrew},
     TITLE = {Modular Elliptic Curves and Fermat's Last Theorem},
   JOURNAL = {Annals of Mathematics},
    VOLUME = {141},
      YEAR = {1995},
    NUMBER = {},
     PAGES = {443-551},
      ISSN = {},
       DOI = {},
       URL = {},
}
\bib{xar}{article}{
    AUTHOR = {Xarles, Xavier},
     TITLE = {Comparison Theorems between Crystalline and \'Etale Cohomology: A Short Introduction (Preliminary Version)},
   JOURNAL = {},
    VOLUME = {},
      YEAR = {},
    NUMBER = {},
     PAGES = {},
      ISSN = {},
       DOI = {},
       URL = {},
}

\end{biblist}
\end{bibdiv}
\end{document}